\documentclass[a4paper,11pt,twoside]{article}
\usepackage{amsfonts}
\usepackage{mathrsfs}
\usepackage{amsmath,amssymb,amsthm}
\usepackage{latexsym}
\usepackage{indentfirst}
\usepackage[top=1in, bottom=1in, left=1in, right=1in]{geometry}
\usepackage{color}
\usepackage{footmisc}
\usepackage{endnotes}
\usepackage{float}
\usepackage{authblk}
\usepackage{longtable}
\usepackage{graphicx}
\usepackage{epstopdf}
\usepackage[round,authoryear]{natbib}

\usepackage{url}
\usepackage{hyperref}

 \DeclareMathOperator*{\var}{Var}
\DeclareMathOperator*{\Ep}{{\rm E}}
\DeclareMathOperator*{\dimm}{dim}
\DeclareMathOperator*{\diag}{diag}
\DeclareMathOperator*{\vecc}{vec}
\DeclareMathOperator*{\dett}{det}
\DeclareMathOperator*{\BF}{BF_q}
\DeclareMathOperator*{\PO}{PO_q}
\DeclareMathOperator*{\poh}{PO_{\hat \theta}}
\DeclareMathOperator*{\bfh}{BF_{\hat \theta}}

\DeclareMathOperator*{\tr}{tr}

\DeclareMathOperator*{\add}{add}

\DeclareMathOperator*{\full}{full}
\def\om{\boldsymbol{\Omega}}
\def\si{\boldsymbol{\Sigma}}
\def\bxi{\boldsymbol{\Xi}}
\newcommand{\ud}{\mathrm{d}}

\newtheorem{theorem}{Theorem}
\newtheorem{lemma}{Lemma}

\theoremstyle{definition}
\newtheorem{remark}{Remark}
\newtheorem{example}{Example}
\allowdisplaybreaks

\def\A{\boldsymbol{A}}
\def\B{\boldsymbol{B}}
\def\C{\boldsymbol{C}}
\def\D{\boldsymbol{D}}
\def\E{\boldsymbol{E}}
\def\F{\boldsymbol{F}}
\def\G{\boldsymbol{G}}
\def\H{\boldsymbol{H}}
\def\I{\boldsymbol{I}}
\def\J{\boldsymbol{J}}
\def\K{\boldsymbol{K}}
\def\P{\boldsymbol{P}}
\def\R{\boldsymbol{R}}
\def\S{\boldsymbol{S}}

\def\V{\boldsymbol{V}}
\def\X{\boldsymbol{X}}

\def\e{\epsilon_n}
\def\M{\mathcal{M}}

\begin{document}

%

\title{\textbf{On Oracle Property and Asymptotic Validity of Bayesian Generalized Method of Moments}\footnote{Published version in Journal of Multivariate Analysis at \href{http://dx.doi.org/10.1016/j.jmva.2015.12.009}{http://dx.doi.org/10.1016/j.jmva.2015.12.009}
\textcopyright $<$2016$>$. This manuscript version is made available under the CC-BY-NC-ND 4.0 license \newline \href{http://creativecommons.org/licenses/by-nc-nd/4.0/}{http://creativecommons.org/licenses/by-nc-nd/4.0/} }}

\author[1]{Cheng Li \thanks{cl332@stat.duke.edu}}
\author[2]{Wenxin Jiang \thanks{wjiang@northwestern.edu}}

\affil[1]{Department of Statistical Science, Duke University}
\affil[2]{Department of Statistics, Northwestern University}

\date{}
\maketitle

\begin{abstract}
Statistical inference based on moment conditions and estimating equations is of substantial interest when it is difficult to specify a full probabilistic model. We propose a Bayesian flavored model selection framework based on (quasi-)posterior probabilities from the Bayesian Generalized Method of Moments (BGMM), which allows us to incorporate two important advantages of a Bayesian approach: the expressiveness of posterior distributions and the convenient computational method of Markov Chain Monte Carlo (MCMC). Theoretically we show that BGMM can achieve the posterior consistency for selecting the unknown true model, and that it possesses a Bayesian version of the oracle property, i.e. the posterior distribution for the parameter of interest is asymptotically normal and is as informative as if the true model were known. In addition, we show that the proposed quasi-posterior is valid to be interpreted as an approximate posterior distribution given a data summary. Our applications include modeling of correlated data, quantile regression, and graphical models based on partial correlations. We demonstrate the implementation of the BGMM model selection through numerical examples.
\end{abstract}


\fontsize{10.95}{14pt plus.8pt minus .6pt}\selectfont

\section{Introduction}
We consider the estimation problem based on the following unconditional moment restrictions
\begin{equation}\label{uncond}
\Ep \left\{g(D,\theta)\right\}=0
\end{equation}
where $D$ is a set of random variables with domain $\mathcal {D}$, $\theta$ is a $p$-dimensional vector of parameters to be estimated, and $g$ is a $m$-dimensional mapping from $\mathcal {D}\times \mathbb{R}^p$ to $\mathbb{R}^m$. Typically it is necessary to have $m\geq p$ for the point identification of $\theta$. Given an i.i.d. or stationary realization $\D=\{D_1,\ldots,D_n\}$ of $D$, one can estimate $\theta$ directly from such a set of $m$ moment functions, without needing to fully specify the underlying data generating process of $D$. In this paper, we consider the case where in \eqref{uncond}, the true parameter $\theta_0$ could possibly lie in a lower dimensional subspace. Our goal is to consistently select the relevant variables and estimate their effects, namely the nonzero components of $\theta_0$, when the specification of full probabilistic model is unavailable but a sufficient number of moment conditions are present.\\

We consider a Bayesian-flavored approach, where a quasi-posterior can be derived from a prior distribution and a quadratic form of moment restrictions. This enables us to accommodate two important advantages of the Bayesian approach: the expressiveness of the posterior distributions and the convenient computational method of MCMC. These are particularly useful for the model selection problem that we study.  We are able to report the most probable model, the second most probable model and so on, together with their quasi-posterior probabilities, which are shown to be asymptotically valid in large samples. We can also use the reversible jump MCMC algorithm (\citealt{green95}, \citealt{della02}) to traverse the space of different models and simulate the quasi-posterior probabilities.\\

For this framework of moment-based Bayesian method of model selection and model averaging, our paper will prove several appealing fundamental theorems. They will address model selection consistency, oracle property, and valid interpretation of the quasi-posterior distribution. In the following, we will first review the related works and then describe in detail the contributions of our current paper.

\subsection{GMM and BGMM}
The moment based estimation problem (\ref{uncond}) is important and has been extensively studied in econometrics and statistics. Well known methods include the generalized method of moments (GMM, \citealt{hansen82}, \citealt{hansen96}, \citealt{newey04}), the empirical likelihood (EL, \citealt{owen88}, \citealt{qin94}), the exponential tilting (ET, \citealt{kita97}), the exponential tilted empirical likelihood (ETEL, \citealt{schen05, schen07}) and the generalized empirical likelihood (GEL, \citealt{neweysmith04}). Essentially they all share the same first order efficiency of optimally weighted GMM estimator, and have been applied to independent data, time series data and panel data in econometrics. On the other hand, researchers in statistics also use the moment based methods for constructing efficient estimators, especially for clustered and correlated longitudinal data. For example, \citet{qu00} proposed a GMM type estimator to avoid the inefficiency from misspecified working correlation matrices in generalized estimating equations (GEE) for longitudinal data. \citet{wang10} considered the EL approach to address the within-subject correlation structure. Recently frequentist penalization methods have been proposed to accommodate increasing dimension $p$. See for example \citet{wang12}, \citet{leng12}, \citet{cho13}, \citet{caner13}, etc. In general, the moment based estimation methods only require information on the low order moments of $D$ and are therefore more flexible, efficient and robust to model misspecification, as long as the moment conditions are correctly specified.\\

Our work focuses on the Bayesian inference of $\theta$ under the moment constraint \eqref{uncond}. Compared to the abundance of frequentist literature, the development of Bayesian methods on this problem still remains limited. One difficulty that hinders the fully probabilistic Bayesian modeling is that some prior distribution on both the distribution of $D$ (denoted as $P_D$) and the parameter $\theta$ needs to be specified, such that the pair $(P_D,\theta)$ satisfies the set of restrictions \eqref{uncond}. Recent progress in this direction includes \citet{kita11} and \citet{florens12}. \citet{kita11} tried to minimize the Kullback-Leibler divergence of $P_D$ to a Dirichlet process, which leads to an ET type likelihood function that computationally requires optimizations within each MCMC iteration step. \citet{florens12} exploited the Gaussian process prior and required a functional transformation of the data that is only asymptotically Gaussian, which still leads to a misspecified likelihood function in finite samples. Besides, both methods have only been tested on simple examples that involve a few parameters and moments.
Instead, another analytically simpler Bayesian way of modeling \eqref{uncond} is {\it the Bayesian generalized method of moments} (BGMM), first proposed and studied by \citet{kim02} and \citet{ch03}, which constructs the simple quasi-likelihood function
\begin{equation}\label{qlik}
 q(\D|\theta)=\frac{1}{\dett\left(2\pi \V_n /n\right)^{\frac{1}{2}}}\exp\left\{ -\frac{n}{2} \bar
g(\D,\theta)^\top  \V_n^{-1} \bar g(\D,\theta)\right\},
\end{equation}
where $\bar g(\D,\theta)$ is the sample average of $g(D_i,\theta),i=1,\ldots,n$, $\V_n$ is a $m\times m$ positive definite matrix that could possibly depend on the data $\D$, and $\dett(\A)$ denotes the determinant of a matrix $\A$. Hereafter we use the symbol ``$q$" to denote the quasi-likelihood function and the quasi-posterior. This quasi-likelihood function has been studied under a Bayesian framework in \citet{kim02} and is named \textit{the limited information likelihood (LIL)}, which minimizes the Kullback-Leibler divergence of the true data generating process $P_D$ to the set of all distributions satisfying the less restrictive asymptotic constraint $\lim_{n\to \infty} \Ep \left\{n \bar g(\D,\theta_0)^\top  \V_n^{-1} \bar g(\D,\theta_0)\right\}/m=1$. This relation holds when we choose $\V_n$ to be a consistent estimator of the covariance matrix $\var(g(D,\theta_0))$. Given a prior distribution $\pi(\theta)$, the quasi-posterior takes the form
\begin{equation}\label{qpost}
 q(\theta|\D)\propto \frac{1}{\dett\left(2\pi \V_n /n\right)^{\frac{1}{2}}}\exp\left\{ -\frac{n}{2} \bar
g(\D,\theta)^\top  \V_n^{-1} \bar g(\D,\theta)\right\}\pi(\theta).
\end{equation}
By using $q(\D|\theta)$ in the Bayesian model, we only need to specify a prior on $\theta$ and thus circumvent the difficulty of directly assigning a prior on the pair $(P_D,\theta)$ with constraints \eqref{uncond}. In the computational aspect, $q(\D|\theta)$ takes an explicit analytical form that allows straightforward MCMC updating for the corresponding Bayesian posterior without any iterative optimization steps (\citealt{ch03}). Furthermore, when $\V_n$ is chosen as a consistent estimator of $\var(g(D,\theta_0))$, the exponential part of $q(\D|\theta)$ resembles the optimally weighted GMM criterion function (\citealt{hansen82}), which in large samples can be viewed as a second order approximation to the true negative log-likelihood function that follows a chi-square distribution with $p$ degrees of freedom if $m=p$ and both are fixed (\citealt{yin09}). \\

The theoretical properties of BGMM have been investigated extensively in \citet{ch03} and \citet{bc09}, who show that a Bernstein-von Mises theorem holds, i.e. the posterior distribution converges asymptotically to normal. The computational aspects of BGMM with no model selection have been investigated in \citet{yin09} and \citet{yin11}. \citet{kim14} has established the pairwise consistency theoretically when each candidate model is compared to the true model separately, and has used MCMC in simulations for such model comparison. \citet{hong12} has discussed a more general Bayesian model selection framework including BGMM as well as Bayesian GEL, and has studied the consistency of Bayes factors and Bayesian information criterion (BIC) under both nested and nonnested scenarios (see a more detailed comparison later in Remark 3 in Section 2). Other applications of BGMM include the moment inequality models (\citealt{liao10}) and the nonparametric instrumental regression (\citealt{liao11}, \citealt{kato12}). However, theoretical properties of BGMM, such as the limiting distribution and the posterior interpretation, have not been systematically studied in the context of model selection with increasing dimensionality.

\subsection{Contributions of current paper}

We study theoretical properties of BGMM in the context of model selection. The detailed contributions of the current paper include the following: \\

\noindent 1. We prove that BGMM automatically achieves the ``global model selection consistency" (see, e.g., \citealt{johnson12}) under some regularity conditions on the moment function $g(D,\theta)$ and the prior. This is to say that the BGMM posterior probability of the true model converges to 1 with high probability. \\

\noindent 2. We derive an oracle property for the BGMM procedure, which states that the BGMM posterior distribution converges in total variation norm to a normal distribution concentrated on the true model space with an efficient variance, as if the true model were known. This oracle property is the Bayesian analog of the frequentist post-model-selection oracle property of \citet{fanli01}, and is comparable to the Bayesian oracle property proposed by \citet{ir11}. While \citet{ir11} showed this oracle property only for the posterior mean estimator in the normal linear model, our version of Bayesian oracle property studies the global asymptotic concentration behavior of the whole posterior for the general form of moment conditions. We apply BGMM to our motivating examples in Section 1.3 and show that the model selection consistency and oracle property hold under mild regularity conditions on the data and the moments. \\

\noindent 3. Our theory for BGMM allows the number of parameters $p$ to increase with the sample size $n$. This is technically challenging because the number of candidate models $2^p$ will increase exponentially fast with $n$. Although \citet{hong12} and \citet{kim14} have established model selection consistency for BGMM with a fixed number of models, their techniques based on pairwise model comparison are not sufficient for showing the global model selection consistency under our increasing dimensional setup. Our theoretical results accommodate an increasing dimension $p$ that satisfies $p^4 /n \rightarrow 0$ up to some logarithm factors, which is the same as the growth rate in \citet{bc09}, who studied BGMM without model selection.\\

\noindent 4. We present a novel interpretation of the BGMM quasi-posterior, as an approximate posterior conditional on a data summary that is equivalent to the GMM estimator. Particularly for model selection, we derive the convergence rates of Bayes factors for the BGMM method and the fully Bayesian method given the GMM estimator, and show that they have similar asymptotic behavior. Therefore, the model posterior probabilities from BGMM are asymptotically valid and can be used directly for comparing different models.\\

\noindent 5. Our numerical experiments provide practical guidance on the MCMC computation in a complicated setup with $2^p$ candidate models. The previous works on BGMM computation either have not considered the model selection problem (\citealt{yin11}), or have only considered pairwise model comparison using MCMC (e.g. \citealt{kim14}). We implement the reversible jump MCMC algorithm and demonstrate BGMM as a practically feasible and efficient alternative to the frequentist regularization methods.\\

Below we provide some motivating examples that involve the moment condition \eqref{uncond} and can be easily incorporated into the BGMM framwork.

\subsection{Three motivating examples}
\label{sec.examples}
The moment condition model (\ref{uncond}) is much more general than probabilistic models such as the normal linear model and generalized linear models, since one could set the moment function to be $g(D,\theta)=\partial_\theta\ln p(D|\theta)$, where $p(D|\theta)$ is the probability density of $D$.  For example, in the Poisson regression model $D=(Y,X^\top)^\top$ and $Y|X \sim \mathcal{P} (e^{X^\top\theta})$ for some covariates $X$, we can use the moment function $g(D,\theta)=X\left(Y-e^{X^\top \theta}\right)$ for quasi-posterior based inference from \eqref{qpost}, although the likelihood based inference would be more straightforward in this case. In fact, the proposed moment based method has more flexibility when only the lower order moments or quantiles are specified rather than the complete probabilistic model, as described in the following examples.

\begin{example}
\textit{Correlated longitudinal data.} In longitudinal studies, suppose the $j$th observation for the $i$th subject is
a scalar response variable $Y_{ij}$ and a $p$-dimensional covariate vector $X_{ij}$. For simplicity, we assume that
each subject has the same number of observations, i.e. $j=1,\ldots,s$ and $i=1,\ldots,n$. Let $Y_i=(Y_{i1},\ldots,Y_{is})^\top$,
$\X_{i}=(X_{i1},\ldots,X_{is})^\top$, and $\Ep(Y_i|\X_{i})=\mu_i(\theta)$, where $\mu_i(\theta)=(\mu(X_{i1}^\top\theta),\ldots,\mu(X_{is}^\top \theta))^\top$ and $\mu(\cdot)$ is a monotone link function. To account for the heteroscedasticity, we assume the conditional variance of $Y_{ij}$ given $X_{ij}$ is a function of the single index $X_{ij}^\top \theta$, i.e. $\var(Y_{ij}|X_{ij})=\phi(X_{ij}^\top \theta)$.
 Then the frequentist GEE method estimates $\theta$ by solving equations
\begin{equation}\label{glm1}
n^{-1}\sum_{i=1}^n \frac{\partial \mu_i(\theta)^\top}{\partial \theta} \S_i^{-1}(Y_i- \mu_i(\theta))=0
\end{equation}
where $\S_i=\A_i^{1/2}\R\A_i^{1/2}$, $\A_i=\A_i(\theta)=\diag \big\{\phi(X_{i1}^\top\theta),\ldots, \phi(X_{is}^\top\theta)\big\}$ is the diagonal matrix with the conditional variance of $Y$ given $X$ and $\R$ is a working correlation matrix. If we denote the data as $D_i=(Y_i,\X_i)^\top$, then the moment function is defined by
\begin{equation}\label{cldmom}
g(D_i,\theta)=\frac{\partial \mu_i(\theta)^\top}{\partial \theta} \S_i^{-1}(Y_i- \mu_i(\theta)).
\end{equation}
And the moment condition \eqref{uncond} is satisfied.
\end{example}

\begin{example}
\textit{Quantile regression.} Suppose that $Y$ is a continuously distributed response variable, and $X$ is a $p$-dimensional predictor vector for the $\tau$-th quantile ($\tau \in (0,1)$) of $Y$. The conditional quantile function of $Y$ given $X$ is specified by $F^{-1}_{Y|X}(\tau)=X^\top \theta$, where $F^{-1}_{Y|X}$ is the generalized inverse of conditional distribution function of $Y$ given $X$. Then let $D=(Y,X^\top)^\top$ and we can construct $p$ moment functions as
\begin{equation}\label{ivqrmom}
g(D,\theta)=X\left\{1(Y-X^\top \theta\leq 0)-\tau\right\},
\end{equation}
where $1(\cdot)$ is the indicator function.
\end{example}

\begin{example}
\textit{Partial correlation selection.} The partial correlation structure of a $s$-dimensional random vector $Y$ is specified by its precision matrix $\om=\si^{-1}$, where $\si=\Ep\left\{(Y-\Ep Y)(Y-\Ep Y)^\top\right\}$ is the covariance matrix of $Y$. Hereafter without loss of generality, we assume that $Y$ is centered such that $\Ep Y=0$. The partial correlation between the $i$th and the $j$th components of $Y$ is defined by $\rho_{ij}=-\omega_{ij}/\sqrt{\omega_{ii}\omega_{jj}}$, where $\omega_{ij}$ denotes the $(i,j)$th entry of $\om$. $\omega_{ij}=0$ implies zero partial correlation between the $i$th and the $j$th components of $Y$ given all the other components. For multivariate Gaussian random vector, there is an equivalence between the conditional independence and the zero partial correlation. In the general case where multivariate Gaussian assumption is not satisfied, we can still use the second moment of $Y$ to identify the zero entries in $\om$. Let $\theta$ be the vectorized upper triangle part of $\om$. Then we can define the moment function
\begin{equation}\label{pcsmom}
g_{ij}(Y,\theta)=Y_iY_j-(\om^{-1})_{ij},
\end{equation}
for $1\leq i\leq j\leq s$, and the stacked moment vector $g(Y,\theta)$ satisfies \eqref{uncond}. We have $\dimm(\theta)=\dimm(g)=s(s+1)/2=:p$ where $\theta$ is just identifiable. The model selection problem for partial
correlation has been studied in, for example, \citet{drton04}, \citet{jiang04}, etc.
\end{example}

\subsection{Organization of the paper}

The rest of the paper is organized as follows. In Section 2.2, we derive the oracle properties for BGMM model selection based on a set of high level assumptions. In Section 2.3, we discuss the validity of the proposed BGMM quasi-posterior. Section 3 provides the algorithm we use for BGMM and numerical experiments to illustrate the empirical performance of BGMM model selection. Section 4 includes further discussions. We check these assumptions for the three motivating examples in Section \ref{sec.examples} and include the technical proofs of all theorems in the supplementary material. A real data application can be found in the online technical report \citet{lj14}.

\subsection{Some useful notation}
We define some useful notation. Let $|\cdot|_k$ denote the
$L_k$ norm for $k\in [0,\infty]$ and $\|\cdot\|$ be the Euclidean
norm ($L_2$ norm). For any generic square matrix $\C$, let $\underline \lambda(\C)$,
$\bar \lambda(\C)$ denote the smallest and the largest eigenvalues
of a square matrix $\C$. Let $\|\C\|=\sqrt{\bar \lambda(\C^\top\C)}$
be the matrix operator norm. For two stochastic sequence $\{a_n\}$ and
$\{b_n\}$, let $a_n\prec b_n$, $a_n\succ b_n$ and $a_n \asymp b_n$
denote $a_n=o(b_n)$, $b_n=o(a_n)$ and $a_n,b_n$ having the same
order as $n\to \infty$. $a\vee b=\max(a,b)$ and $a\land b=\min(a,b)$. The notations $o_p$ and $O_p$ always refer to
the probability measure $P_{\D}$ of the sample $\D$. We use ``$C$" to
denote any generic constant whose value can change in different places. We use the statement ``the event $A$ happens w.p.a.1 as $n\to \infty$" as an abbreviation for the statement ``the event $A$ happens with $P_{\D}$ probability approaching 1 as $n\to \infty$ ", i.e. $\lim_{n\to \infty} P_{\D}(A)=1$.\\

\section{Theoretical Properties of Bayesian GMM Model Selection}
The Bayesian model selection problem has been extensively studied, but mostly for normal linear regression models and generalized linear models. See for example, \citet{chipman01}, \citet{smith96}, \citet{ir05}, \citet{jiang07}, \citet{liang08}, \citet{johnson12}, \citet{liang13}, etc. Our Bayesian model selection is substantially different from all these papers. Instead of having a probabilistic model such as the simple normal linear model, we work with the moment conditions \eqref{uncond} and do model selection using BGMM. Our true parameter $\theta_0$ is the unique solution of \eqref{uncond} and possibly lies in a lower dimensional subspace of the whole parameter space $\Theta\subseteq \mathbb{R}^p$. We restrict $\Theta$ to be a compact and connected set in $\mathbb{R}^p$, with finite $L_2$ radius $R=\sup_{\theta\in \Theta}\|\theta\|$ for some large constant $R>0$.\\

Without loss of generality, in the following we will consider models generated by all the possible coordinate subspaces of $\mathbb{R}^p$, which leads to a total of $2^p$ different models $\M$ and the parameter space partition $\Theta=\bigcup_{|\M|\leq p} \Theta(\M)$. Let $k=|\M|$ ($0\leq k\leq p$) be the size of a generic model $\M$, which is the number of nonzero components in any $\theta\in \M$. Suppose $\M_0$ is the true model space that contains $\theta_0$, and $k_0=|\M_0|$ is the dimension of $\theta_0$. For a given model $\M$ and a generic $\theta$, let $\theta=(\theta_1^\top, \theta_2^\top)^\top$ where $\theta_1 \in \mathbb{R}^k$ and $\theta_2 \in \mathbb{R}^{p-k}$ correspond to the components that lie in and outside $\Theta(\M)$, respectively. So $\theta_2=0$ if $\theta \in \Theta(\M)$. We emphasize that the meaning of subscripts ``1" and ``2" can change with the model index $\M$.\\

For such a model selection setup, the prior distribution can be written in the hierarchical structure $\pi(\theta)=\sum_M \pi(\theta|\M)\pi(\M) =\sum_{\M,k} \pi(\theta|\M) \pi(\M||\M|=k)\pi(k)$ for $k=0,1,\ldots,p$. If a model $\M$ does not contain all the nonzero components for a given $\theta$, then $\pi(\theta|\M)=0$. We assume that each $\pi(\theta|\M)$ has a density function. For two different models $\M_1$ and $\M_2$, the (quasi-) Bayes factor of $\M_1$ with respect to $\M_2$ is defined as
\begin{equation}\label{bfac}
\BF[\M_1:\M_2]=\frac{q(\D |\M_1)}{q(\D|\M_2)}=\frac{\int_{\Theta(\M_1)}
q(\D|\theta,\M_1)\pi(\theta|\M_1)\ud\theta}{\int_{\Theta(\M_2)}
q(\D|\theta,\M_2)\pi(\theta|\M_2)\ud\theta}
\end{equation}
and accordingly the (quasi-) posterior odds is the product of the Bayes
factor and the prior odds
\begin{equation}\label{podds}
\PO[\M_1:\M_2]=\frac{q(\D |\M_1)}{q(\D|\M_2)}\cdot\frac{\pi(\M_1)}{\pi(\M_2)}=\frac{\int_{\Theta(\M_1)}
q(\D|\theta,\M_1)\pi(\theta|\M_1)\ud\theta}{\int_{\Theta(\M_2)}
q(\D|\theta,\M_2)\pi(\theta|\M_2)\ud\theta} \cdot \frac{\pi(\M_1)}{\pi(\M_2)}
\end{equation}

The model selection consistency we are going to establish is the global model selection consistency (\citealt{johnson12}), in the sense that asymptotically the true model $\M_0$ will not only be the MAP model (\textit{maximum a posteriori}) but also have posterior probability tending to 1. Equivalently, we will show that the sum of all posterior odds $\PO[\M:\M_0]$ with $\M\neq \M_0$ converges to zero in probability. This strongest mode of consistency implies that the posterior mass will be concentrated around the true model and most of the $2^p$ models receive negligible probabilities. This is a desirable property in practice for interpretation, since commonly used Bayesian estimation procedures such as model averaging will then involve only a few models instead of many candidate models.

\subsection{Assumptions}
The set of assumptions below follows closely the set of conditions for Z-estimation in \citet{bc09}. They are high level assumptions imposed on the data generation process, the model parameters, the moment conditions and the priors. For a specific model, these assumptions are not necessarily in the most general form, but they do cover a wide class of moment condition models in practice and are sufficient for illustrating the theoretical properties of BGMM. \\

For the data generation process and the true parameter $\theta_0$, we make the following assumptions.\\

\noindent \textit{Assumption 1} (Data Generation Process) $\{D_i, i=1,\ldots,n\}$ is an i.i.d. sequence. $\Ep g(D,\theta_0)=0$ for some $\theta_0\in \Theta$. $\Theta$ is a compact and connected set with $L_2$ radius $R$ for some large constant $R>0$, and it  contains an open neighborhood of $\theta_0$.\\

\noindent \textit{Assumption 2} (Dimension) Let $\dimm(\theta)=p$ and $\dimm(g)=m$. Assume that $p\leq m$, $p\asymp m$, $p^4\ln^2 n/n\to 0$ and $p^{2+\alpha}\ln n/n^\alpha \to 0$, where $\alpha$ is defined in Assumption 4. \\

\noindent \textit{Assumption 3} (Beta-min) Let $\e=\sqrt{p/n}$. Assume $1\succeq \min_{j\in \M_0}|\theta_{0,(j)}|\succ \sqrt{\ln n}\e$, where $\theta_{0,(j)}$'s for $j\in \M_0$ denote the nonzero components of the true parameter $\theta_0$.\\

The i.i.d. assumption in Assumption 1 can be possibly relaxed to a weakly dependent stationary process using more involved techniques. The compactness assumption for the parameter space $\Theta$ is standard and mainly for technical convenience, and it can be relaxed to the full space of $\mathbb{R}^p$ if we can control the tail behavior of the prior (see the discussion after Assumptions 7 and 8).  Assumption 2 allows increasing dimension $p$, and the growth rate of $p$ is comparable with those in \citet{bc09}, \citet{cho13}, \citet{wang11}, \citet{leng12}, etc.  The beta-min condition in Assumption 3 is commonly used in the frequentist GEE literature (see e.g. \citealt{wang12}, \citealt{leng12}, \citealt{cho13}). It gives the minimal magnitude of nonzero coefficients that could be detected by BGMM.\\

Let $B_0(\epsilon)=\{\theta\in \Theta: \|\theta-\theta_0\|<\epsilon\}$ for any $\epsilon>0$. We make the following assumptions on the moment conditions. \\

\noindent \textit{Assumption 4} (Moment) (i) The moment function $g(D,\theta)$ satisfies the continuity property
 $$\sup_{\eta\in \mathbb{R}^m,\|\eta\|=1} \left[\Ep \left\{(\eta^\top (g(D,\theta)-g(D,\theta_0)))^2\right\}\right]^{1/2}\leq O\left((\sqrt{p}\|\theta-\theta_0\|)^{\alpha}\right),$$
  uniformly in $\theta\in \Theta$ for some constant $\alpha \in (0,1]$.\\
(ii) The class of functions $\mathcal {F} =\left\{\eta^\top(g(D,\theta)-g(D,\theta_0)),\theta\in \Theta, \eta\in \mathbb{R}^m, \|\eta\|=1\right\}$ has an envelope function $F$ almost surely bounded in $L_2$ norm $\|\cdot\|_{P_D,2}$ as order $O(\sqrt{p})$. The $L_2$ uniform covering number $N\left(\epsilon\|F\|_{P_D,2},\mathcal{F},L_2(P_D)\right)$ satisfies that for any small $\epsilon>0$,
$$\ln N\left(\epsilon\|F\|_{P_D,2},\mathcal{F},L_2(P_D)\right)=O\Big(p\ln \Big(\frac{n}{\epsilon}\Big) \Big).$$
\noindent \textit{Assumption 5} (Linearization) (i) $\|\Ep g(D,\theta)\|\geq \delta_0 \land \left(\delta_1
\|\theta-\theta_0\|\right)$ uniformly on $\Theta$ for some positive constants $\delta_0,\delta_1$. \\
(ii) $\G:= \nabla_{\theta}\Ep g(D,\theta_0)$ exists, and the eigenvalues of $\G^\top \G$ are bounded from below and above as $n\to \infty$.\\
(iii) $\H(\theta):= \nabla^2_{\theta\theta^\top }\Ep g(D,\theta)$ exists for
$\theta\in B_0(C\epsilon_n)$, and uniformly over
$\theta\in B_0(C\epsilon_n)$ for any fixed $C>0$,
$\sup_{\|u\|=1,\|v\|=1,u,v \in \mathbb{R}^p}\left\|\H(\theta)(u,v)\right\|=O(\sqrt{p})$.\\

Assumptions 4 and 5 on moment function $g(D,\theta)$ parallel the conditions ZE.1 and ZE.2 in \citet{bc09} respectively. The continuity index $\alpha$ in Assumption 4(i) satisfies $\alpha=1$ for the mean regression, such as the examples of correlated longitudinal data and partial correlation selection, and $\alpha=1/2$ for the quantile regression model. The entropy condition in Assumption 4(ii) controls the complexity of the class of moment functions $g(D,\theta)$. Assumption 5(i) guarantees the point identification of the true parameter $\theta_0$, and part (ii) and (iii) impose mild assumptions on the first and second derivatives of $\Ep g(D,\theta)$ around $\theta_0$. These regularity conditions are mainly used to derive large deviation bounds via empirical process results, and they will be verified later for our motivating examples. Note that unlike \citet{wang12}, \citet{leng12} and \citet{cho13}, we do not require the moment function $g(D,\theta)$ itself to be differentiable. This allows more general applications to discontinuous $g(D,\theta)$, such as in the case of quantile regression.\\

\noindent \textit{Assumption 6} (Variance) $\V_n$ is a positive definite matrix for all $n$, and converges in the matrix operator norm to $\V=\var\left\{g(D,\theta_0)\right\}$. The eigenvalues of $\V_n$ and $\V$ are bounded below and above for some positive constants $\underline \lambda$ and $\bar \lambda$ w.p.a.1 as $n\to \infty$.\\

Assumption 6 assumes that the positive definite weighting matrix $\V_n$ is a consistent estimator of the covariance matrix of $g(D,\theta)$ at $\theta_0$, similar to the preliminary estimator of the optimal weighting matrix used in the two step GMM estimation. Although this consistency of $\V_n$ to $\V$ is not required for the model selection consistency, it is necessary for the valid posterior inference such as posterior credible sets. Essentially $\V_n$ needs to satisfy the generalized information inequality (\citealt{ch03}), such that the LIL asymptotically satisfies the second Bartlett identity as a true likelihood function does. Such consistent estimator $\V_n$ usually exists for all our motivating examples.\\

Finally we impose the following assumptions on the prior.\\

\noindent \textit{Assumption 7} (Prior on $\theta$) (i) $\pi(\theta|\M)$ has a density function restricted to $\Theta$, and is bounded above by a constant $c_{\pi}$ uniformly over all model spaces $\M$. \\
(ii) Suppose $\theta=(\theta_1^\top,\theta_2^\top)^\top$ is decomposed according to the model $\M$. Then uniformly over all models $\M\supseteq \M_0$, for any given $C>0$, $|\ln\pi(\theta_1|\M)- \ln\pi(\theta_{0,\M,1}|\M)|=o(1)$ as $n\to\infty$ if $\theta=(\theta_1^\top,0^\top)^\top\in \Theta(\M) \cap B_0(C\e)$, where $\theta_{0,\M,1}$ is the subvector of the true parameter $\theta_0$ restricted to $\Theta(\M)$. \\
(iii) Uniformly for all models $\M \supseteq \M_0$, there exist constants $c_0,c_1>0$, such that $\pi(\theta_0|\M)\geq e^{-c_0|\M|}$ and $|\ln\pi(\theta_{0,\M,1}|\M)- \ln\pi(\theta_0|\M_0)|\leq c_1(|\M|-|\M_0|)$.

\noindent \textit{Assumption 8} (Prior on models) The model prior $\pi(\M)$ satisfies: \\
(i) $\lim_{n\to\infty} \sup_{\{\M:\M \supset \M_0\} }
\left(\pi(\M)/\pi(\M_0)\right)\left(p/\sqrt{n}\right)^{|\M |-| \M_0|} =0$.\\
(ii)  $\sup_{\{\M: \ \M_0\backslash \M \neq \emptyset\}} \pi(\M)/\pi(\M_0) \preceq e^{r_1p\ln n}$ for some constant $r_1>0$ as $n\to\infty$. \\

Our Assumption 7(i) that the prior is restricted on a compact set $\Theta$ is mainly for technical simplification, which can be relaxed as long as the tail probability of $\pi(\theta|\M)$ decays sufficiently fast on each model $\M$. The idea is that the one can divide the possibly noncompact support into a compact set $\Theta_n$ with radius increasing sufficiently slowly with $n$, and let the prior mass outside $\Theta_n$ be negligible as in the large $n$ asymptotics (\citealt{jt08}). Other than this, Assumption 7(i) is mild and encompasses most of the commonly used priors truncated on $\Theta$. Assumption 7 (ii) and (iii) for the priors on parameters $\pi(\theta|\M)$ are satisfied by, for example, a uniform prior on the model space $\Theta(\M)$, or a truncated multivariate normal prior on $\Theta(\M)$.

Assumption 8 requires that the models larger than the true model $\M_0$ do not receive overly large prior mass, and the prior on the true model cannot be exponentially small compared to any other models. This is automatically true when $p$ does not increase with $n$, provided that $\pi(\M_0)$ is a positive constant. With increasing dimensions, these requirements can be satisfied by, for example, a prior where each coordinate enters the model independently with a fixed probability $\nu\in (0,1)$, which includes the uniform prior as a special case if $\nu=0.5$.  Other examples include priors that propose a model size $|\M|$ according to Poisson or geometric distributions upper truncated at $p$, while all models of the same size are equally likely. A detailed verification of Assumption 8 for these priors can be found in Section 4 of the supplementary material.

\subsection{Oracle Properties of BGMM}
With all these assumptions, we now state the main results as follows. The proof of Theorem \ref{main} is given in the supplementary material.

\begin{theorem} \label{main}
Suppose Assumptions 1-8 hold. Then\\
(i) (Model Selection Consistency)
\begin{displaymath}
q(\M_0|\D)\to 1, \text{   w.p.a.1 as $n\to \infty$  }
\end{displaymath}
that is, the quasi-posterior probability of the true model converges to 1, w.p.a.1  as $n\to \infty$.

\vspace{+0.3cm}

\noindent (ii) (Posterior Asymptotic Normality) Given a model $\M$, let $\G_{\M}$ be the submatrix of the derivative matrix $\G$ with respect to the subvector $\theta_1$ in $\theta=(\theta_1^\top,\theta_2^\top)^\top$. Let $ \bar \theta_{\M_0,1} = \theta_{0,\M_0,1}-(\G_{\M_0}^\top\V_n^{-1}\G_{\M_0})^{-1}\G_{\M_0}^\top \V_n^{-1}\bar g(\D,\theta_0)$, where $\theta_{0,\M_0,1}$ is the subvector of $\theta_0$ restricted to $\Theta(\M_0)$. Then w.p.a.1 as $n\to \infty$,
\begin{displaymath}
\sup_{A\subseteq \Theta} \Bigg| \int_A q(\theta|\D)\ud \theta - \int_{A\cap \Theta(\M_0)} \phi\left(\theta_1;\bar \theta_{\M_0,1}, (\G_{\M_0}^\top \V_n^{-1} \G_{\M_0})^{-1}/n\right)\ud \theta_1\Bigg|\to 0,
\end{displaymath}
where $\phi(\cdot;\mu,\si)$ is the normal density with mean $\mu$ and covariance matrix $\si$, and $\theta=(\theta_1^\top,\theta_2^\top)^\top$ is decomposed according to the true model $\M_0$.
\end{theorem}

Part (i) of Theorem \ref{main} establishes the global model selection consistency of BGMM, similar to previous Bayesian results from \citet{johnson12} and \citet{liang13} for the normal linear model and the generalized linear models. Based on the BGMM posterior, the zero components of the true parameter $\theta_0$ are estimated to be zero with $P_{\D}$ probability approaching 1. It also implies that asymptotically the MAP model $\hat \M$ converges to the true model $\M_0$ in $P_{\D}$-probability. This parallels the frequentist model selection results via penalization for moment based models and estimating equations (\citealt{wang12}, \citealt{leng12}, \citealt{cho13}, etc.)\\

Part (ii) of Theorem \ref{main} establishes an asymptotic normality result, in the sense that the total variation difference between the BGMM posterior measure and a $k_0$-dimensional normal distribution concentrated on the true model converges to zero in probability as the sample size increases. This is a direct extension of the Bayesian CLT result in \citet{ch03} and \citet{bc09} from a single full model space to the joint of all submodel spaces. Because the BGMM posterior is a mixture distribution on $2^p$ model spaces $\Theta(\M)$ with different dimensions, we do not present result using the $L_1$ distance between two densities $q(\theta|\D)$ and $\phi(\theta_1;\bar \theta_{\M_0,1}, (\G_{\M_0}^\top \V_n^{-1} \G_{\M_0})^{-1}/n)$. The asymptotic mean of the normal distribution $\bar \theta_{\M_0, 1}$ is the first order approximation to $\hat \theta_{\M_0,1} = \arg\min_{\theta \in \Theta(\M_0)} \bar g(\D,\theta)^\top \V_n^{-1} \bar g(\D,\theta)$, i.e. the GMM estimator restricted to the subspace $\Theta(\M_0)$. Furthermore, given Assumption 6, the generalized information equality is satisfied (\citealt{ch03}), and the asymptotic variance of the limiting normal distribution is the same as the corresponding frequentist variance of the GMM estimator $\hat \theta_{\M_0, 1}$.
\begin{remark}
{\it Bayesian oracle property.} The conclusion of Theorem \ref{main} can be written heuristically as follows: Let $\theta=(\theta_1^\top,\theta_2^\top)^\top$ be decomposed according to the true model $\M_0$, then

\vspace{+0.3cm}

(i) $\theta_{2}|\D\approx 0$ w.p.a.1 as $n\to\infty$ ;

\vspace{+0.2cm}

(ii) $q(\theta_{1}|\D)\approx \mathcal{N} \left(\bar \theta_{\M_0,1},(\G_{\M_0}^\top \V^{-1} \G_{\M_0})^{-1}/n\right)$ w.p.a.1 as $n\to\infty$ .

\vspace{+0.3cm}

The zero components in $\theta_0$ are estimated to be zero given the data $\D$ with large probability, and the nonzero components in $\theta_0$ almost follow a normal distribution centered at the first order approximation to the GMM estimator under the model $\M_0$, with the same optimal GMM asymptotic variance matrix, as if the true model $\M_0$ were known. We call this the {\it Bayesian oracle property} for model selection, which resembles the frequentist oracle property for penalized likelihood in \citet{fanli01}. Theorem \ref{main} guarantees that the BGMM posterior will automatically identify the unknown true model, and automatically converges to an asymptotic normal distribution centered around the unknown true parameter with the optimal GMM variance, as if the true model were known. Compared to the oracle property in \citet{ir11}, our version is much stronger in two aspects: 1. Our model assumptions are based on the general form of moment conditions \eqref{uncond} and are therefore more general than the normal linear regression model in \citet{ir11}; 2. Our oracle property characterizes the overall shrinkage of posterior distribution to an asymptotic normal distribution on the true model, while \citet{ir11} only considered the asymptotics of the posterior mean estimator.
\end{remark}

\begin{remark}\label{oraclecenter}
We explain why the oracle center $\bar\theta_{\M_0,1}$ (of the asymptotic normal approximation to the quasi posterior) is a desirable result. Roughly speaking, this oracle center will be often close to the unknown nonzero components of the true parameter $\theta_0$ in large samples, since their difference has the order $O_p\left(\sqrt{p/n}\right)$. This oracle center is also similar to the center of Bayesian CLT in \citet{bc09} and it applies to all our motivating examples in Section \ref{sec.examples}. Furthermore, in many cases we have the higher order approximation from $\bar \theta_{\M_0,1}$ to the GMM estimator $\hat \theta_{\M_0,1}$, with the difference $\|\bar \theta_{\M_0,1}-\hat \theta_{\M_0,1}\|=O_p(n^{-1})$, following the stochastic expansion of GMM estimator in \citet{neweysmith04}. When this high order approximation holds, the oracle center $\bar \theta_{\M_0,1}$ in Theorem \ref{main} (ii) is equivalent to and can be replaced by the oracle GMM estimator $\hat \theta_{\M_0,1}$.
\end{remark}

\begin{remark}\label{comparehong}
Our work in model selection of BGMM may be regarded as a more detailed study of a special case of \citet{hong12}. They have considered model selection in a more general framework, which allows general objective functions, including the GMM and GEL criterion functions. In addition, they allow multiplicity in the set of ``best models" which could be mutually nonnested (see their Section 4.2.2). Their results indicate that model selection consistency can hold in the nested case but fail in the nonnested case.
Regarding such opposite conclusions in these two cases, we have benefited from an anonymous referee on clarifying this point, who noted that consistent model selection has two meanings in \citet{hong12}. The first meaning is that a consistent model selection procedure selects the set of ``best" models w.p.a.1. The second meaning is that if there is multiplicity in the set of best models, a consistent model selection procedure should pick the most parsimonious model among the best models w.p.a.1. Our result in model selection consistency of BGMM is obtained in the nested case, since we have considered all $2^p$ coordinate subspaces of $\mathbb{R}^p$. Therefore, the oracle property of BGMM in Theorem 1 fulfills both two meanings of consistent model selection described in \citet{hong12}.

Although the nested case we have considered is not as general as \citet{hong12}, we have allowed the dimension $p$ to increase with $n$, which is new for BGMM model selection and also technically challenging. Because the number of candidate models is $2^p$, which increases exponentially fast in $p$ and hence in $n$, the previously studied pairwise model comparison using posterior odds or Bayes factors between one candidate model and the true model (such as \citealt{hong12}, \citealt{kim14}) is insufficient to show the global model selection consistency. In addition to the increasing dimensionality and the more detailed study on the limiting distribution, we will also discuss below the asymptotic validity and interpretation of the BGMM quasi-posterior, which is new in the literature.
\end{remark}

\subsection{Asymptotic Validity of the BGMM Posterior}
As shown in \citet{kim02}, the limited information likelihood we have used for BGMM provides a large sample approximation to the true likelihood function of $\theta$ given the moment restrictions $\Ep g(D,\theta)=0$. One may ask about how well this approximation could be. For the validity of usual Bayesian inference, such as constructing the Bayesian credible sets, it is necessary and sufficient to impose Assumption 6 that $\V_n$ consistently estimates $\V$, i.e. $\V_n$ satisfies the generalized information equality as in \citet{kim02} and \citet{ch03}. However, due to the limited information contained in $\Ep g(D,\theta)=0$, in general one cannot expect the LIL $q(\D|\theta)$ to coincide with the true likelihood function $p(\D|\theta)$. Instead the quasi-posterior $q(\theta|\D)$ can be used to approximate the posterior of $\theta$ given some summary statistic from the sample. Let $\hat \theta$ be the minimizer of the GMM criterion function $\bar g(\D,\theta)^\top \V^{-1} \bar g(\D,\theta)$ over the full $p$-dimensional model space. So $\hat \theta$ is implicitly a statistic of the sample $\D$, and it does not depend on $\theta_0$ and the unknown true model $\M_0$. Since
the asymptotic center of the BGMM posterior is the first order approximation to the GMM estimator, one can expect that the LIL $q(\D|\theta)$ approximates the density $p(\hat \theta|\theta)$ of $\hat \theta$. Accordingly, the BGMM posterior $q(\theta|\D)$ approximates the posterior $p(\theta|\hat \theta)$ of $\theta$ given $\hat \theta$, at least asymptotically. In the following, we formalize this idea and show more general results under the model selection setup. \\

For two generic models $M_1$ and $M_2$, we define the Bayes factor based on $p(\hat \theta|\theta)$ as
\begin{equation}
\bfh[\M_1:\M_2]=\frac{p(\hat \theta |\M_1)}{p(\hat \theta|\M_2)}=\frac{\int_{\Theta(\M_1)}
p(\hat \theta |\theta)\pi(\theta|\M_1)\ud\theta}{\int_{\Theta(\M_2)}
p(\hat \theta |\theta)\pi(\theta|\M_2)\ud\theta} \nonumber
\end{equation}

For theory development, in this section we focus on the situation with a nonincreasing dimension $p$. We make the following extra assumption.\\

\noindent \textit{Assumption 9} (i) $\dimm(\theta)=p$ and $1\leq p \leq \bar p$, for some large fixed integer $\bar p$.\\
\noindent (ii) $\min_{j\in \M_0}|\theta_{0,(j)}|\geq \underline \theta$ for some small constant $\underline \theta >0$. \\
\noindent (iii) Let $\V(\theta)=\var\left\{g(D,\theta)\right\}$ and $\G(\theta)=\nabla_{\theta}\Ep g(D,\theta)$. Then the elements of $\V(\theta)$ and $\G(\theta)$ are continuous functions of $\theta$, and the eigenvalues of $\G(\theta)^\top\G(\theta)$ and $\V(\theta)$ are uniformly bounded below and above for all $\theta \in \Theta$.\\
\noindent (iv) For any two models $\M_1$ and $\M_2$, there exists a constant $r>0$ such that $\pi(\M_2)/\pi(\M_1)\leq r$.\\
\noindent (v) $\big\|\hat \theta-\bar \theta\big\|=O_p(1/n)$, where $\hat \theta$ is the GMM estimator on the full model space, and $\bar \theta = \theta_0 - (\G^\top\V^{-1}\G)^{-1} \G^\top\V^{-1}\bar g(\D,\theta_0)$.\\

The strengthened beta-min condition in (ii) is to emphasize the difference between the models that make the type I error and the type II error. According to theorems we are going to present below, the models in the former group have an exponentially small $\BF[\M:\M_0]$, while the models in the latter group have a polynomially small $\BF[\M:\M_0]$. This is also the essential behavior from the Bayesian hypothesis test, which favors the true alternative hypothesis more. We will show that similar behavior is also shared by $\bfh[\M:\M_0]$, and hereby establish a correspondence between the BGMM method and the exact Bayasian method given $\hat \theta$.\\

Part (iii) assumes the continuity of the matrices in $\theta$ and also the uniform bound for eigenvalues. This is a mild assumption given the compactness of $\Theta$. Part (iv) has strengthened Assumption 8 and required that no model should be assigned extremely large or small prior. Part (v) is about the high order approximation of $\bar\theta$ to the GMM estimator $\hat \theta$ on the full model space, similar to the discussion in Remark \ref{oraclecenter}, which usually holds when the moment condition $g(D,\theta)$ is continuously differentiable in $\theta$ (\citealt{neweysmith04}) and hence may not apply to the example of quantile regression.\\

Let $\F(\theta)$ be a $p\times p$ matrix such that $\F(\theta)^\top\F(\theta)= \G(\theta)^\top\V(\theta)^{-1}\G(\theta)$. Define $Z=\sqrt{n}\F(\theta)(\hat \theta-\theta)$. Then $Z$ is asymptotically $p$-dimensional standard normal if the true parameter is $\theta_0=\theta$. We impose the following high level assumption on the difference between the exact density function $p_{Z}(z)$ of $Z$ and the normal density.\\

\noindent \textit{Assumption 10} (Uniform Bound) As $n\to \infty$,
    \begin{equation}
\sup_{\theta \in \Theta}\sup_z \left(1+\|z\|^{p+1} \right)\left| p_{Z}(z|\theta) - \phi(z;0,\I_p) \right| = \tau_n ,\nonumber
    \end{equation}
where $\tau_n=o(1)$ does not depend on $z$ and $\theta$, and $\I_p$ is the $p\times p$ identity matrix.\\

Assumption 10 claims that the difference between the density of the normalized GMM estimator $Z$ and its asymptotic limit of normal density can be uniformly bounded by an integrable function $c(\|z\|)=1/(1+\|z\|^{p+1})$, and the uniformity is for both the value of $z$ and the parameter $\theta$ in the compact space $\Theta$. This is a high level condition that originates from the Condition E in \citet{yc04}. We do not intend to give a full proof of it under low level assumptions, but we explain why it is a reasonable assumption below.\\

Consider the case where the $(p+1)$-th moment of $g(D,\theta)$ exists. To show Assumption 10, we proceed in several steps. First, under similar regularity conditions that make Assumption 9(v) hold, one can see that for a fixed $\theta$, the density of $Z$ is asymptotically uniformly close to the density of the normalized first order approximation $\bar Z = \sqrt{n}\F(\theta) (\bar \theta-\theta)$, up to the order $O(1/\sqrt{n})$, where $\bar \theta$ is defined in Theorem \ref{main} (ii). See \citet{kr12, kr13} for the formal proofs of a general class of nonlinear estimators, which can also be applied to the GMM estimator. Second, due to the sample average form of $\bar \theta$ and hence $\bar Z$, one can use Proposition 1 in \citet{yc04} and take $c(x)=1/(1+x^{p+1})$. This proposition provides a bound for the difference between the density of $\bar Z$ and its limiting normal density, which holds uniformly for all $\theta \in \Theta$. Its proof involves the techniques in Chapter 19 of \citet{brr} about the uniform convergence of continuous characteristic functions in the compact set $\Theta$. Third, one can show that in Proposition 1 of \citet{yc04}, the summation of the Edgeworth series beyond the leading normal density term has the order $o_p(1)$. This is due to the finite moments of $g(D,\theta)$ up to the $(p+1)$-th order, as well as the boundedness of multivariate Hermite polynomials. Finally we combine all these pieces and conclude that the uniform deviation in Assumption 10 holds with some $\tau_n=o(1)$.\\

The next theorem provides a comparison between the convergence rates for the Bayes Factors $\bfh[\M:\M_0]$ from the likelihood given the statistic $\hat \theta$ with $\BF[\M:\M_0]$ from the BGMM method.

\begin{theorem} \label{bfqhat}
(Equivalence of Bayes Factors) Suppose Assumptions 1-10 hold, and the true model size is $|\M_0|=k_0$. Then under the same prior $\pi(\theta|\M)$ and $\pi(\M)$, w.p.a.1 as $n\to \infty$,\\
\noindent (i) For any model $\M$ with $\M \supseteq \M_0$,
\begin{align*}
&\frac{\BF [\M:\M_0]}{\bfh [\M:\M_0]}  \to 1; \\
& \BF [\M:\M_0] \asymp \bfh [\M:\M_0] \asymp n^{-\frac{|\M|-k_0}{2}} \succeq n^{-\frac{p-k_0}{2}};
\end{align*}
\noindent (ii) For any model with $\M$ with $\M_0\backslash \M\neq \emptyset$, there exists a constant $C>0$, such that
\begin{align*}
&\BF[\M:\M_0]  \leq   \exp\big(-Cn\underline \theta^2\big) \prec n^{-\frac{p-k_0+1}{2}};\\
&\bfh[\M:\M_0] \leq\exp\big(-Cn\underline \theta^2\big) \vee \tau_n n^{-\frac{p-k_0+1}{2}} \prec  n^{-\frac{p-k_0+1}{2}}.
\end{align*}
\end{theorem}

Theorem \ref{bfqhat} compares the Bayes factors from BGMM and $p(\hat \theta|\theta)$, for the models that make a type I error (Part ii) and a type II error (Part i). The theorem has at least two direct implications. First, for the models that make a type II error (including more components of $\theta$ than necessary), the Bayes factors are asymptotically equal, and both decrease polynomially in the sample size $n$. The polynomial index reflects the difference in dimensions between $\M$ and $\M_0$. Second, for the models that make a type I error (missing at least one nonzero component in $\theta_0$), the Bayes factor from BGMM decreases exponentially fast in $n$. For the Bayes factor from $p(\hat \theta|\theta)$, we have obtained an upper bound for its rate, which also depends on the rate $\tau_n$ in Assumption 10 besides the usual exponential rate. Because $\tau_n=o(1)$ by Assumption 10, we can see clearly that there exists at least a $n^{-1/2}$ gap between the convergence rates of Bayes factors for the models with type I and type II errors. The threshold rate is $n^{-(p-k_0)/2}$, which depends on the unknown dimension $k_0$ of the true model $\M_0$. In general, the posterior probabilities of the models with type I errors converge faster to zero than the posterior of the models with type II errors.\\

This extra part $\tau_n n^{-(p-k_0+1)/2}$ for the Bayes factor in (ii) arises mainly technically from our Assumption 10. Usually, the order $\tau_n=o(1)$ in Assumption 10 is tight and cannot be improved. However, we conjecture that it could be removed by making stronger assumptions on the density function $p(\hat \theta|\theta)$, or the density $p_Z(z|\theta)$ of the normalized statistic $Z$. For example, one can assume that $p_Z\left(\sqrt{n}\F(\theta)(\hat \theta-\theta)\big|\theta\right)$ decreases exponentially fast in $n$ as $\theta$ moves away from the true parameter $\theta_0$. However, we note that usually it is difficult to verify such assumptions because $\hat \theta$ does not have an explicit density, except for a few special cases where $\hat \theta$ comes from the exponential family. We also note that such compromised rate also shows up in Lemma 1 of \citet{marin13}, where they studied the convergence rates of Bayes factors given a general statistic. Although typically one cannot obtain the exact form of the density $p(\hat \theta|\theta)$ and its posterior $p(\theta|\hat \theta)$, Theorem \ref{bfqhat} provides some evidence that in the asymptotic sense, the Bayes factors from BGMM behave very similarly to the Bayes factors from $p(\hat \theta|\theta)$, indicating the validity of using $\BF[\M:\M_0]$ for model selection purpose.

\begin{remark}
In principle, Theorem \ref{bfqhat} provides a guideline to interpret the BGMM posterior probabilities of different models. For simplicity, suppose that all models receive the uniform prior $\pi(\M) \propto 1$. Then since $q(\M_0|\D)\to 1$, the posterior
$q(\M|\D)$ is roughly the same as $\BF[\M:\M_0]$. Because of the gap between the polynomial rate in (i) and the exponential rate in (ii), we can choose any rate in between as a threshold, for example $e^{-\sqrt{n}}$. If a model $\M$ has $q(\M|\D) \geq e^{-\sqrt{n}}$, then we can approximately regard $q(\M|\D)$ as the true posterior probability $p(\M|\hat \theta)$ and consider $\M$ as a model with nonnegligible posterior. This fits well with the common practice that we rank the models according to their posterior probabilities and only study the models on top of the list.
\end{remark}

Based on Theorem \ref{bfqhat}, we can further show that the BGMM posterior $q(\theta|\D)$ and the exact posterior $p(\theta|\hat \theta)$ are close in the total variation distance asymptotically.

\begin{theorem} \label{qphat}
Suppose Assumptions 1-10 hold. Let the full model be $\M_{\full}$. Then under the same prior $\pi(\theta|\M)$ and $\pi(\M)$, w.p.a.1 as $n\to \infty$, \\
\noindent (i) (Model Selection Convergence Rate)
If $\M_0\neq \M_{\full}$, then
\begin{align*}
& \frac{q(\M:\M\neq \M_0 |\D) }{p(\M:\M\neq \M_0 | \hat \theta)} \to 1; \\
& q(\M:\M\neq \M_0 |\D)  \asymp   p(\M:\M\neq \M_0 | \hat \theta) \asymp n^{-\frac{1}{2}} \to 0;
\end{align*}
If $\M_0 = \M_{\full}$, then for some constant $C>0$,
\begin{align*}
& q(\M:\M\neq \M_0 |\D)  \leq  \exp\big(-Cn\underline \theta^2\big) \to 0;\\
& p(\M:\M\neq \M_0 |\hat \theta) \leq  \exp\big(-Cn\underline \theta^2\big) \vee \tau_n n^{-\frac{p-k_0+1}{2}} \to 0.
\end{align*}
\noindent (ii) (Asymptotic Posterior Validity)
\begin{displaymath}
\sup_{A \subseteq \Theta} \Bigg|\int_A q(\theta|\D)\ud \theta - \int_A p(\theta | \hat \theta)\ud \theta \Bigg|\to 0.
\end{displaymath}
\end{theorem}

Part (i) of the theorem is a direct corollary from Theorem \ref{bfqhat}. It implies that the posterior probability of the true model $\M_0$ converges to 1 at exactly the same rate using either the BGMM or $p(\hat \theta|\theta)$, when the true model is a strict submodel of the full model. When the true model is exactly the same as the full model, we have only upper bounds for the model selection convergence rates, as they usually decrease exponentially fast, but again the rate is compromised by $\tau_n$ from Assumption 10 when we consider $p(\M:\M\neq \M_0 |\hat \theta)$. In either scenario, we have the global model selection consistency for both the BGMM posterior and the posterior given $\hat \theta$.\\

Part (ii) gives the asymptotic validity of the BGMM posterior, in the sense that it provides the same asymptotic inference as the exact posterior of $\theta$ given the statistic $\hat \theta$. It has the immediate implication that the posterior credible sets for the parameters constructed from the BGMM posterior are asymptotically valid. It is worth noting that the conclusion of (ii) is only related to the global model selection consistency for both the BGMM posterior and the posterior of $p(\theta|\hat\theta)$, and does not depend on the exact convergence rates of model selection in Part (i). In fact, Part (ii) also holds for general non-model selection prior $\pi(\theta)$ as long as it has a bounded continuous density on $\Theta$. This can be obtained from combining Theorem 1 in \citet{ch03} and Theorem 2 in \citet{yc04} (where $T_n=\hat \theta$), under Condition E in \citet{yc04}. Our proof of Theorem \ref{qphat} follows a similar route by using Assumption 10, but has accommodated the nature of model selection priors $\pi(\theta,\M)= \pi(\theta|\M)\pi(\M)$.

\begin{remark}
We have discussed the asymptotic closeness of the BGMM posterior to the posterior given the GMM estimator $\hat \theta$. One can further explore the higher order asymptotics of $q(\theta|\D)$ and $p(\theta|\hat \theta)$, for example expanding both posterior densities as Edgeworth series of the asymptotic pivotal quantity $\sqrt{n}\F(\theta)(\theta-\hat \theta)$. In this sense, our result in Part (ii) of Theorem \ref{qphat} only captures the leading order closeness from $q(\theta|\D)$ to $p(\theta|\hat \theta)$. However, we conjecture that in general the higher order terms of $q(\theta|\D)$ and $p(\theta|\hat \theta)$ do not match with each other, since the LIL takes a quadratic form of the moment conditions while the true density of $\hat \theta$ depends on other features of $P_{\D}$, such as the high order moments. Similar work in this direction includes \citet{fang06}, where they have shown by a simple example of sample mean that the Edgeworth expansions from the empirical likelihood and the density of the sample average do not agree in high order terms.
\end{remark}

\section{Numerical Study}
\subsection{Algorithm}
Because the LIL \eqref{qlik} allows any form of moment function $g(D,\theta)$, usually one cannot derive an analytical close form for the BGMM model posterior $q(\M|\D)$. Therefore, we adopt a reversible jump MCMC algorithm with Metropolis moves both between models and within a model to explore the joint posterior of $q(\theta,\M|\D)$, similar in spirit to the MCMC algorithm for the Gibbs posterior model selection (\citealt{chen10}, \citealt{jt08}), and also the PAC-Bayesian model selection (\citealt{alquier13}, \citealt{guedj13}). In the $i$th iteration, the between-model steps either add a new component to the nonzero part of $\theta^{(i)}$, or remove an existing component in the nonzero part of $\theta^{(i)}$, each with probability 0.5. When we add a new component, the parameter value for this new component is sampled from $\mathcal{N}(0,\sigma_{\add}^2)$, while the values of the existing components in $\theta^{(i)}$ are retained. Both the ``add" and the ``remove'' operations will be accepted or rejected with a probability based on the ratio of the posteriors evaluated at the new proposed parameter and the current parameter. This between-model step is then followed by a within-model step, in which we draw a new parameter value in the same model as $\theta^{(i)}$ from a proposal distribution. In practice, to efficiently explore each model space, we use a normal distribution as a proposal distribution, with mean zero and a properly chosen variance $c\cdot \bxi_{\M}$. Here $\bxi_{\M}$ is the submatrix of $\bxi$ with rows and columns corresponding to the model $\M$, and $\bxi$ is an estimated covariance matrix for the GMM estimator $\hat \theta$, which can be obtained numerically by inverting the Hessian matrix at the preliminary one-step GMM estimator $\tilde \theta$ on the full model space. We set $c=2.4^2$ as suggested in \citet{gelman} to achieve the ideal acceptance rate for within-model Metropolis moves. We also run pilot chains to tune the value of $\sigma_{\add}$ for better mixing of the Markov chain. As a result, the Markov chain consists of $\theta^{(i)}$ drawn from the full BGMM posterior across different model spaces.


\subsection{Example: Correlated Binary Responses}
The conditional mean $\mu_{ij}(\theta)=\Ep (Y_{ij}|X_{ij})$ of the longitudinal binary response $Y_{ij}$ is given by
\begin{equation}\label{logitmodel}
\ln \frac{\mu_{ij}(\theta)}{1-\mu_{ij}(\theta)} = X_{ij}^\top\theta,
\end{equation}
where $i=1,\ldots,n$ and $j=1,\ldots,s$. In the following simulations, we first fix the sample size $n=400$ and the cluster size $s=10$, in order to compare with the similar simulation setups in \citet{wang12}. For $X_{ij}=(X_{ij1},\ldots, X_{ijp})^\top$, we consider two situations with $p=50$ and $p=100$. $X_{ij1},\ldots,X_{ijp}$ are generated independently from a uniform distribution on $[-1,1]$. We also consider two sets of true parameter values,
\begin{align*}
\theta_0 &=(1.5,-1.5,1,-1,0.5,-0.5,0,0,\ldots,0)\\
\theta_0 &=(1.5,-1.5,1.5,-1.5,1,-1,1,-1,0.5,-0.5,0.5,-0.5,0,0,\ldots,0)
\end{align*}
with the number of nonzero components $k_0=6$ and $k_0=12$ respectively. Note that $\theta_0$ contains weak signals $0.5$ and $-0.5$ and more nonzero components in the second setting. Similar to \citet{wang12} and \citet{cho13}, we use the R package {\tt mvtBinaryEP} to generate the correlated binary responses $(Y_{i1},\ldots,Y_{is})^\top$ for each $i=1,\ldots,n$ with an exchangeable correlation structure with correlation coefficient $\rho=0.3$. \\

Since this is a special case of the first motivating example in Section \ref{sec.examples}, we examine the performance of BGMM using the moment function $g(D,\theta)$ defined in \eqref{cldmom}. We compare the BGMM method to the frequentist penalized GEE method (PGEE) proposed by \citet{wang12} which is used to fit high dimensional longitudinal data. Let $\theta_{(k)}$ be the $k$th component of $\theta$. The PGEE solves a similar estimating equation to \eqref{glm1}
$$n^{-1}\sum_{i=1}^n \frac{\partial \mu_i(\theta)^\top}{\partial \theta} \S_i^{-1}(Y_i- \mu_i(\theta)) - P_{\lambda_n}(\theta)=0,$$
with an additional SCAD penalty $P_{\lambda}(\theta) = (P_{\lambda}(\theta_{(1)}),\ldots,P_{\lambda}(\theta_{(p)}))^\top$ and for $k=1,\ldots,p$,
$$P_{\lambda}(\theta_{(k)}) = \lambda_n\left\{1(\theta_{(k)}\leq \lambda_n) +1\left(\lambda_n<\theta_{(k)}\leq a\lambda_n\right)\frac{a\lambda_n-\theta_{(k)}}{(a-1)\lambda_n}\right\}.$$
The PGEE can be solved by an iterative Newton-Raphson algorithm as described in \citet{wang12}. In our simulations, we perform in the same way as \citet{cho13}, fix $a=3.7$ and truncate the estimated coefficients to zero if $|\hat \theta_{(k)}|\leq 10^{-3}$ ($k=1,\ldots,p$). $\lambda_n$ is selected from the grid set $\{0.01,0.02,\ldots,0.2\}$ by 5-fold cross validation. We use an estimated correlation matrix for $\R$ based on the sample, instead of varying the correlation structures in \citet{wang12}. In fact, the finite sample estimates of $\R$ are quite precise for the true $\R$ in our $p<n$ case.\\

For the BGMM method, the prior on $\theta$ given a model $\M$ is the product of independent normal densities
\begin{equation}\label{pit1}
\pi(\theta|\M)=\prod_{j\in \M}\frac{1}{\sqrt{2\pi}\sigma_{\theta}}e^{-\frac{\theta_{(j)}^2}{2\sigma^2_{\theta}}},
\end{equation}
where we choose $\sigma_{\theta}=10$ for a large prior spread. Note that although theoretically this prior is not truncated on a compact set $\Theta$ as in Assumption 7, in practice this has no influence on in our experimental results.\\

The prior on the model $\M$ is specified as follows:
\begin{equation}\label{pim1}
\pi(\M)\propto \nu^{|\M|}(1-\nu)^{p-|\M|},
\end{equation}
which means that each component of $\theta$ independently enters the model $\M$ with probability $\nu\in (0,1)$. When $\nu=0.5$, this is the same as the uniform prior over all $2^p$ models. When $\nu$ moves towards zero, the prior gradually induces more sparsity on $\theta$ and favors more parsimonious models, which imposes a further penalization on the model size besides the incorporated BIC-type penalization in BGMM. It can be verified (see Section 4 of the supplementary material) that the prior \eqref{pim1} satisfies Assumption 8 when $\nu\in (0,1)$ is either fixed or $\nu=n^{-c}$ for some $c>0$, and it satisfies Assumption 9(iv) if $\nu$ is fixed. \\

In our simulation, for each simulated datasets, we run one single Markov chain with the length $3\times 10^4$, and drop the first $10^4$ iterations as burnin. We consider two choices of the tuning parameter $\nu$ in \eqref{pim1}. In the first case (referred to as BGMM1), we fix $\nu=0.5$ throughout the Markov chain for the next $2\times 10^4$ iterations. In the second case (referred to as BGMM2), we adopt a two-step tuning strategy in an effort to make the value of $\nu$ more adaptive to the sparsity level of the true model. We estimate the posterior average model size $\hat \Ep_{\M|\D}|\M|$ using the first $10^4$ MCMC runs after the burnin, and then reset $\nu=\hat \Ep_{\M|\D}|\M|/p$ for the next $10^4$ MCMC runs. Finally for both chains of BGMM1 and BGMM2, we keep $N=10^3$ MCMC samples from the last $10^4$ runs for every 10 iteration. The variance of proposal normal density described in Section 3.1 is fixed at $\sigma_{\add}=0.2$. Our experiments with other values of $\sigma_{\add}$ (such as $0.05,0.1,0.15,0.25$) show that $\sigma_{\add}=0.2$ is sufficient for exploring the full posterior of $\theta$, and the MCMC results such as the MAP models and posterior distributions of parameters are not sensitive to different values of $\sigma_{\add}$. \\

As a benchmark, the PGEE method and the BGMM method are compared together with the naive method and the oracle method. The naive method estimates $\theta$ by usual GEE without doing model selection, while for the oracle method, the true model is pretendedly known and $\theta$ is estimated only on the nonzero components. We apply each method to the same dataset and repeat this process for 100 Monte Carlo replications. We compare three aspects of these methods: the model selection, the parameter estimation, and the prediction.\\

To evaluate the model selection performance, we consider the model selected by PGEE and the MAP model from BGMM, and report the proportion of times the method exact selecting (EX), underselecting (UN) and overselecting (OV)
the nonzero components of $\theta_0$. We also report the true positives (TP, the average number of correctly selected nonzero components in $\theta_0$), and the false positives (FP, the average number of selected nonzero components that are actually zero in $\theta_0$).\\

For the estimation accuracy, similar to \citet{cho13}, we report the estimated mean square error (MSE) $\sum_{m=1}^{100} \|\hat \theta_m-\theta_0\|^2/(100k_0)$, where $\hat \theta_m$ is the $m$th estimated parameter vector. This MSE is calculated for the naive method, the oracle method, the PGEE method, and the posterior mean of $\theta$ from the BGMM method. \\

For the prediction accuracy, we calculate the average MSE for the conditional mean $\mu_{ij}$ (denoted by pMSE), defined as $\sum_{i=1}^n\sum_{j=1}^s \left(\mu_{ij}(\hat \theta)-\mu_{ij}(\theta_0)\right)^2/(ns)$ for the naive, the oracle, and the PGEE method. For the BGMM method, we use the pMSE averaged over the posterior sample $\sum_{i=1}^n\sum_{j=1}^s \sum_{k=1}^N \left(\mu_{ij}( \theta^{(k)})-\mu_{ij}(\theta_0)\right)^2/(Nns)$, where $\theta^{(1)},\ldots, \theta^{(N)}$ are the MCMC draws of $\theta$. \\

As Table \ref{logittab} indicates, both the frequentist PGEE method and our BGMM method have always successfully identified the nonzero components of $\theta_0$ with no underselection. However, the PGEE performs much more conservative and has a serious overselection problem in all the simulations settings, which is consistent with the findings in \citet{cho13}. It selects the true model for $33\%$ of all time when $p=50$, and only $16\%$ of all time when $p=100$ and $k_0=12$. Meanwhile PGEE overselects about $4\sim 6$ extra redundant variables on average. In contrast, the BGMM MAP models have much higher probability of exactly selecting the true model, and have much smaller false positives. We also note that the extra two-step tuning of $\nu$ in BGMM2 has brought significant advantage over the uniform model prior with fixed $\nu=0.5$ in BGMM1. When $p=100$, the performance of BGMM1 deteriorates as the probability of exact selection drops to about $50\%$, but BGMM2 still maintains a high accuracy with over $90\%$ of exact model selection. This is because that in BGMM2, the first step of $10^4$ runs has consistently estimated the true proportion of nonzero components in $\theta$, and then the second step of $10^4$ runs can learn the sparsity of the model space better with $\nu$ roughly equal to the true average marginal inclusion probability.\\

For the estimation and prediction, it is clear that the naive GEE estimator with no model selection performs poorly in MSE and pMSE compared to the oracle estimator. Figure 1 and Figure 2 show that the MSE and pMSE for the BGMM method are comparable to those from the oracle and the PGEE method, as their boxplots largely overlap with each other. Also it seems that BGMM tends to have smaller variation across difference simulations than PGEE. The averaged levels of three MSEs from both BGMM methods are also slightly smaller than those from the PGEE estimator in most of the cases (Table \ref{logittab}), and they are all close to the MSE and pMSE from the oracle estimator. Overall, BGMM2 seems to be the best of all these methods besides the oracle. This has partly supported our theoretical results about the oracle properties of the BGMM method, in the sense that the posterior variance of BGMM is asymptotically the same as the variance of the oracle GMM estimator.\\

Finally, we vary the sample size $n$ among $200,400,800,1200,2000$ and compare the performance of PGEE, BGMM1 and BGMM2 for the model \eqref{logitmodel} averaged over 20 simulated datasets. Figure 3 and Figure 4 plot their exact model selection probabilities (the same as the EX in Table \ref{logittab}), MSEs and pMSEs. Overall, BGMM2 has the best performance of exact model selection, and BGMM1 tends to perform better as $n$ increases. All three methods have poor model selection accuracy for $p=100,n=200$ due to the relative high dimension and the small sample size. As the sample size $n$ increases, the differences between their MSEs and pMSEs become negligible, as they all perform similarly to the oracle estimator.

\newpage
\begin{table}[H]
\centering
\caption{Comparison of BGMM with PGEE for Correlated Binary Responses. $k_0$ is the number of nonzero components in the true parameter $\theta_0$. $p$ is the dimension of $\theta_0$. $n$ is the sample size. Standard errors are shown in the parentheses. EX: exact selection; UN: under selection; OV: over selection; TP: true positives; FP: False positives. MSE: mean square error of $\theta$; pMSE: prediction mean square error of $\mu_{ij}(\theta)$.}
\small
\label{logittab}
\begin{tabular}{cccccccc}
\hline
   \hline
   \multicolumn{8}{c}{$k_0=6$, $p=50$, $n=400$}\\
\hline
 & EX & UN & OV  & TP & FP & MSE ($\times 10^{-3}$) & pMSE ($\times 10^{-4}$)\\
\hline
Naive  & 0 & 0 & 1 & 6 & 44 & 24.81 (0.57) & 21.77 (0.43) \\
Oracle & 1 & 0 & 0 & 6 & 0 & 4.15 (0.35) & 2.66 (0.16)\\
PGEE  & 0.33 & 0 & 0.67 & 6 & 5.18 & 7.02 (0.57) &  5.05 (0.41) \\
BGMM1 & 0.80 & 0 & 0.20 & 6 & 0.21 & 5.58 (0.40) & 3.85 (0.23)\\
BGMM2 & 0.98 & 0 & 0.02 & 6 & 0.02 & 5.06 (0.38) & 3.33 (0.20) \\
\hline
   \hline
 \multicolumn{8}{c}{$k_0=12$, $p=50$, $n=400$}\\
\hline
 & EX & UN & OV & TP & FP & MSE ($\times 10^{-3}$) & pMSE ($\times 10^{-4}$)\\
\hline
Naive  & 0 & 0 & 1 & 12 & 38 & 13.09 (0.35) & 20.33 (0.40) \\
Oracle & 1 & 0 & 0 & 12 & 0 & 5.53 (0.49) & 4.82 (0.19)\\
PGEE  & 0.33 & 0 & 0.67 & 12 & 4.08 & 6.13 (0.58) &  6.61 (0.38) \\
BGMM1 & 0.88 & 0 & 0.12 & 12 & 0.13 & 5.13 (0.29) & 5.91 (0.23)  \\
BGMM2 & 0.94 & 0 & 0.06 & 12 & 0.07 & 4.89 (0.27) & 5.63 (0.22)  \\
\hline
   \hline
     \multicolumn{8}{c}{$k_0=6$, $p=100$, $n=400$}\\
\hline
 & EX & UN & OV  & TP & FP & MSE ($\times 10^{-3}$) & pMSE ($\times 10^{-4}$)   \\
\hline
Naive  & 0 & 0 & 1 & 6 & 94 & 49.28 (0.67) & 42.53 (0.61) \\
Oracle & 1 & 0 & 0 & 6 & 0 & 6.72 (0.68) & 3.09 (0.17) \\
PGEE  & 0.28 & 0 & 0.72 & 6 & 5.27 & 9.16 (1.37) &  5.21 (0.51) \\
BGMM1 & 0.55 & 0 & 0.45 & 6 & 0.61 & 8.39 (0.51) & 5.32 (0.27) \\
BGMM2 & 0.93 & 0 & 0.07 & 6 & 0.07 & 7.32 (0.51) & 4.05 (0.21) \\
\hline
   \hline
        \multicolumn{8}{c}{$k_0=12$, $p=100$, $n=400$}\\
\hline
 & EX & UN & OV  & TP & FP & MSE ($\times 10^{-3}$) & pMSE ($\times 10^{-4}$)  \\
\hline
Naive  & 0 & 0 & 1 & 12 & 88 & 27.39 (0.55) & 42.74 (0.63) \\
Oracle & 1 & 0 & 0 & 12 & 0 & 8.09 (0.64) & 5.85 (0.22) \\
PGEE  & 0.16 & 0 & 0.84 & 12 & 6.35 & 10.20 (1.40) &  8.62 (0.63) \\
BGMM1 & 0.59 & 0 & 0.41 & 12 & 0.55 & 8.24 (0.52) & 8.82 (0.36) \\
BGMM2 & 0.90 & 0 & 0.10 & 12 & 0.10 & 7.76 (0.52) & 7.68 (0.29) \\
\hline
   \hline
\end{tabular}
\end{table}

\begin{figure}[H]
  \centering
\includegraphics[width=11.5cm,height=9.75cm]{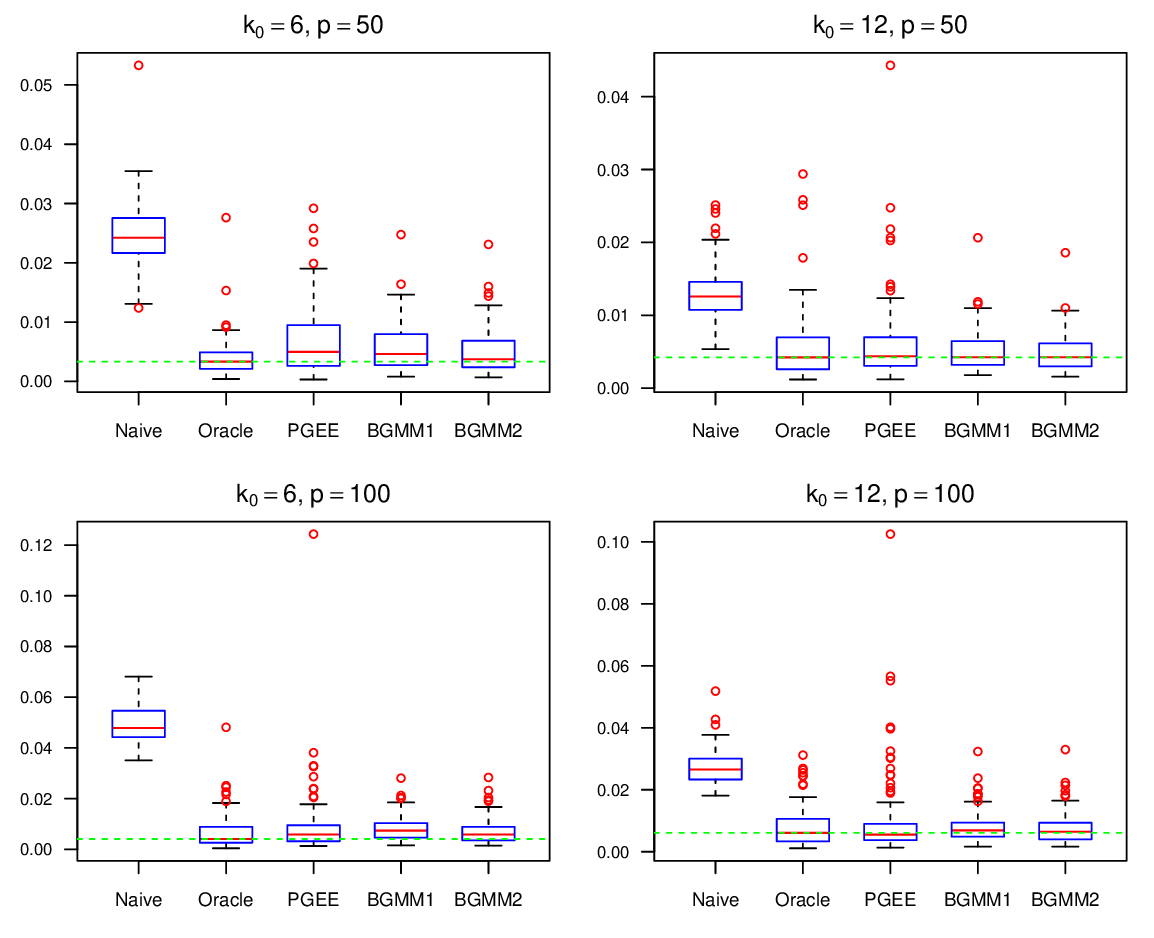}
  \caption{Boxplots for the MSE of $\theta$ over 100 simulated datasets. }\label{fig1}
\end{figure}

\begin{figure}[H]
  \centering
\includegraphics[width=11.5cm,height=9.75cm]{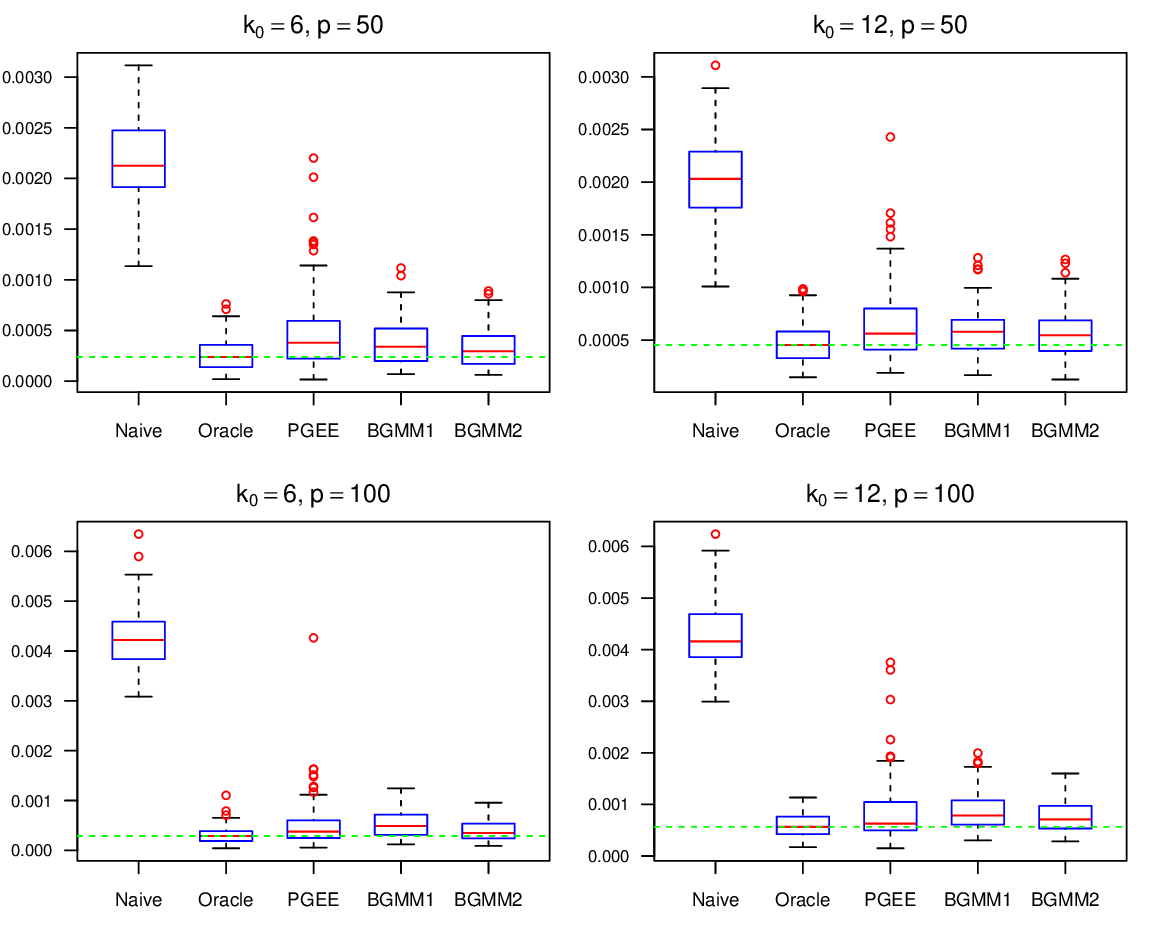}
  \caption{Boxplots for the MSE of $\mu_{ij}(\theta)$ over 100 simulated datasets. }\label{fig2}
\end{figure}

\begin{figure}[H]
  \centering
\includegraphics[width=12.5cm,height=9.75cm]{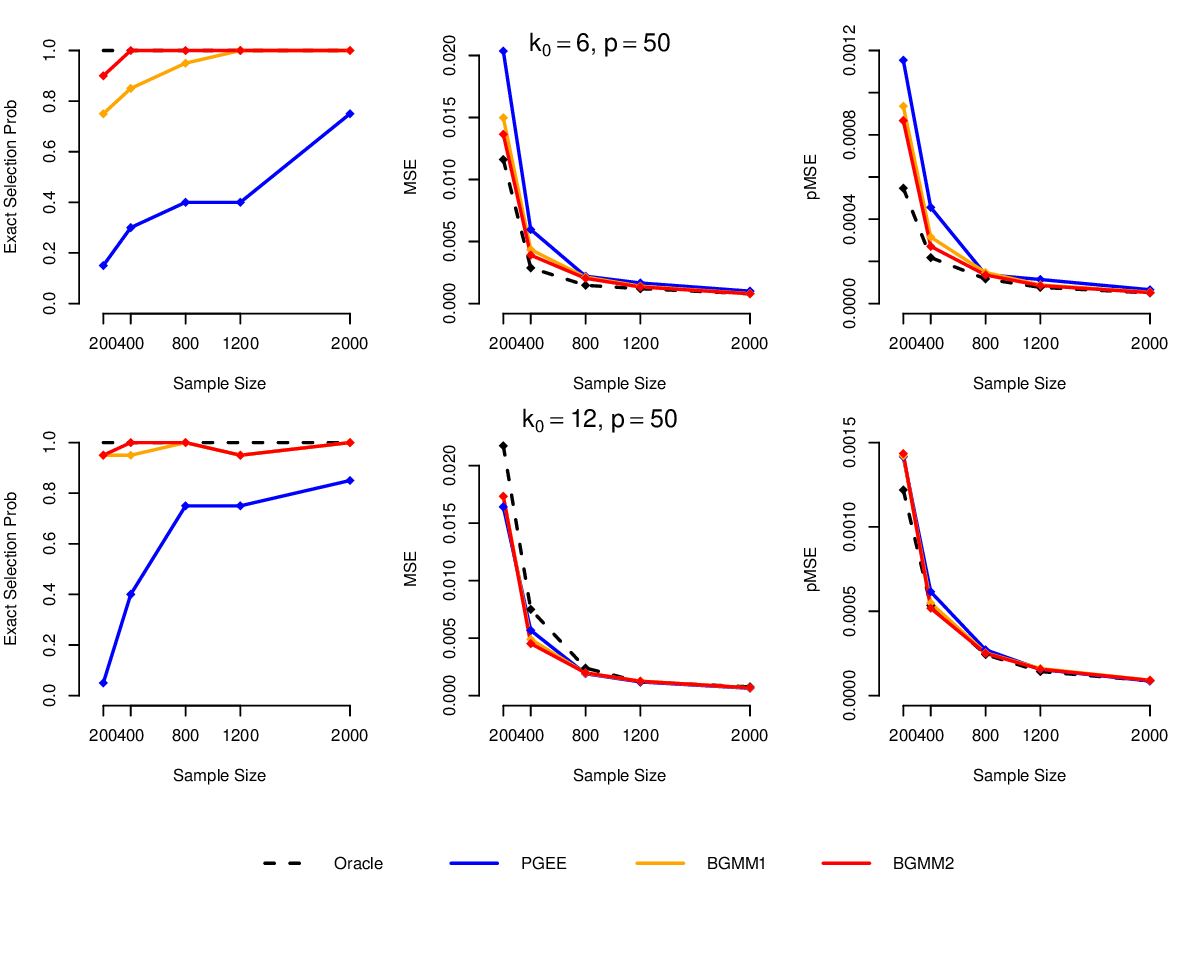}
  \caption{Exact selection probability, MSE, and prediction MSE for $p=50$ over 20 simulated datasets. }\label{fig3}
\end{figure}

\begin{figure}[H]
  \centering
\includegraphics[width=12.5cm,height=9.75cm]{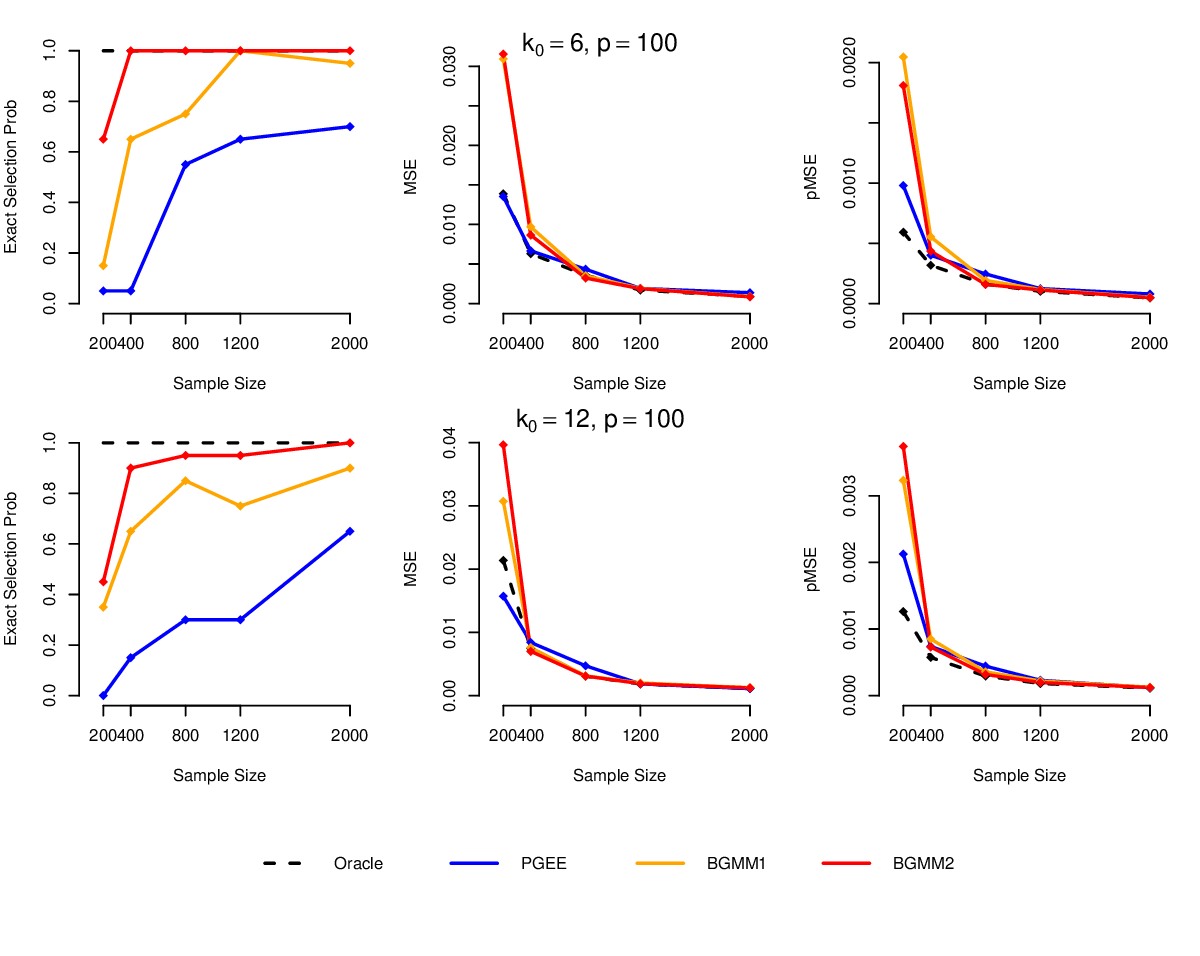}
  \caption{Exact selection probability, MSE, and prediction MSE for $p=100$ over 20 simulated datasets. }\label{fig4}
\end{figure}

\section{Discussions}
In this paper, we have studied some theoretical properties and applications of a Bayesian moment based model selection method. As we have commented, this method combines advantages of a Bayesian approach, such as the expressiveness of the posterior distribution and convenient MCMC algorithms for computation, with the model robustness of the moment based methods. We have formulated and proved the Bayesian oracle property of the proposed model selection method, which guarantees efficient posterior inference as if we knew which variables are truly relevant. We have studied the meaning of the quasi-posterior probabilities used in BGMM, which can be interpreted as the leading order large sample approximation to the true posterior probabilities conditional on the observed GMM estimator. The empirical performance of BGMM has been demonstrated by numerical experiments.\\

We have only considered quasi-posterior constructed from the GMM based quasi-likelihood function. Many other alternatives, such as EL, GEL, and ETEL, can be formulated under a similar Bayesian framework, with possible interpretations of the induced quasi-Bayesian posterior. See for example, \citet{ch03}, \citet{lazar03}, \citet{schen05}, etc. We conjecture that similar Bayesian asymptotic properties for model selection can be derived for these quasi-likelihoods. \\

{\bf Acknowledgment}: We thank the Associate Editor and the two anonymous Referees from Journal of Multivariate Analysis for their helpful suggestions on improving the paper.\\

\newpage

\begin{center}
{\bf \Large Supplementary Materials}
\end{center}

\vspace{1cm}

\setcounter{equation}{0}
\setcounter{lemma}{0}
\renewcommand{\theequation}{A.\arabic{equation}}
\renewcommand{\thelemma}{A.\arabic{lemma}}

This document consists of four parts. Part 1 includes the proof Theorem 1 in the main paper. Part 2 includes the proofs of Theorem 2 and 3 in the main paper. In Part 3 we prove the oracle property for the three motivating examples mentioned in Section 1.3 of the main paper. In Part 4 we provide a proposition that verifies Assumption 8 for several examples of priors on the models.

 \vspace{+1cm} \noindent {\bf
\Large 1. Proofs for the Bayesian Oracle Property of BGMM}

\vspace{+0.2cm}

We first prove some useful lemmas. For a generic square matrix $\C$, let $\tr(\C)$ be the trace of $\C$, and $\|\C\|_F=\sqrt{\tr (\C^\top\C)}$ be the Frobenius norm of $\C$. In the following, we use the statement ``the event $A$ happens w.p.$1-\eta$" to denote the relation $P_{\D}(A)\geq 1-\eta$.

\begin{lemma} \label{wn}(\citealt{bc09})
Let
$$W_n(\D,\theta)=\bar g(\D,\theta)-\Ep g(D,\theta)-\left(\bar g(\D,\theta_0)-\Ep g(D,\theta_0)\right).$$
Then under Assumptions 2 and 4, uniformly for all $\theta \in \Theta$,
$$\|W_n(\D,\theta)\|=O_p\left(\sqrt{\frac{p\ln n}{n}}\|\theta-\theta_0\|^\alpha + n^{-1}p^{3/2}\ln n\right).$$
\end{lemma}
\noindent \textbf{Proof:}
The proof can be found in (A.10) of \citet{bc09}. The $L_2$ norm of the deviation $W_n(\D,\theta)$ can be controlled using an empirical process result when the moment $g(D,\theta)$ satisfies Assumption 4. The only adaptation here is that the condition on VC dimension in ZE.1 of their paper is now replaced by the condition on the uniform covering number in Assumption 4(ii). The conclusion still holds according to the proof of Lemma 16 in \citet{bci11}.\hfill$\blacksquare$

\vspace{+0.5cm}

\begin{lemma}\label{g11}

Given a model $\M$, let $\theta=(\theta_1^\top,\theta_2^\top)^\top$ be decomposed according to $\M$, and let $\G_{\M}$ be the same as defined in Theorem 1. Define
\begin{align*}
S_{\M}(\D) &=\exp\left\{-\frac{n}{2}\bar g(\D,\theta_0)^\top \left(\V_n^{-1}-\V_n^{-1}\G_{\M}(\G_{\M}^\top\V_n^{-1}\G_{\M})^{-1}\G_{\M}^\top\V_n^{-1}\right) \bar g(\D,\theta_0)\right\},\\
\bar \theta_{\M,1} &= \theta_{0,\M,1}-(\G_{\M}^\top\V_n^{-1}\G_{\M})^{-1}\G_{\M}^\top \V_n^{-1}\bar g(\D,\theta_0),
\end{align*}

where $\theta_{0,\M,1} \in \mathbb{R}^{|\M|}$ be the subvector of $\theta_0$ restricted to $\Theta(\M)$.  Then under Assumptions 1-8, uniformly for all spaces $\M \supseteq \M_0$, for any fixed constant $C>0$,
\begin{align*}
& \int_{B_0(C\e)\cap \Theta(\M)}e^{-\frac{n}{2}\bar g(\D,\theta)^\top \V_n^{-1} \bar g(\D,\theta)} \pi(\theta|\M) \ud\theta \\
&= (1+o_p(1))S_{\M}(\D)\int_{B_0(C\e)\cap \Theta(\M)}e^{-\frac{n}{2}(\theta_1-\bar \theta_{\M,1})^\top \G_{\M}^\top\V_n^{-1}\G_{\M}(\theta_1-\bar \theta_{\M,1})}\pi(\theta|\M)\ud \theta_1.
\end{align*}
\end{lemma}

\noindent \textbf{Proof:}
First of all, the $L_2$ norm of $\bar g(\D,\theta_0)$ satisfies for any $C>0$,
\begin{equation}\label{g0}
\Pr\Bigg(\|\bar g(\D,\theta_0)\|\geq C\sqrt{\frac{p}{n}}\Bigg)\leq \frac{n\Ep \|\bar g(\D,\theta_0)\|^2}{C^2p}=\frac{\tr(\var (g(D,\theta_0)))}{C^2p}\leq \frac{\bar \lambda(\V)}{C^2},
\end{equation}
which implies that $\|\bar g(\D,\theta)\|=O_p(\sqrt{p/n})$ since the eigenvalues of $\V=\var\left(g(D,\theta_0)\right)$ are bounded above according to Assumption 6.

Second, for $\theta\in B_0(C\e)\cap \Theta(\M)$, let $r_{\M}(\D,\theta)=\bar g(\D,\theta)- \bar g(\D,\theta_0) - \G_{\M}(\theta_1-\theta_{0,\M,1})$. Then using second order Taylor expansion of $\Ep g(D,\theta)$ at $\theta_0$ for $\theta \in B_0(C\e)\cap \Theta(\M)$ and with all zero components of $\theta$ excluded, we have
$$r_{\M}(\D,\theta)=\frac{1}{2}\H_{\M}(\tilde \theta_1)(\theta_1-\theta_{0,\M,1},\theta_1-\theta_{0,\M,1})+W_n(\D,\theta),$$
where $\tilde \theta_1$ is between $\theta_1$ and $\theta_{0,\M,1}$ and $\tilde \theta=(\tilde \theta_1^\top,0)^\top$, $\H_{\M}$ is the submatrix of the second order derivative matrix $\H$ restricted to $\Theta(\M)$. By Assumption 5(iii), we have that
$$ \left\|\H_{\M}(\tilde \theta_1)(\theta_1-\theta_{0,\M,1},\theta_1-\theta_{0,\M,1})\right\|
\leq  \left\|\sup_{\|u\|=1,\|v\|=1}\H(\tilde \theta)(u,v)\right\|\left\|\theta_1-\theta_{0,\M,1}\right\|^2
\leq  O(\sqrt{p}\e^2).$$
Therefore, using Lemma \ref{wn}, we obtain that the order of $r_{\M}(\D,\theta)$ on $\theta\in B_0(C\e)\cap \Theta(\M)$ uniformly for all $\M \supseteq \M_0$ is,
$$\|r_{\M}(\D,\theta)\| \leq O_p\left(n^{-1}p^{3/2}+ (p/n)^{\frac{\alpha+1}{2}}\sqrt{\ln n} + n^{-1}p^{3/2}\ln n\right)= o_p\big((pn)^{-1/2}\big),$$
where the last equality also holds if $p$ increases as in the footnote of Assumption 2. Then using decomposition
$\bar g(\D,\theta)=\G_{\M}(\theta_1-\theta_{0,\M,1})+\bar g(\D,\theta_0)+ r_{\M}(\D,\theta)$ we have
\begin{align}\label{gdecomp}
&  \frac{n}{2}\bar g(\D,\theta)^\top \V_n^{-1} \bar g(\D,\theta)\nonumber \\
={} & \frac{n}{2}\Bigg\{ (\theta_1-\bar \theta_{\M,1})^\top (\G_{\M}^\top \V_n^{-1}\G_{\M} )(\theta_1-\bar \theta_{\M,1})\Bigg\}\nonumber \\
+{} &\frac{n}{2}\Bigg\{ \bar g(\D,\theta_0)^\top \left(\V_n^{-1}-\V_n^{-1}\G_{\M}(\G_{\M}^\top\V_n^{-1}\G_{\M})^{-1}\G_{\M}^\top\V_n^{-1}\right) \bar g(\D,\theta_0)\Bigg\}\nonumber \\
+{} & \frac{n}{2}\Bigg\{r_{\M}(\D,\theta)^\top \V_n^{-1}r_{\M}(\D,\theta)+ 2r_{\M}(\D,\theta)^\top \V_n^{-1} \G_{\M}(\theta_1-\theta_{0,\M,1}) + 2r_{\M}(\D,\theta)^\top\V_n^{-1} \bar g(\D,\theta_0)\Bigg\}\nonumber \\
\end{align}
where $\bar \theta$ is defined in the lemma. By Assumptions 5(ii) and 6, the eigenvalues of $\V_n$ and $\G^\top\G$ are bounded above and below w.p.a.1 as $n\to \infty$, so are the eigenvalues of any $\G_{\M}^\top\G_{\M}$ since $\G_{\M}$ is a submatrix of $\G$. Therefore on $B_0(C\e)\cap \Theta(\M)$, $\|r_{\M}(\D,\theta)^\top \V_n^{-1}r_{\M}(\D,\theta)\|\leq \underline \lambda(\V_n)^{-1}\|r_{\M}(\D,\theta)\|^2=o_p((pn)^{-1})$, $\|2r_{\M}(\D,\theta)^\top \V_n^{-1} \G_{\M}(\theta_1-\theta_{0,\M,1})\|\leq 2\underline \lambda(\V_n)^{-1}\bar \lambda (\G_{\M}^\top\G_{\M})\cdot \|r_{\M}(\D,\theta)\|\cdot \|\theta-\theta_0\|\leq o_p((pn)^{-1/2}\e)=o_p(n^{-1})$, $\|2r_{\M}(\D,\theta)^\top\V_n^{-1} \bar g(\D,\theta_0)\|\leq 2\underline \lambda(\V_n)^{-1}\cdot\|r_{\M}(\D,\theta)\| \cdot\|\bar g(\D,\theta_0)\|=o_p((pn)^{-1/2}(p/n)^{1/2})=o_p(n^{-1})$. These together imply that the last term in \eqref{gdecomp} is of order $o_p(1)$, and this holds uniformly for all $\M \supseteq \M_0$. The conclusion then follows if we define $S_{\M}(\D)$ as in the lemma, which does not depend on $\theta$ and can be moved outside the integral. \hfill$\blacksquare$

\vspace{+0.5cm}

\begin{lemma}\label{g12}
Under Assumptions 1-8, there exists a constant $C_1>0$, such that uniformly for all spaces $\M \supseteq \M_0$, for any fixed constant $C\geq C_1$ and all sufficiently large $n$,
\begin{align}\label{laplace}
& \int_{B_0(C\e)\cap \Theta(\M)}e^{-\frac{n}{2}(\theta_1-\bar \theta_{\M,1})^\top \G_{\M}^\top\V_n^{-1}\G_{\M}(\theta_1-\bar \theta_{\M,1})}\pi(\theta|\M)\ud \theta_1 \nonumber \\
 ={}& (2\pi/n)^{|\M|/2}\big\{\dett (\G_{\M}^\top\V_n^{-1}\G_{\M})\big\}^{-1/2}\pi(\theta_0|\M)\big(1+o_p(1)\big),
\end{align}
where $\bar \theta_{\M,1}$ is defined in Lemma \ref{g11} and $\theta=(\theta_1^\top,\theta_2^\top)^\top$ is decomposed according to $\M$.
\end{lemma}

\noindent \textbf{Proof:}
First we let $\P_{\M}=\V_n^{-1/2}\G_{\M}(\G_{\M}^\top\V_n^{-1}\G_{\M})^{-1}\G_{\M}^\top \V_n^{-1/2}$, where
$\V_n^{1/2}$ is the symmetric positive definite square root of $\V_n$. Then $\P_{\M}$ is idempotent and has eigenvalues 0 and 1. The difference between $\bar \theta_{\M,1}$ and $\theta_{0,\M,1}$ can be controlled by
\begin{align*}
&  \|\bar \theta_{\M,1}-\theta_{0,\M,1}\|^2=\|(\G_{\M}^\top\V_n^{-1}\G_{\M})^{-1}\G_{\M}^\top \V_n^{-1}\bar g(\D,\theta_0)\|^2\\
\leq {}&  \underline \lambda \big(\G_{\M}^\top\V_n^{-1}\G_{\M} \big)^{-1} \cdot \bar g(\D,\theta_0)^\top \V_n^{-1/2} \P_{\M} \V_n^{-1/2}  \bar g(\D,\theta_0)\\
\leq {}& \bar \lambda (\V_n) \underline \lambda (\G^\top\G)^{-1} \cdot \underline \lambda (\V_n)^{-1} \|\bar g(\D,\theta_0)\|^2.
\end{align*}

Since $\|\bar g(\D,\theta_0)\|=O_p(\sqrt{p/n})$, we know that $ \|\bar \theta_{\M,1}-\theta_{0,\M,1}\|$ is also $O_p(\sqrt{p/n})$ since all the eigenvalues here are bounded. So for any small $\eta>0$, we can pick $C'$ sufficiently large, such that $ \|\bar \theta_{\M,1}-\theta_{0,\M,1}\|\leq C'\sqrt{p/n}$ w.p.$1-\eta$ and uniformly for all $\M \supseteq \M_0$.

Next we evaluate the integral on the left hand side of \eqref{laplace} on $B_{\M}(C\e):=\{\theta=(\theta_1^\top,0)^\top\in \Theta(\M):\|\theta_1-\bar \theta_{\M,1}\|\leq C\e\}$ for a fixed $C>0$. We observe that the integral takes the same form as a Gaussian random vector centered at $\bar \theta_{\M,1}$. Define $U\sim \mathcal{N}(0,\G_{\M}^\top\V_n^{-1}\G_{\M})$. Then we have that there exists a large $C''$, such that when $C\geq C''$, w.p.$1-2\eta$ and uniformly for all $\M \supseteq \M_0$,
\begin{align}\label{laphat}
&  \int_{B_{\M}(C\e)}e^{-\frac{n}{2}(\theta_1-\bar \theta_{\M,1})^\top \G_{\M}^\top\V_n^{-1}\G_{\M}(\theta_1-\bar \theta_{\M,1})}\pi(\theta|\M)\ud \theta_1 \nonumber \\
 ={}& (2\pi/n)^{|\M|/2}\big\{\dett (\G_{\M}^\top\V_n^{-1}\G_{\M})\big\}^{-1/2} \Pr\big(\|U\|\leq C\sqrt{p}\big)\cdot \pi(\theta_0|\M)(1+o(1)) \nonumber \\
= {}& (2\pi/n)^{|\M|/2}\big\{\dett (\G_{\M}^\top\V_n^{-1}\G_{\M})\big\}^{-1/2} \pi(\theta_0|\M)(1+o(1)),
\end{align}
where the $o(1)$ depends on $\eta$ (and hence on $C$ and $n$). In the first equality above, we have used Assumption 7(ii) to obtain that $\pi(\theta|\M) = \pi(\theta_0|\M) (1+o(1))$ as $n\to \infty$ uniformly over all $\M \supseteq \M_0$ and all $\theta \in B_0((C+C')\e)$, since w.p.$1-\eta$, $B_{\M}(C\e)\subseteq B_0((C+C')\e)$. In the second equality, we used the fact that the eigenvalues of $\G_{\M}^\top\V_n^{-1}\G_{\M}$ are bounded in probability using Assumptions 5 and 6. Hence by Chebyshev's inequality, $\|U\|=O_p(\sqrt{p})$. Hence we can pick a large $C''$ such that for $C\geq C''$, w.p.$1-\eta$, $\Pr(\|U\|\leq C\sqrt{p})=1+o(1)$.

Finally we set $C_1=C'+C''$. Then for $C\geq C_1$, $B_{\M}((C-C')\e)\subseteq B_0(C\e)\cap \Theta(\M) \subseteq B_{\M}((C+C')\e)$, and $C-C'\geq C''$ guarantees that \eqref{laphat} is satisfied w.p.$1-2\eta$ and uniformly for all $\M \supseteq \M_0$. Therefore \eqref{laplace} follows since the integral on the left hand side of \eqref{laplace} can be bounded between the integrals on $B_{\M}((C-C')\e)$ and $B_{\M}((C+C')\e)$, and both integrals satisfy \eqref{laphat}. \hfill$\blacksquare$

\vspace{+0.5cm}
\begin{lemma}\label{g13}
Under Assumptions 1-8, there exists a constant $C_2>0$, such that uniformly for all spaces $\M \supseteq \M_0$, for all large constant $C\geq C_2$, w.p.a.1 as $n\to\infty$,
\begin{align*}
& \int_{ \Theta(\M)\backslash B_0(C\e)}e^{-\frac{n}{2}\bar g(\D,\theta)^\top \V_n^{-1} \bar g(\D,\theta)} \pi(\theta|\M) \ud\theta\\
\leq {}&  c_{\pi} \left(\frac{4\pi \bar\lambda} {n\delta_1^2}\right)^{|\M|/2}\exp\left(-\frac{C^2\delta_1^2}{16\bar \lambda}p\right)+ \exp\left(-\frac{n}{4}\bar \lambda^{-1}\delta_0 ^2\right),
\end{align*}
where $c_{\pi}$ is from Assumption 7(i), and $\delta_0, \delta_1$ are from Assumption 5(i).
\end{lemma}
\noindent \textbf{Proof:} The proof uses similar techniques to the proof of Lemma 8 in \citet{bc09}.
Here we first directly cite part of the results from \citet{bc09}, since they remain valid under our Assumptions 2, 4, and 5.

\vspace{+0.5cm}

\noindent 1. For any small $\eta>0$, there exists a large $C'>0$, such that $\|\Ep g(D,\theta)\|> 8\|\bar g(\D,\theta_0)\|$ uniformly on $\Theta\backslash B_0(C'\e)$ w.p. $1-\eta$. $C'$ depends on $\delta_0,\delta_1$ in Assumption 5(i) and $\bar \lambda$ in Assumption 6.

\vspace{+0.5cm}

\noindent 2. $\|W_n(\D,\theta)\|=o_p(\|\Ep g(D,\theta)\|)$ uniformly on $\Theta\backslash B_0(C'\e)$. So for $n$ sufficiently large,
$\|W_n(\D,\theta)\|\leq \|\Ep g(D,\theta)\|/8$ for all $\theta \in \Theta\backslash B_0(C'\e)$ w.p. $1-\eta$.

\vspace{+0.5cm}

Note that the two results above hold uniformly for all $\theta \in \Theta(\M)\backslash B_0(C'\e)$ and for all $\M \supseteq \M_0$. Therefore, we have
\begin{align*}
\|\bar g(\D,\theta)\| &=  \left\|\Ep g(D,\theta)+\bar g(\D,\theta_0)+W_n(\D,\theta)\right\|\\
&\geq  \Big|  \|\Ep g(D,\theta)\|-\|\bar g(\D,\theta_0)\|-\|W_n(\D,\theta)\| \Big| \\
& \geq  \frac{3}{4}  \left\|\Ep g(D,\theta)\right\|
\end{align*}
uniformly for all $\theta \in \Theta(\M)\backslash B_0(C'\e)$, all $\M \supseteq \M_0$, and all sufficiently large $n$ w.p. $1-2\eta$.

Therefore, for $C>C'$, by Assumptions 1-8, we have
\begin{align}\label{out}
&  \int_{ \Theta(\M)\backslash B_0(C\e)}e^{-\frac{n}{2}\bar g(\D,\theta)^\top \V_n^{-1} \bar g(\D,\theta)} \pi(\theta|\M) \ud\theta \nonumber \\
\leq {}& \int_{ \Theta(\M)\backslash B_0(C\e)} \exp\left\{-\frac{n}{2}\bar \lambda (\V_n)^{-1}\cdot \frac{9}{16}\left\|\Ep g(D,\theta)\right \|^2  \right\}\pi(\theta|\M)\ud \theta\nonumber \\
\leq {}& c_{\pi}\int_{ \Theta(\M)\backslash B_0(C\e)}  \exp\left\{-\frac{n}{4}\bar \lambda ^{-1} \delta_1^2 \left\|\theta-\theta_0\right\|^2  \right\}\ud \theta \nonumber \\
& + \int_{ \Theta(\M)\backslash B_0(C\e)}
\exp\Big\{-\frac{n}{4}\bar \lambda ^{-1} \delta_0 ^2 \Big\}\pi(\theta|\M)\ud \theta \nonumber \\
\leq {}& c_{\pi}\left(\frac{4\pi\bar \lambda} {n\delta_1^2}\right)^{|\M|/2}\Pr(\|U\|\geq C\sqrt{p}) + \exp\left(-\frac{n}{4}\bar \lambda ^{-1} \delta_0 ^2\right) \nonumber \\
\leq {} &c_{\pi} \left(\frac{4\pi\bar \lambda }{n\delta_1^2}\right)^{|\M|/2}\exp\left(-\frac{C^2\delta_1^2}{16\bar \lambda }p\right)+ \exp\left(-\frac{n}{4}\bar \lambda^{-1}\delta_0 ^2\right)
\end{align}
where in the second inequality we used Assumption 5(i), 6 and 7(i) and required $n$ to be sufficiently large, in the third inequality we let $U \sim \mathcal{N}\left(0,(2\bar \lambda /\delta_1^2) \I_{|\M|}\right)$, applied the Gaussian concentration inequality and required $C > C_2=\max(2\sqrt{2\bar \lambda}/\delta_1,C')$. The whole inequality holds uniformly for all $\M \supseteq \M_0$, and all sufficiently large $n$, w.p. $1-2\eta$. \hfill$\blacksquare$

\vspace{+0.5cm}

\begin{lemma}\label{g1}
Suppose Assumptions 1-8 holds. Then w.p.a.1 as $n\to \infty$, uniformly for all $\M \supseteq \M_0$,
\begin{equation}\label{bfq1}
\BF[\M:\M_0] \asymp  \left(\frac{2\pi}{n}\right)^{-\frac{|\M|-|\M_0|}{2}} \frac{S_{\M}(\D)}{S_{\M_0}(\D)} \cdot
\frac{\left\{\dett (\G_{\M}^\top\V_n^{-1}\G_{\M})\right\}^{-1/2}\pi(\theta_0|\M)}
{\left\{\dett (\G_{\M_0}^\top\V_n^{-1}\G_{\M_0})\right\}^{-1/2}\pi(\theta_0|\M_0)},
\end{equation}
where $S_{\M}(\D)$ is defined in Lemma \ref{g11}. Moreover,
\begin{equation}\label{logss}
0<\ln \frac{S_{\M}(\D)}{S_{\M_0}(\D)} \asymp (|\M|-|\M_0|).
\end{equation}
\end{lemma}
\noindent \textbf{Proof:}
We first establish an approximation of the integral on the true model space $\Theta(\M_0)$. From Lemma \ref{g11}, \ref{g12} and \ref{g13}, we can pick a large constant $C>\max(C_1,C_2)$ such that for all sufficiently large $n$,
\begin{align}\label{approx}
& \int_{ \Theta(\M_0)}e^{-\frac{n}{2}\bar g(\D,\theta)^\top \V_n^{-1} \bar g(\D,\theta)} \pi(\theta|\M_0) \ud\theta \nonumber\\
={} & \big(1+o_p(1)\big ) S_{\M_0}(\D) (2\pi/n)^{|\M_0|/2}\big\{\dett (\G_{\M_0}^\top\V_n^{-1}\G_{\M_0})\big\}^{-1/2}\pi(\theta_0|\M_0)
\end{align}
This is because the density at $\theta_0$ on $\M_0$ is lower bounded by $e^{-c_0k_0}$ by Assumption 7(iii), and hence the upper bound in Lemma \ref{g13} is of smaller order compared to the right hand side of \eqref{laplace}, which implies that the integral on the space $\Theta(\M_0)$ is mostly concentrated on the neighborhood $B(C\e)$ and the outside part is negligible. \\

For any $\M \supseteq \M_0$ and $\M\neq \M_0$, we can decompose the Bayes factor in two parts:
\begin{align}\label{po1}
 \BF[\M:\M_0] & =  \frac{\int_{ \Theta(\M)\cap B_0(C\e)}e^{-\frac{n}{2}\bar g(\D,\theta)^\top \V_n^{-1} \bar g(\D,\theta)} \pi(\theta|\M) \ud\theta}{\int_{ \Theta(\M_0)}e^{-\frac{n}{2}\bar g(\D,\theta)^\top \V_n^{-1} \bar g(\D,\theta)} \pi(\theta|\M_0) \ud\theta} \nonumber \\
& + \frac{\int_{ \Theta(\M)\backslash B_0(C\e)}e^{-\frac{n}{2}\bar g(\D,\theta)^\top \V_n^{-1} \bar g(\D,\theta)} \pi(\theta|\M) \ud\theta}{\int_{ \Theta(\M_0)}e^{-\frac{n}{2}\bar g(\D,\theta)^\top \V_n^{-1} \bar g(\D,\theta)} \pi(\theta|\M_0) \ud\theta}\nonumber \\
:= {} & I_1+I_2
\end{align}
Based on Lemma \ref{g11}, Lemma \ref{g12} and \eqref{approx}, $I_1$ can be bounded by
\begin{align}\label{I11}
I_1 &= \frac{\big(1+o_p(1)\big ) S_{\M}(\D) (2\pi/n)^{|\M|/2}\big\{\dett (\G_{\M}^\top\V_n^{-1}\G_{\M})\big\}^{-1/2}\pi(\theta_0|\M)}{\big(1+o_p(1)\big ) S_{\M_0}(\D) (2\pi/n)^{|\M_0|/2}\big\{\dett (\G_{\M_0}^\top\V_n^{-1}\G_{\M_0})\big\}^{-1/2}\pi(\theta_0|\M_0)},
\end{align}
where the $o_p(1)$ holds uniformly for all $M\supset M_0$. Now we analyze the term $S_{\M}(\D)/S_{\M_0}(\D)$ and prove \eqref{logss}. According to the definition of $S_{\M}(\D)$ in Lemma \ref{g11},
\begin{align}\label{ss1}
\frac{S_{\M}(\D)}{S_{\M_0}(\D)}&= \exp\Bigg\{\frac{n}{2}\bar g(\D,\theta_0)^\top \big(\V_n^{-1}\G_{\M}(\G_{\M}^\top\V_n^{-1}\G_{\M})^{-1}\G_{\M}^\top\V_n^{-1} \nonumber \\
&\quad  - \V_n^{-1}\G_{\M_0}(\G_{\M_0}^\top\V_n^{-1}\G_{\M_0})^{-1}\G_{\M_0}^\top\V_n^{-1} \big)\bar g(\D,\theta_0)\Bigg\}\nonumber \\
&= \exp\Bigg\{\frac{n}{2}\bar g(\D,\theta_0)^\top \V_n^{-1/2}(\P_{\M}-\P_{\M_0})\V_n^{-1/2}\bar g(\D,\theta_0)\Bigg\},
\end{align}
where $\P_{\M}$ is the projection matrix defined at the beginning of the proof of Lemma \ref{g12}. Given $M\supset M_0$, $\P_{\M}-\P_{\M_0}$ is semi-positive definite and idempotent, with trace $|\M|-|\M_0|$. Since by CLT, $\sqrt{n} \V_n^{-1/2} \bar g(\D,\theta_0)$ converges in distribution to $\mathcal{N}(0,I_p)$, it follows that $2\ln \left(S_{\M}(\D)/S_{\M_0}(\D)\right)$ is asymptotically a $\chi^2_{|\M|-|\M_0|}$ random variable. Hence \eqref{logss} is proved.

For $I_2$, Lemma \ref{g13}, \eqref{approx} and Assumption 7 together yield
\begin{align}
I_2 &\leq \frac{c_{\pi} \Big(\frac{4\pi\bar \lambda }{n\delta_1^2}\Big)^{|\M|/2}\exp\Big(-\frac{C^2\delta_1^2}{16\bar \lambda }p\Big)+ \exp\Big(-\frac{n}{4}\bar \lambda ^{-1}\delta_0 ^2\Big)}{\big(1+o_p(1)\big ) S_{\M_0}(\D) (2\pi/n)^{|\M_0|/2}\big[\dett (\G_{\M_0}^\top\V_n^{-1}\G_{\M_0})\big]^{-1/2}\pi(\theta_0|\M_0)}.\nonumber
\end{align}
Hence we take the ratio of $I_2$ to $I_1$ in \eqref{I11}, and have
\begin{align}\label{I1I2}
\frac{I_2}{I_1}
&\leq \frac{c_{\pi} \Big(\frac{4\pi\bar \lambda }{n\delta_1^2}\Big)^{|\M|/2}\exp\Big(-\frac{C^2\delta_1^2}{16\bar \lambda }p\Big)+ \exp\Big(-\frac{n}{4}\bar \lambda ^{-1}\delta_0 ^2\Big)}
{\big(1+o_p(1)\big ) S_{\M}(\D) (2\pi/n)^{|\M|/2}\big\{\dett (\G_{\M}^\top\V_n^{-1}\G_{\M})\big\}^{-1/2}\pi(\theta_0|\M)}\nonumber\\
& \leq\frac{ \big(1+o_p(1)\big )c_{\pi} \Big(\frac{2\bar \lambda}{\delta_1^2}\Big)^{|\M|/2}\exp\Big(-\frac{C^2\delta_1^2}{16\bar \lambda }p\Big)\cdot \Big(\frac{\bar \lambda(\G^\top\G)}{\underline \lambda (\V_n)}\Big)^{|\M|/2}}{\exp\Big\{-\frac{n}{2}\bar g(\D,\theta_0)^\top \V_n^{-1}\bar g(\D, \theta_0)\Big\} \cdot e^{-c_0|\M|}},
\end{align}
where in the second inequality, we applied Assumption 7(iii), and also a lower bound on $S_{\M}(\D)$ using its definition in Lemma \ref{g11}. In fact, by \eqref{g0}
$$\frac{n}{2}\bar g(\D,\theta_0)^\top \V_n^{-1}\bar g(\D, \theta_0)=O_p(p).$$
So one can see that in \eqref{I1I2}, w.p.a.1 as $n\to\infty$, we can pick $C$ sufficiently large, such that $I_2/I_1$ is arbitrarily small, since the exponential index $|\M|$ cannot exceed $p$. Therefore, in \eqref{po1}, $\BF[\M:\M_0]=(1+o_p(1))I_1$, and the conclusion of \eqref{bfq1} follows from this and \eqref{I11}. \hfill$\blacksquare$

\vspace{+0.5cm}

\begin{lemma}\label{g2}
Suppose Assumptions 1-8 holds. Then w.p.a.1 as $n\to\infty$, there exists a large constant $C_1>0$, such that uniformly over all $\M$ with $\M_0\backslash \M\neq \emptyset$,
\begin{equation}\label{bfq2}
\BF[\M:\M_0]\leq \exp\left(-C_1 n\min_{j\in \M_0} \theta_{0,(j)}^2\right).
\end{equation}
\end{lemma}
\noindent \textbf{Proof:}
First we observe that if $\M_0\backslash \M\neq \emptyset$, i.e. $\M$ misses at least one component of the true model $\M_0$, then for any $\theta\in \Theta(\M)$, it must hold that $\|\theta-\theta_0\|\geq \min_{j:\theta_{0,(j)}\neq 0} |\theta_{0,(j)}|$. By Assumption 3, there exists a sequence $t_n\to \infty$ such that $\min_{j\in \M_0} |\theta_{0,(j)}|= \sqrt{\ln n} t_n\e$. Therefore for $n$ sufficiently large, the whole space $\Theta(\M)$ is outside the neighborhood $B_0(\sqrt{\ln n} t_n\e)$. Similar to the derivation of \eqref{out}, we can bound the marginal probability $q(\D|\M)$ by
\begin{align}\label{out2}
&  \int_{ \Theta(\M)}e^{-\frac{n}{2}\bar g(\D,\theta)^\top \V_n^{-1} \bar g(\D,\theta)} \pi(\theta|\M) \ud\theta \nonumber \\
= {} &  \int_{ \Theta(\M)\backslash B_0(\sqrt{\ln n} t_n\e)}e^{-\frac{n}{2}\bar g(\D,\theta)^\top \V_n^{-1} \bar g(\D,\theta)} \pi(\theta|\M) \ud\theta \nonumber \\
\leq {}& c_{\pi} \Big(\frac{4\pi\bar \lambda }{n\delta_1^2}\Big)^{|\M|/2}\exp\left(-\frac{C^2\delta_1^2}{16\bar \lambda}t_n^2 p\ln n\right)+ \exp\left(-\frac{n}{4}\bar \lambda ^{-1}\delta_0 ^2\right)
\end{align}
for sufficiently large $C$.

Therefore, by using the approximation \eqref{approx}, we have that w.p.a.1 as $n\to\infty$,
\begin{align*}
& \BF[\M:\M_0]\\
\leq {}& \frac{c_{\pi} \Big(\frac{4\pi\bar \lambda }{n\delta_1^2}\Big)^{|\M|/2}\exp\Big(-\frac{C^2\delta_1^2}{16\bar \lambda}t_n^2 p\ln n\Big)+ \exp\Big(-\frac{n}{4}\bar \lambda ^{-1}\delta_0 ^2\Big)}
{\big(1+o_p(1)\big ) S_{\M_0}(\D) (2\pi/n)^{|\M_0|/2}\big\{\dett (\G_{\M_0}^\top\V_n^{-1}\G_{\M_0})\big\}^{-1/2}\pi(\theta_0|\M_0)}\\
\leq {}& \frac{\big(1+o_p(1)\big )c_{\pi} \Big(\frac{2\bar \lambda }{\delta_1^2}\Big)^{|\M|/2}\exp\Big(-\frac{C^2\delta_1^2}{16\bar \lambda}t_n^2 p\ln n\Big)}
{\exp\Big\{-\frac{n}{2}\bar g(\D,\theta_0)^\top \V_n^{-1}\bar g(\D,\theta_0)\Big\} \Big(\frac{\underline \lambda(\V_n)}{\bar \lambda (\G^\top\G)}\Big)^{k_0/2} e^{-c_0k_0}},
\end{align*}
where in the last inequality, we did the same as in \eqref{I1I2}, and the second term on the numerator is absorbed into the first term because $t_n^2 p\ln n = n\min_{j\in \M_0}\theta^2_{0,j}\preceq n $  by Assumption 3. Now because the term $t_n^2 p\ln n$ dominates all the other terms in the exponential, the conclusion follows by choose appropriate $C_1>0$.  \hfill$\blacksquare$

\vspace{+0.5cm}

\noindent {\bf Proof of Theorem 1 (i):}\\
We combine the results from Lemma \ref{g1} and Lemma \ref{g2}, and show that $\sum_{\M\neq \M_0} \PO[\M:\M_0]\to 0$ w.p.a.1 as $n\to \infty$, since this is equivalent to $q(\M_0|\D)\to 1$ w.p.a.1 as $n\to \infty$. To prove this, it suffices to show $\sum_{\M:\M_0\backslash \M\neq \emptyset}\PO[\M:\M_0]\to 0$ and $\sum_{\M:\M \supseteq \M_0} \PO[\M:\M_0]\to 0$ respectively. For those models with $\M_0\backslash \M\neq \emptyset$, Assumption 8(ii) implies that\\
$\sup_{\{\M: \ \M_0\backslash \M \neq \emptyset\}} \pi(\M)/\pi(\M_0) \leq e^{r_1p\ln n}$ for sufficiently large $n$, and Assumption 3 implies that $2r_1 p\ln n \leq C_1 n\min_{j\in \M_0} \theta_{0,(j)}^2$ for sufficiently large $n$. Therefore it follows from Lemma \ref{g2} that as $n\to\infty$,
\begin{align*}
&  \sum_{\M:\M_0\backslash \M\neq \emptyset} \PO[\M:\M_0]  \\
\leq {} &  \sum_{\M:\M_0\backslash \M\neq \emptyset} \frac{\pi(\M)}{\pi(\M_0)} \exp\left(-C_1 n\min_{j\in \M_0} \theta_{0,(j)}^2\right)\\
\leq {} &  2^p \cdot\exp\left(r_1p\ln n -C_1 n\min_{j\in \M_0} \theta_{0,(j)}^2 \right) \\
\leq {}& \exp\left(p\ln 2 + r_1p\ln n -2r_1 p\ln n\right) \to 0.
\end{align*}

For those models with $\M \supseteq \M_0$, we use the conclusion of Lemma \ref{g1}. Note that on the right-hand side of \eqref{bfq1}, $\pi(\theta_0|\M)/\pi(\theta_0|\M_0)\leq \exp(c_1(|\M|-|\M_0|))$ by Assumption 7(iii). Based on Assumptions 5, 6, \eqref{bfq1} and \eqref{logss}, we can pick a large constant $C_2>0$, such that uniformly over all these models, w.p.a.1 as $n\to \infty$,
$$\BF[\M:\M_0]\leq (C_2n)^{-\frac{|\M|-|\M_0|}{2}}.$$
By Assumption 8(i), we know that for any $\M\supset \M_0$, $\frac{\pi(\M)}{\pi(\M_0)} = |o(1)| \cdot \left(\sqrt{n}/p\right)^{|\M |-| \M_0|}$ as $n\to \infty$, where $o(1)$ is uniform for all such $\M$. Therefore we have
\begin{align*}
&  \sum_{\M:\M \supset \M_0} \PO[\M:\M_0]  \\
\leq {} & \sum_{\M:\M \supset \M_0} \frac{\pi(\M)}{\pi(\M_0)} (C_2n)^{-\frac{|\M|-|\M_0|}{2}}\\
\leq {} & |o(1)|\cdot \sum_{k=k_0+1}^{p} \binom{p-k_0}{k-k_0} \left(\sqrt{n}/p\right)^{k-k_0} (C_2n)^{-\frac{k-k_0}{2}}\\
\leq {} & |o(1)|\cdot \sum_{k=k_0+1}^{p} \binom{p-k_0}{k-k_0}  \left(\frac{1}{\sqrt{C_2} p}\right)^{k-k_0}\\
\leq {} & |o(1)|\cdot \left\{ \left(1+\frac{1}{\sqrt{C_2} p} \right)^{p-k_0} - 1 \right\}\\
\leq {}&  |o(1)|\cdot (e^{\frac{1}{\sqrt{C_2}}}-1) \to  0.
\end{align*}
So the proof is complete. \hfill$\blacksquare$

\vspace{+0.5cm}

\noindent \textbf{Proof of Theorem 1 (ii):}\\
In the conclusion of part (ii), the first integral can be rewritten as
\begin{align}
& \int_{A}q(\theta|\D)\ud \theta \nonumber \\
= {}& \frac{\sum_{\M} \pi(\M)\int_{A\cap \Theta(\M)} q(\D|\theta,\M)\pi(\theta|\M)\ud \theta}{\sum_{\M} \pi(\M)\int_{\Theta(\M)} q(\D|\theta,\M)\pi(\theta|\M)\ud \theta} \nonumber \\
= {}& \frac{\sum_{\M\neq \M_0} \pi(\M)\int_{A\cap \Theta(\M)} q(\D|\theta,\M)\pi(\theta|\M)\ud \theta + \pi(\M_0)\int_{A\cap \Theta(\M_0)} q(\D|\theta,\M)\pi(\theta|\M_0)\ud \theta}{\sum_{\M\neq \M_0} \pi(\M)\int_{ \Theta(\M)} q(\D|\theta,\M)\pi(\theta|\M)\ud \theta + \pi(\M_0)\int_{\Theta(\M_0)} q(\D|\theta,\M)\pi(\theta|\M_0)\ud \theta}. \nonumber
\end{align}
Therefore if we divide the numerator and denominator by $\pi(\M_0)\int_{\Theta(\M_0)} q(\D|\theta,\M)\pi(\theta|\M_0)\ud \theta$, we have
$$\frac{\int_{A\cap \Theta(\M_0)} \tilde q(\theta_1|\D)\ud \theta_1 }{\sum_{\M\neq \M_0} \PO[\M:\M_0] + 1} \leq \int_{A}q(\theta|\D)\ud \theta \leq \frac{\sum_{\M\neq \M_0} \PO[\M:\M_0] + \int_{A\cap \Theta(\M_0)} \tilde q(\theta_1|\D)\ud \theta_1 }{\sum_{\M\neq \M_0} \PO[\M:\M_0] + 1},$$
where $\tilde q(\theta_1|\D)=\frac{q(\D|\theta,\M)\pi(\theta|\M_0)}{\int_{\Theta(\M_0)} q(\D|\theta,\M)\pi(\theta|\M_0)\ud \theta}$ is the conditional quasi-posterior on the true model space $\M_0$ and $\theta_1$ represents the nonzero components of $\theta$ in the true model $\M_0$. According to the model selection consistency of Theorem 1 part (i), $\sum_{\M\neq \M_0} \PO[\M:\M_0] \to 0$ w.p.a.1 as $n\to \infty$. Hence
$$\frac{\int_{A\cap \Theta(\M_0)} \tilde q(\theta_1|\D)\ud \theta_1 }{ 1+|o_p(1)|} \leq \int_{A}q(\theta|\D)\ud \theta \leq \frac{ \int_{A\cap \Theta(\M_0)} \tilde q(\theta_1|\D)\ud \theta_1 +|o_p(1)|}{ 1+|o_p(1)| }.$$
Since $\int_{\Theta(\M_0)} \tilde q(\theta_1|\D)\ud \theta_1=1$ and $\int_{\Theta} q(\theta|\D)\ud \theta=1$, and the $o_p(1)$ does not depend on the set $A$, the inequality above implies that for all set $A\subseteq \Theta$, w.p.a.1 as $n\to \infty$.
$$\Big| \int_{A}q(\theta|\D)\ud \theta - \int_{A\cap \Theta(\M_0)} \tilde q(\theta_1|\D)\ud \theta_1\Big|\to 0 .$$
Therefore, to show part (ii) of Theorem 1, it suffices to show that w.p.a.1. as $n\to \infty$,
\begin{equation}\label{eii}
\sup_{A\subseteq \Theta} \Big| \int_{A\cap \Theta(\M_0)} \tilde q(\theta_1|\D)\ud \theta_1 - \int_{A\cap \Theta(\M_0)} \phi\big(\theta_1;\bar \theta_{\M_0,1}, (\G_{\M_0}^\top \V_n^{-1} \G_{\M_0})^{-1}/n\big) \ud \theta_1\Big|\to 0.
\end{equation}
Note that the densities $\tilde q$ and $\phi$ are defined on the same support $\Theta(\M_0)$, so the rest is a standard proof of Bayesian CLT similar to \citet{bc09}. Using the decomposition \eqref{gdecomp} in Lemma \ref{g11}, Lemma \ref{g12} and Lemma \ref{g13}, we have
\begin{align}\label{qtilde}
 & \int_{A\cap \Theta(\M_0)} \tilde q(\theta_1|\D)\ud \theta_1 \nonumber \\
 ={}&\frac{S_{\M_0}(\D)\int_{A\cap \Theta(\M_0)} e^{-\frac{n}{2} (\theta_1-\bar  \theta_{\M_0,1})^\top \G_{\M_0}^\top\V_n^{-1}\G_{\M_0}(\theta_1-\bar \theta_{\M_0,1})+o_p(1)} \pi(\theta_1|\M_0)\ud \theta_1 } {\big(1+o_p(1)\big ) S_{\M_0}(\D) (2\pi/n)^{|\M_0|/2}\big\{\dett (\G_{\M_0}^\top\V_n^{-1}\G_{\M_0})\big\}^{-1/2}\pi(\theta_0|\M_0)} \nonumber \\
 ={}&\frac{ \big(1+o_p(1)\big )\int_{A\cap B_0(C\e)\cap \Theta(\M_0)} e^{-\frac{n}{2} (\theta_1-\bar \theta_{\M_0,1})^\top \G_{\M_0}^\top\V_n^{-1}\G_{\M_0}(\theta_1 -\bar \theta_{\M_0,1})} \pi(\theta_1|\M_0)\ud \theta_1 }{(2\pi/n)^{|\M_0|/2}\big\{\dett (\G_{\M_0}^\top\V_n^{-1}\G_{\M_0})\big\}^{-1/2}\pi(\theta_0|\M_0)} \nonumber \\
= {}&\big(1+o_p(1)\big )\int_{A\cap B_0(C\e)\cap \Theta(\M_0)} \phi\big(\theta_1;\bar \theta_{\M_0,1},(\G_{\M_0}^\top\V_n^{-1}\G_{\M_0})^{-1}/n \big)\ud \theta_1 \nonumber \\
 ={}& \int_{A\cap B_0(C\e)\cap \Theta(\M_0)} \phi\big(\theta_1;\bar \theta_{\M_0,1},(\G_{\M_0}^\top\V_n^{-1}\G_{\M_0})^{-1}/n \big)\ud \theta_1 + o_p(1)
\end{align}
The numerator of the second equality shrinks the range of integral to within the neighborhood $B_0(C\e)$ because the integral outside $B_0(C\e)$ is of $o_p(1)$ compared to the denominator, according to the approximation in \eqref{approx}. In the third equality, we used Assumption 7(ii) and have that on $B_0(C\e)\cap\Theta(\M_0)$, $\pi(\theta|\M_0)/\pi(\theta_0|\M_0)=1+o(1)$. The $o_p(1)$ in the last expression does not depend on the set $A$. Therefore \eqref{eii} holds and this completes the proof. \hfill$\blacksquare$

\vspace{+1.5cm}

\noindent{\bf \Large 2. Proofs for the Asymptotic Validity of BGMM}

\vspace{+0.5cm}

In this section, we give the proofs of Theorem 2 and Theorem 3. For the ease of notation, let $\si(\theta) = \G(\theta)^\top\V(\theta)^{-1}\G(\theta) $, where $\G(\theta)$ and $\V(\theta)$ are defined in Assumption 9. In the following, for any matrix $\A(\theta)$ that depends on $\theta$, we use the notation ``$\A$" to refer to the matrix evaluated at $\theta_0$. For example, $\si=\G^\top\V^{-1}\G$, where $\G$ and $\V$ are defined in Assumption 5 and 6, i.e. $\G(\theta)$ and $\V(\theta)$ evaluated at $\theta=\theta_0$, respectively. For any model $\M$, one can partition any matrix $\G(\theta)$ into $\G(\theta)=(\G_{\M}(\theta),\G_{\M^c}(\theta))$, according to the partial derivative with respect to the components of $\theta$ in either $\M$ or $\M^c$. Let
\begin{align*}
\si_{11}(\theta) &= \G_{\M}(\theta)^\top\V(\theta)^{-1}\G_{\M}(\theta)\\
\si_{12}(\theta) &= \G_{\M}(\theta)^\top\V(\theta)^{-1}\G_{\M^c}(\theta)\\
\si_{22}(\theta) &= \G_{\M^c}(\theta)^\top\V(\theta)^{-1}\G_{\M^c}(\theta)\\
\J_{\M}(\theta) &=  \I_p-\V(\theta)^{-1}\G_{\M}(\theta) (\G_{\M}(\theta)^\top \V(\theta)^{-1} \G_{\M}(\theta))^{-1}\G_{\M}(\theta)^\top
\end{align*}


We have the following lemma about the quadratic term in the asymptotic normal density of the GMM estimator $\hat \theta$. The proof is straightforward algebra.

\begin{lemma}\label{matalg}
Under Assumptions 1-10,
\begin{align}
&\frac{n}{2}(\theta-\hat \theta)^\top \si(\theta)(\theta-\hat \theta)  = \frac{n}{2}(\theta_1-\xi_1(\theta))^\top \si_{11}(\theta) (\theta_1-\xi_1(\theta)) + T_{\M}(\theta) \nonumber \\
&\xi_1(\theta)  := \hat \theta_1 + \si_{11}(\theta)^{-1} \si_{12}(\theta) \hat \theta_2 \nonumber \\
&T_{\M}(\theta)  := \frac{n}{2} \hat \theta_2^\top \G_{\M^c}(\theta)^\top\J_{\M}(\theta) \V(\theta)^{-1}\G_{\M^c}(\theta)\hat \theta_2 \nonumber
\end{align}
where $\hat \theta=(\hat \theta_1^\top,\hat \theta_2^\top)^\top$ is the GMM estimator on the full model space in $\mathbb{R}^p$ and is decomposed according to the model $\M$.
\end{lemma}

\vspace{+.6cm}

\begin{lemma}\label{xi1}
Suppose Assumptions 1-10 hold. Then uniformly for all $\theta \in B_0(C\e)\cap \Theta(\M)$ and all $\M\supseteq \M_0$, with any fixed constant $C>0$, for $\xi_1(\theta)$ in Lemma \ref{matalg},
$$\xi_1(\theta) = \bar \theta_{\M,1} + o_p\big(1/\sqrt{n}\big),$$
where $\bar \theta_{\M,1} = \theta_{0,\M,1}-(\G_{\M}^\top\V_n^{-1}\G_{\M})^{-1}\G_{\M}^\top \V_n^{-1}\bar g(\D,\theta_0)$ is the same as in Lemma \ref{g11}. Therefore,
$$\frac{n}{2}(\theta-\hat \theta)^\top \si(\theta)(\theta-\hat \theta)  = \frac{n}{2}(\theta_1- \bar \theta_{\M,1})^\top \si_{11}(\theta) (\theta_1-\bar \theta_{\M,1}) + T_{\M}(\theta) +o_p(1).$$
\end{lemma}
\noindent {\bf Proof:}
First we use the continuity of $\G^\top(\theta)\G(\theta)$ and $\V(\theta)$ in $\theta$ from Assumption 9(iii), and replace $\si_{11}(\theta)$ and $\si_{12}(\theta) $ in the expression of $\xi_1(\theta)$ by $\si_{11}$ and $\si_{12}$, respectively. This is because we are considering $\theta \in B_0(C\e)$, and this leads to (note that $\hat \theta_2 = O_p(1/\sqrt{n})$)
$$\xi_1(\theta) = \hat \theta_1 + (\si_{11}^{-1} \si_{12}+o_p(1)) \hat \theta_2 =\hat \theta_1 + \si_{11}^{-1} \si_{12} \hat \theta_2 + o_p(1/\sqrt{n}).$$
Next we use Assumption 9(v), and replace $\hat \theta_1$ and $\hat \theta_2$ with their first order approximations. For the $\bar \theta $ in Assumption 9(v), suppose we decompose it into $\bar \theta =(\bar \theta_1^\top, \bar \theta_2^\top)^\top$ according to a given model $\M$. Then by Assumption 9(v), we have $\|\hat \theta_1 - \bar \theta_1\|=O_p(1/n)$ and $\|\hat \theta_2 - \bar \theta_2\|=O_p(1/n)$. Furthermore, through pure matrix algebra, we can derive that
\begin{align*}
\xi_1(\theta) & = \hat \theta_1 + \si_{11}^{-1} \si_{12} \hat \theta_2 + o_p(1/\sqrt{n})\\
& = \bar \theta_1 + \si_{11}^{-1} \si_{12} \bar  \theta_2 + o_p(1/\sqrt{n}) + O_p(1/n) \\
&= (\I_{|\M|},\si_{11}^{-1} \si_{12})\bar \theta + o_p(1/\sqrt{n}) \\
& = \bar \theta_{\M,1}  + o_p(1/\sqrt{n}),
\end{align*}
where in the last display, $ \bar \theta_{\M,1}$ is defined in Lemma \ref{g11}. Therefore from Lemma \ref{matalg}, for any $\theta \in B_0(C\e)$
\begin{align*}
 & \frac{n}{2}(\theta-\hat \theta)^\top \si(\theta)(\theta-\hat \theta)  \\
 ={}& \frac{n}{2}(\theta_1-\xi_1(\theta))^\top \si_{11}(\theta) (\theta_1-\xi_1(\theta)) + T_{\M}(\theta) \\
 ={}& \frac{n}{2}(\theta_1-\bar \theta_{\M,1}  + o_p(1/\sqrt{n}))^\top \si_{11}(\theta) (\theta_1-\bar \theta_{\M,1}  + o_p(1/\sqrt{n})) + T_{\M}(\theta) \\
 ={} & \frac{n}{2}(\theta_1-\bar \theta_{\M,1} )^\top \si_{11}(\theta) (\theta_1-\bar \theta_{\M,1} ) + T_{\M}(\theta) +o_p(1),
\end{align*}
which completes the proof.  \hfill$\blacksquare$

\vspace{+1cm}

\begin{lemma}\label{tm}
Suppose Assumptions 1-10 hold. Let $T_{\M}(\theta)$ be given in Lemma \ref{matalg}. For any generic model $\M$, w.p.a.1 as $n\to \infty$, \\
\noindent (i) Uniformly for all $\theta \in B_0(C\e)\cap \Theta(\M)$ and all models $\M \supseteq \M_0$,  $T_{\M}(\theta)-T_{\M}(\theta_0)=o_p(1)$, given any fixed constant $C>0$;\\
\noindent (ii) For any $\M \supseteq \M_0$, $T_{\M_0}(\theta_0)-T_{\M}(\theta_0)=\ln \left(S_{\M}(\D)/S_{\M_0}(\D)\right)+o_p(1)$, where $S_{\M}(\D)$ is defined in Lemma \ref{g11}.
\end{lemma}
\noindent {\bf Proof:}
First, we can use Assumption 9 and the uniform boundedness of the eigenvalues of $\G(\theta)^\top\G(\theta)$ and $\V(\theta)$ for $\theta \in \Theta$, and express $T_{\M}(\theta)$ as
\begin{align}\label{tmbar}
&T_{\M}(\theta) \nonumber \\
={}& \frac{n}{2} \bar \theta_2^\top \G_{\M^c}(\theta)^\top\J_{\M}(\theta)\V(\theta)^{-1} \G_{\M^c}(\theta)\bar \theta_2 + n \bar \theta_2^\top \G_{\M^c(\theta)}^\top\J_{\M}(\theta) \V(\theta)^{-1}\G_{\M^c}(\theta) (\hat \theta_2-\bar \theta_2) \nonumber \\
& + \frac{n}{2}(\hat \theta_2 -\bar \theta_2 )  ^\top \G_{\M^c}(\theta)^\top\J_{\M}(\theta)\V(\theta)^{-1} \G_{\M^c} (\theta)(\hat \theta_2 -\bar \theta_2)\nonumber\\
 ={}& \frac{n}{2} \bar \theta_2^\top \G_{\M^c}(\theta)^\top\J_{\M}(\theta)\V(\theta)^{-1} \G_{\M^c}(\theta)\bar \theta_2 + O_p(1/\sqrt{n}) + O_p(1/n) \nonumber\\
={}& \frac{n}{2} \bar \theta_2^\top \G_{\M^c}(\theta)^\top\J_{\M} (\theta)\V(\theta)^{-1}\G_{\M^c}(\theta)\bar \theta_2 +o_p(1),
\end{align}
where $\bar \theta$ is the same as defined in Theorem 1 and $\bar \theta_2$ is the subvector of $\bar \theta$ with all those components not contained in model $\M$.\\
Furthermore, if $\M \supseteq \M_0$, then $\theta_{0,\M,2}=0$. It can be shown by straightforward matrix algebra that in this case,
\begin{equation}\label{tbar2}
\bar \theta_2 = (\G_{\M^c}^\top \J_{\M} \V^{-1} \G_{\M^c})^{-1} \G_{\M^c}^\top \J_{\M} \V^{-1} \bar g(\D,\theta_0).
\end{equation}
Given the expression of $T_{\M}(\theta)$ in \eqref{tmbar} and $\bar \theta_2$ in \eqref{tbar2}, it is clear that the dependence of $T_{\M}(\theta)$ on $\theta$ is only through the weighting matrix $\G_{\M^c}(\theta)^\top\J_{\M}(\theta) \V(\theta)^{-1}\G_{\M^c}(\theta)$ up to an error of $o_p(1)$ that is uniform for all $\M$. Since Assumption 9(iii) has assumed the continuity of $\G(\theta)$ and $\V(\theta)$ with respect to $\theta$, it follows that for any fixed model $\M$, $\G_{\M^c}(\theta)^\top\J_{\M}(\theta) \V(\theta)^{-1}\G_{\M^c}(\theta)$ is also continuous in $\theta$. Moreover, due to the boundedness of $p$ in Assumption 9(i), we have at most $2^{\bar p}$ models $\M$, so the continuity is uniform in $\M$. Therefore, in the shrinking neighborhood $\theta \in B_0(C\e)$ where $\e \to 0$, we have that uniformly over all models $\M\supseteq \M_0$, $T_{\M}(\theta)-T_{\M}(\theta_0)=o_p(1)$, which has proved (i).\\

By \eqref{tmbar} and \eqref{tbar2}, for any $\M\supseteq \M_0$, we have $T_{\M}(\theta_0)=\frac{n}{2}\bar g(\D,\theta_0)^\top \K(\M) \bar g(\D,\theta_0)+o_p(1)$, where
$$\K(\M)=\V^{-1} \J_{\M}^\top \G_{\M^c}( \G_{\M^c}^\top \J_{\M} \V^{-1} \G_{\M^c})^{-1} \G_{\M^c}^\top \J_{\M} \V^{-1}.$$
By pure matrix algebra, it can be shown that
\begin{align*}
& \K(\M_0)-\K(\M) \\
={}& \V^{-1}\G_{\M}(\G_{\M}^\top\V^{-1}\G_{\M})^{-1}\G_{\M}^\top\V^{-1} - \V^{-1}\G_{\M_0}(\G_{\M_0}^\top\V^{-1}\G_{\M_0})^{-1}\G_{\M_0}^\top\V^{-1}
\end{align*}
Because $\V_n$ is a consistent estimator for $\V$ and their eigenvalues are bounded, we can replace $\V$ in the display above by $\V_n$, which will incur an error of $o_p(1)$. Then it follows that
\begin{align*}
& T_{\M_0}(\theta_0)-T_{\M}(\theta_0)  \\
={}& \frac{n}{2}\bar g(\D,\theta_0)^\top  \Big\{\V_n^{-1}\G_{\M}(\G_{\M}^\top\V_n^{-1}\G_{\M})^{-1}\G_{\M}^\top\V_n^{-1} \\
& - \V_n^{-1}\G_{\M_0}(\G_{\M_0}^\top\V_n^{-1}\G_{\M_0})^{-1}\G_{\M_0}^\top\V_n^{-1}\Big\}\bar g(\D,\theta_0)+o_p(1).
\end{align*}
We compare this with \eqref{ss1} and obtain that $T_{\M_0}(\theta_0)-T_{\M}(\theta_0)=\ln \left(S_{\M}(\D)/S_{\M_0}(\D)\right)+o_p(1)$. (ii) is proved.\hfill$\blacksquare$

\vspace{+0.5cm}

\begin{lemma}\label{lemdif}
Under Assumptions 1-11, uniformly for all models $\M$,
\begin{align}\label{diff}
& \int_{\Theta(\M)}\Big|p(\hat \theta|\theta_1)-\phi(\hat \theta;\theta_1, \si(\theta)^{-1}/n) \Big| \pi(\theta_1|\M)\ud \theta_1= O_p\Big(\tau_n n^{\frac{p-|\M|}{2}}\Big),
\end{align}
where $\theta=(\theta_1^\top,\theta_2^\top)$ is decomposed according to $\M$, $p(\hat\theta|\theta_1)$ denotes the conditional density of $\hat\theta$ given $\theta=(\theta_1^\top,0)^\top$, and
\begin{align}\label{phiexpr}
&\phi(\hat \theta;\theta_1, \si(\theta)^{-1}/n) \nonumber \\
&= \left(\frac{2\pi}{n}\right)^{-\frac{p}{2}}\dett (\si(\theta))^{1/2}\exp\left\{-\frac{n}{2}(\hat \theta_1^\top-\theta_1^\top, \hat \theta_2^\top)\si(\theta)\binom{\hat \theta_1-\theta_1}{\hat \theta_2} \right\},
\end{align}
where $\hat\theta = (\hat\theta_1^\top,\hat\theta_2^\top)^\top$ is decomposed according to $\M$.
\end{lemma}
\noindent \textbf{Proof:}
Let $c(x)=1/(1+x^{p+1})$. First one can do a variable transformation from $\hat \theta$ to $Z=\sqrt{n}\F(\theta)(\hat \theta-\theta)$, where $\theta=(\theta_1^\top,0)^\top$ and $\F(\theta)^\top\F(\theta)= \G(\theta)^\top\V(\theta)^{-1}\G(\theta)=\si(\theta)$. The densities have the relation $p(\hat \theta|\theta_1)=n^{p/2}\dett(\F(\theta))p_Z(z|\theta_1)$ and $\phi\left(\hat \theta;\theta_1, \si(\theta)^{-1}/n\right) = n^{p/2}\dett(\F(\theta)) \phi(z;0, \I_p)$. Note that this transformation is in $\mathbb{R}^p$ and does not directly involve the integration with respect to $\theta_1$. Then using Assumption 10, the left-hand side of \eqref{diff} can be bounded by
\begin{align}\label{dif1}
& \int_{\Theta(\M)}\Big|p(\hat \theta|\theta_1)-\phi(\hat \theta;\theta_1, \si(\theta)^{-1}/n) \Big| \pi(\theta_1|\M)\ud \theta_1 \nonumber \\
= {} & n^{p/2} \sup_{\theta\in \Theta}\dett(\F(\theta)) \int_{\Theta(\M)} \Big|p_Z(z|\theta_1)-\phi(z;0, \I_p) \Big|\pi(\theta_1|\M) \ud \theta_1 \nonumber \\
\leq {} & n^{p/2}\sup_{\theta\in \Theta}\left(\frac{\bar \lambda(\G(\theta)^\top\G(\theta))}{\underline \lambda(\V(\theta))} \right)^{p/2}\cdot \tau_n \int_{\Theta(\M)}  c(\|z\|)  \pi(\theta_1|\M) \ud \theta_1 \nonumber \\
\leq {} & n^{p/2}\sup_{\theta\in \Theta}\left(\frac{\bar \lambda(\G(\theta)^\top\G(\theta))}{\underline \lambda(\V(\theta))} \right)^{p/2}\cdot \tau_n c_{\pi} \int_{\Theta(\M)}  c(\|z\|) \ud \theta_1,
\end{align}
where in the last inequality we have used Assumption 7(i) that $\pi(\theta_1|\M)\leq c_{\pi}$ uniformly for all $\M$. We now transform $\theta_1$ into $\beta=\sqrt{n}(\hat \theta_1-\theta_1)$ in the foregoing integration about the function $c(\|z\|)$. We have
$$\|z\|^2=n(\hat \theta-\theta)^\top \si(\theta) (\hat \theta-\theta)\geq \underline\lambda(\si(\theta)) \|\sqrt{n}(\hat \theta -\theta)\|^2,$$
and
$$\left\|\sqrt{n}(\hat \theta -\theta)\right\|^2=\left\|\sqrt{n}(\hat \theta_1-\theta_1)\right\|^2+ n\left\|\hat \theta_2\right\|^2 = \|\beta\|^2 + n\left\|\hat \theta_2\right\|^2 \geq \|\beta\|^2.$$
Furthermore, the eigenvalue satisfies
$$\inf_{\theta\in \Theta} \underline\lambda(\si(\theta)) \geq \inf_{\theta\in \Theta} \underline\lambda(\G(\theta)^\top\G(\theta))\bar \lambda(\V(\theta))^{-1},$$
which is lower bounded by constant according to Assumption 9(iii). Therefore, along with the nonincreasing property of the function $c(\cdot)$, we have
\begin{align}\label{dif2}
& \int_{\Theta(\M)}  c(\|z\|) \ud \theta_1 \leq n^{-\frac{|\M|}{2}}\int_{\mathbb{R}^{|\M|}}   c\left( \inf_{\theta\in \Theta} \sqrt{\frac{\underline\lambda(\G(\theta)^\top\G(\theta))}{\bar \lambda(\V(\theta))} } \|\beta\| \right)  \ud \beta \nonumber \\
\leq {} & n^{-\frac{|\M|}{2}} \Bigg(\inf_{\theta\in \Theta}\frac{\underline\lambda(\G(\theta)^\top\G(\theta))}{\bar \lambda(\V(\theta))} \Bigg)^{-\frac{|\M|}{2}}\int_0^{\infty} \frac{x^{|\M|-1}}{1+x^{p+1}}\ud x = O\left(n^{-\frac{|\M|}{2}}\right)
\end{align}
The conclusion follows from \eqref{dif1}, \eqref{dif2}, and the boundedness of the eigenvalues of $\si(\theta)$. \hfill$\blacksquare$

\vspace{+1cm}


\noindent {\bf Proof of Theorem 2 (i):}\\
We derive an order expression for $\bfh[\M:\M_0]=\frac{\int_{\Theta(\M)}p(\hat \theta|\theta_1)\pi(\theta_1|\M)\ud \theta_1}{\int_{\Theta(\M_0)}p(\hat\theta| \theta_1)\pi(\theta_1|\M_0)\ud \theta_1}$ where $\M\supseteq \M_0$. Here the two $\theta_1$'s in the numerator and the denominator lie in $\M$ and $\M_0$ respectively, possibly with different dimensions if $\M\supset \M_0$. Hereafter, all $o_p$ and $O_p$ hold uniformly over all the models with $\M\supseteq \M_0$. First of all, based on Lemma \ref{lemdif}, we have
\begin{align}\label{numer1}
&\int_{\Theta(\M)} p(\hat\theta| \theta_1)\pi(\theta_1|\M)\ud \theta_1 \nonumber \\
 \leq {}& \int_{\Theta(\M)} \Big| p(\hat\theta| \theta_1) - \phi(\hat \theta;\theta_1, \si(\theta)^{-1}/n) \Big|\pi(\theta_1|\M)\ud \theta_1 +
 \int_{\Theta(\M)} \phi(\hat \theta;\theta_1, \si(\theta)^{-1}/n) \pi(\theta_1|\M)\ud \theta_1 \nonumber \\
 \leq {}& O_p\Big(\tau_n n^{\frac{p-|\M|}{2}}\Big) +  \int_{\Theta(\M)} \phi(\hat \theta;\theta_1, \si(\theta)^{-1}/n) \pi(\theta_1|\M)\ud \theta_1
\end{align}
We claim that the second term in \eqref{numer1} satisfies
\begin{equation}\label{tnot1}
\left|\int_{\Theta(\M)} \phi(\hat \theta;\theta_1, \si(\theta)^{-1}/n) \pi(\theta_1|\M)\ud \theta_1
-  \int_{\Theta(\M)}  \phi(\hat \theta;\theta_1, \si^{-1}/n) \pi(\theta_1|\M) \ud \theta_1\right| =o_p\Big(n^{\frac{p-|\M|}{2}}\Big),
\end{equation}
where we have replaced $\si(\theta)$ with $\si$, i.e. the matrix $\si(\theta)$ evaluated at $\theta=\theta_0$. \\
To show \eqref{tnot1}, we first observe that both integrals in the display can be made to order $o_p\big(n^{\frac{p-|\M|}{2}}\big)$ outside the neighborhood $B_0(C\e)$ for some constant $C>0$. This is because by the decomposition in Lemma \ref{matalg},
\begin{align*}
& \int_{\Theta(\M)\backslash B_0(C\e)} \phi(\hat \theta;\theta_1, \si(\theta)^{-1}/n) \pi(\theta_1|\M)\ud \theta_1\\
= {}& \int_{\Theta(\M)\backslash B_0(C\e)} \left(\frac{2\pi}{n}\right)^{-p/2} \dett(\si(\theta))^{1/2}
e^{-\frac{n}{2}(\theta-\hat \theta)^\top \si(\theta) (\theta-\hat \theta)} \pi(\theta_1|\M)\ud \theta_1 \\
\leq {}& c_{\pi}\Big(\frac{2\pi}{n}\Big)^{-\frac{p}{2}}\sup_{\theta \in\Theta} \dett(\si(\theta))^{1/2} \int_{\Theta(\M)\backslash B_0(C\e)} \exp\left\{-\frac{n \inf_{\theta\in\Theta}\underline \lambda(\si(\theta))}{2}(\theta-\hat \theta)^\top (\theta-\hat \theta)\right\} \ud \theta_1 \\
\leq {}& c_{\pi}\left(\frac{2\pi}{n}\right)^{-\frac{p}{2}}\sup_{\theta \in\Theta} \dett(\si(\theta))^{1/2} \int_{\Theta(\M)\backslash B_0(C\e)} \exp\left\{-\frac{n \inf_{\theta\in\Theta}\underline \lambda(\si(\theta))}{2}(\theta_1-\hat \theta_1)^\top (\theta_1-\hat \theta_1)\right\} \ud \theta_1 \\
\leq {}&  c_{\pi}\left(\frac{2\pi}{n}\right)^{-\frac{p-|\M|}{2}}\frac{\sup_{\theta\in \Theta}\dett(\si(\theta))^{1/2}}{\inf_{\theta\in \Theta}\underline \lambda(\si(\theta))^{|\M|/2}}
\int_{\Theta(\M)\backslash B_0(C\e)} \phi\left(\theta_1;\hat \theta_1,\frac{1}{n\inf_{\theta\in \Theta}\underline \lambda(\si(\theta))}\I_{|\M|}\right) \ud \theta_1
\end{align*}
where we have used Assumption 7(i) that the prior is bounded, and the fact that $T_{\M}(\theta)\geq 0$. All the determinants and the eigenvalues here are bounded from below and above, by Assumption 9. By choosing $C$ sufficiently large, we can make the integral in the last display arbitrarily small, due to the Gaussian concentration inequality. Hence as we choose $C$ arbitrarily large,
\begin{equation}\label{tnot2}
\int_{\Theta(\M)\backslash B_0(C\e)} \phi(\hat \theta;\theta_1, \si(\theta)^{-1}/n) \pi(\theta_1|\M)\ud \theta_1 = o_p\left(n^{\frac{p-|\M|}{2}}\right).
\end{equation}
Similarly we have
\begin{equation}\label{tnot3}
\int_{\Theta(\M)\backslash B_0(C\e)} \phi(\hat \theta;\theta_1, \si^{-1}/n) \pi(\theta_1|\M)\ud \theta_1 = o_p\left(n^{\frac{p-|\M|}{2}}\right).
\end{equation}
Therefore it is sufficient to show \eqref{tnot1} on $\Theta(\M)\cap B_0(C\e)$. We will use the continuity of $\si(\theta)$ with respect to $\theta$ again, in the sense that $\dett(\si(\theta))/\dett(\si)=1+o_p(1)$ and also $\|\si_{11}(\theta)-\si_{11}\|=o_p(1)$ for $\theta \in B_0(C\e)$. It follows from Lemma \ref {xi1}, Lemma \ref{tm} (i) and the Gaussian concentration inequality that
\begin{align}\label{tnot4}
& \int_{\Theta(\M)\cap B_0(C\e)} \phi(\hat \theta;\theta_1, \si^{-1}(\theta)/n) \pi(\theta_1|\M)\ud \theta_1 \nonumber \\
={}& \int_{\Theta(\M)\cap B_0(C\e)} \left(\frac{2\pi}{n}\right)^{-p/2} \dett(\si(\theta))^{1/2}
e^{-\frac{n}{2}(\theta_1-\bar \theta_{\M,1})^\top\si_{11}(\theta)(\theta_1-\bar \theta_{\M,1}) - T_{\M}(\theta)+o_p(1)} \pi(\theta_1|\M)\ud \theta_1  \nonumber \\
={}& \int_{\Theta(\M)\cap B_0(C\e)} \left(\frac{2\pi}{n}\right)^{-p/2} \dett(\si)^{1/2}(1+o_p(1))\cdot \nonumber \\
& \exp\Big\{-\frac{n}{2}(\theta_1-\bar \theta_{\M,1})^\top\si_{11}(\theta_1-\bar \theta_{\M,1})+o_p\left(\frac{n}{2}\left\|\theta_1-\bar \theta_{\M,1}\right\|^2 \right) \nonumber  \\
& - T_{\M}(\theta_0)+o_p(1)\Big\} \pi(\theta_1|\M)\ud \theta_1  \nonumber \\
={}& (1+o_p(1)) \int_{\Theta(\M)\cap B_0(C\e)} \left(\frac{2\pi}{n}\right) ^{-p/2} \dett(\si)^{1/2}\cdot  \nonumber \\
& \exp\Big\{-\frac{n}{2}(\theta_1-\bar \theta_{\M,1})^\top\si_{11}(\theta_1-\bar \theta_{\M,1})- T_{\M}(\theta_0)\Big\} \pi(\theta_1|\M)\ud \theta_1  \nonumber \\
= {}& (1+o_p(1)) \int_{\Theta(\M)\cap B_0(C\e)} \phi\left(\hat \theta;\theta_1,\si^{-1}/n\right) \pi(\theta_1|\M)\ud \theta_1 \nonumber \\
={}& (1+o_p(1))e^{-T_{\M}(\theta_0)} \left(\frac{2\pi}{n}\right)^{-\frac{p-|\M|}{2}}\left(\frac{\dett(\si)}{\dett(\si_{11})}\right)^{1/2}\pi(\theta_{0,\M,1}|\M).
\end{align}
Therefore, \eqref{tnot1} follows immediately from \eqref{tnot2} - \eqref{tnot4} and the fact that $T_{\M}(\theta_0)=O_p(1)$. We can further combine \eqref{numer1} and \eqref{tnot1} and conclude that
$$\int_{\Theta(\M)} p(\hat\theta| \theta_1)\pi(\theta_1|\M)\ud \theta_1 =\int_{\Theta(\M)} \phi(\hat \theta;\theta_1,\si^{-1}/n) \pi(\theta_1|\M)\ud \theta_1 + o_p \left(n^{\frac{p-|\M|}{2}}\right).$$
Given this result, the Bayes factor based on $p(\hat\theta| \theta_1)$ can be directly translated into the Bayes factor based on $ \phi\left(\hat \theta;\theta_1,\si^{-1}/n\right)$. It follows that
\begin{align}
& \bfh[\M:\M_0] \nonumber \\
={}& \frac{\int_{\Theta(\M)}p(\hat \theta|\theta_1)\pi(\theta_1|\M)\ud \theta_1}{\int_{\Theta(\M_0)}p(\hat\theta| \theta_1)\pi(\theta_1|\M_0)\ud \theta_1} \nonumber \\
={}& \frac{\int_{\Theta(\M)} \phi(\hat \theta;\theta_1,\si^{-1}/n) \pi(\theta_1|\M)\ud \theta_1 + o_p \big(n^{\frac{p-|\M|}{2}}\big)}
{\int_{\Theta(\M_0)} \phi(\hat \theta;\theta_1,\si^{-1}/n) \pi(\theta_1|\M_0)\ud \theta_1 + o_p \big(n^{\frac{p-|\M_0|}{2}}\big)}\nonumber \\
\overset{\eqref{tnot2}-\eqref{tnot4}}{=\joinrel=\joinrel=\joinrel=\joinrel=\joinrel=\joinrel=}&
\frac{ (1+o_p(1))e^{-T_{\M}(\theta_0)} \Big(\frac{2\pi}{n}\Big)^{-\frac{p-|\M|}{2}}\big\{\dett (\G_{\M}^\top\V_n^{-1}\G_{\M})\big\}^{-1/2}\pi(\theta_{0,\M,1}|\M) + o_p \big(n^{\frac{p-|\M|}{2}}\big)}
{ (1+o_p(1))e^{-T_{\M_0}(\theta_0)} \Big(\frac{2\pi}{n}\Big)^{-\frac{p-|\M_0|}{2}}\big\{\dett (\G_{\M_0}^\top\V_n^{-1}\G_{\M_0})\big\}^{-1/2}\pi(\theta_{0,\M,1}|\M_0) + o_p \big(n^{\frac{p-|\M_0|}{2}}\big)} \nonumber \\
\overset{\text{Lemma \ref{tm} (ii)}}{=\joinrel=\joinrel=\joinrel=\joinrel=\joinrel=\joinrel=\joinrel=} {}&  (1+o_p(1))\Big(\frac{2\pi}{n}\Big)^{-\frac{|\M|-|\M_0|}{2}} \frac{S_{\M}(\D)}{S_{\M_0}(\D)} \cdot
\frac{\big\{\dett (\G_{\M}^\top\V_n^{-1}\G_{\M})\big\}^{-1/2}\pi(\theta_0|\M)}
{\big\{\dett (\G_{\M_0}^\top\V_n^{-1}\G_{\M_0})\big\}^{-1/2}\pi(\theta_0|\M_0)}.\nonumber
\end{align}
The last display is exactly the expression of $\BF[\M:\M_0]$ in Lemma \ref{g1}. Also note that uniformly over all $\M\supseteq \M_0$, by Assumptions 7 and 9(i), we have $\pi(\theta_0|\M)/\pi(\theta_{0}|\M_0)=1+o(1)$ and $\pi(\theta_0|\M)\geq \exp(-c_0|\M|)\geq \exp(-c_0\bar p)$ which is a constant lower bound. This is why the terms $o_p \big(n^{\frac{p-|\M|}{2}}\big)$ and $o_p \big(n^{\frac{p-|\M_0|}{2}}\big)$ can be absorbed in the last equality. Therefore, $\frac{\BF[\M:\M_0]}{\bfh[\M:\M_0]}\to 1$ has been proved. Since now $p$ is bounded above by $\bar p$ in Assumption 9(i), one can see that the order of $\bfh[\M:\M_0]$ is equal to $n^{-\frac{|\M|-k_0}{2}}$, which completes the proof.  \hfill$\blacksquare$

\vspace{+1cm}

\noindent {\bf Proof of Theorem 2 (ii):}\\
The first inequality directly follows from Lemma \ref{g2} and Assumption 9(ii). For the second one, we have that for $z=\sqrt{n}\F(\theta)(\hat \theta-\theta)$ with any $\theta \in \Theta(\M)$ and $\M_0\backslash \M\neq \emptyset$,
\begin{align*}
& \|z\|^2 = n(\hat \theta-\theta)^\top \si(\theta) (\hat \theta-\theta)\geq \inf_{\theta\in \Theta} \underline\lambda(\si(\theta)) n\|\hat \theta-\theta\|^2\\
&\geq  n\inf_{\theta\in \Theta} \underline\lambda(\si(\theta)) \big(\|\theta-\hat \theta_1\|^2+\|\hat \theta_2\|^2\big)
\end{align*}
Since  $\M_0\backslash \M\neq \emptyset$, w.p.a.1 as $n\to \infty$, $\hat \theta_2 \to \theta_{0,\M,2} \neq 0$. Furthermore, by Assumption 9(ii),
$$\|\hat \theta_2\|^2=\|\theta_{0,\M,2}\|^2+o_p(1)\geq \frac{1}{2}\|\theta_{0,\M,2}\|^2 \geq \frac{1}{2}\underline \theta^2.$$
Therefore, if we let $C_1=\inf_{\theta\in \Theta} \underline\lambda(\si(\theta))$, then w.p.a.1 as $n\to \infty$,
$$\|z\|^2 \geq C_1 n\big(\|\theta-\hat \theta_1\|^2+\|\hat \theta_2\|^2\big)\geq C_1n\|\theta-\hat \theta_1\|^2 +\frac{C_1\underline \theta^2n}{2}.$$
Thus for some constant $C>0$, we can bound the difference between $p(\hat \theta|\theta_1)$ and the normal limit, by
\begin{align}\label{dif3}
& \int_{\Theta(\M)}\Big|p(\hat \theta|\theta_1)-\phi\left(\hat \theta;\theta_1, \si(\theta)^{-1}/n\right) \Big| \pi(\theta_1|\M)\ud \theta_1 \nonumber \\
\leq {} & n^{p/2}\sup_{\theta\in \Theta}(\dett(\si(\theta)))^{1/2} \cdot \tau_n c_{\pi} \int_{\Theta(\M)}  \frac{1}{1+(\|z\|^2)^{\frac{p+1}{2}}} \ud \theta_1 \nonumber \\
\leq {} & C\tau_n n^{p/2} \int_{\Theta(\M)} \frac{1}{\left\{C_1n\|\theta_1-\hat \theta_1\|^2 +\frac{C_1\underline \theta^2n}{2}\right\}^{\frac{p+1}{2}}} \ud \theta_1 \nonumber \\
={}& C\tau_n n^{p/2} \int_{\Theta(\M)} \frac{(C_1 \underline \theta^2 n /2)^{-\frac{p+1}{2}}}{\left(\frac{\|\theta_1-\hat \theta_1\|^2}{\underline \theta^2/2}+1\right)^{\frac{p+1}{2}}} \ud \theta_1 \nonumber \\
\leq{}& C \tau_n n^{-\frac{1}{2}} \int_{\Theta(\M)}\frac{1}{(1+\|u\|^2)^{\frac{p+1}{2}}} \ud u\nonumber \\
\leq{}& C \tau_n n^{-\frac{1}{2}} \int_0^{\infty} \frac{x^{|\M|-1}}{(1+x^2)^{\frac{p+1}{2}}} \ud x \nonumber \\
\leq {}& C \tau_n n^{-\frac{1}{2}},
\end{align}
where $C$ has absorbed all the constant terms. Using \eqref{dif3}, we can bound the marginal probability $\int_{\Theta(\M)}p(\hat \theta|\theta_1)\pi(\theta_1|\M)\ud \theta_1$ as
\begin{align*}
& \int_{\Theta(\M)}p(\hat \theta|\theta_1)\pi(\theta_1|\M)\ud \theta_1 \\
\leq {}& \int_{\Theta(\M)} \Big| p(\hat\theta_1| \theta_1) - \phi\left(\hat \theta;\theta_1, \si(\theta)^{-1}/n\right) \Big|\pi(\theta_1|\M)\ud \theta_1 +
 \int_{\Theta(\M)} \phi(\hat \theta;\theta_1, \si(\theta)^{-1}/n) \pi(\theta_1|\M)\ud \theta_1 \\
\leq {}& C \tau_n n^{-\frac{1}{2}}  + c_{\pi} \int_{\Theta(\M)} \left(\frac{2\pi}{n}\right)^{-\frac{p}{2}} \left\{\dett(\si(\theta))\right\}^{1/2}
e^{-\frac{C_1 n}{2} (\theta_1-\hat \theta_1)^\top (\theta_1-\hat \theta_1)-\frac{C_1n\underline\theta^2}{4} }\ud \theta_1 \\
\leq {} &  C \tau_n n^{-\frac{1}{2}}  + c_{\pi}\dett(\si(\theta))^{1/2}C_1^{-|\M|/2} e^{-C_1n\underline \theta^2/4} \left(\frac{2\pi}{n}\right)^{-\frac{p-|\M|}{2}} \\
\leq{}& C_1 \tau_n n^{-\frac{1}{2}}  + e^{-C_2n\underline \theta^2}
\end{align*}
for some redefined constants $C_1,C_2>0$. Therefore, w.p.a.1 as $n\to \infty$, the Bayes factor can be bounded by
\begin{align*}
& \bfh[\M:\M_0] \\
={} & \frac{\int_{\Theta(\M)}p(\hat \theta|\theta_1)\pi(\theta_1|\M)\ud \theta_1}{\int_{\Theta(\M_0)}p(\hat\theta| \theta_1)\pi(\theta_1|\M_0)\ud \theta_1}\\
={}& \frac{ C_1 \tau_n n^{-\frac{1}{2}}  + e^{-C_2n\underline \theta^2} }
{(1+o_p(1))e^{-T_{\M_0}(\theta_0)} \left(\frac{2\pi}{n}\right)^{-\frac{p-|\M_0|}{2}}\left\{\dett\left(\G_{\M_0}^\top\V^{-1}\G_{\M_0}\right)\right\}^{-1/2}\pi(\theta_{0,\M,1}|\M_0) + o_p \left(n^{\frac{p-|\M_0|}{2}}\right)}\\
\leq{} & \tau_n n^{\frac{k_0-p-1}{2}} +  e^{-Cn\underline \theta^2}
\end{align*}
for some redefined constant $C>0$ and the first constant can be absorbed into $\tau_n$. This completes the proof of Theorem 2 (ii).  \hfill$\blacksquare$

\vspace{+1cm}

\noindent {\bf Proof of Theorem 3 (i):}\\
If $\M_0\neq \M_{\full}$, then there exists at least one model $\M$ such that $\M\supset \M_0$. Also the total number of models is now bounded by $2^{\bar p}$. Hence under the same prior $\pi(\theta,\M)$, Theorem 2 (i) and Assumption 9(iv) imply that $\sum_{\M:\M\supset \M_0}\poh[\M:\M_0]=O_p\left(n^{-1/2}\right)$, $\sum_{\M:\M\supset \M_0}\PO[\M:\M_0]=O_p\left(n^{-1/2}\right)$, and
$$\sum_{\M:\M\supset \M_0}\poh[\M:\M_0]=\left(1+o_p(1)\right)\sum_{\M:\M\supset \M_0}\PO[\M:\M_0].$$
On the other hand, Theorem 2 (ii) implies that $\sum_{\M:\M_0\backslash \M\neq \emptyset}\PO[\M:\M_0] = \exp(-Cn\underline \theta^2)$ for some constant $C>0$, and also
\begin{align*}
& \sum_{\M:\M_0\backslash \M\neq \emptyset}\poh[\M:\M_0]  \\
\leq {}& 2^{p}\exp\left(-Cn\underline \theta^2\right) \vee \tau_n n^{\frac{k_0-p-1}{2}} \\
\leq {} & \exp\left(-Cn\underline \theta^2\right) \vee \tau_n n^{\frac{k_0-p-1}{2}}
\end{align*}
with adjusted $C$ and $\tau_n=o_p(1)$. Therefore, it is clear that
\begin{align*}
& \sum_{\M:\M_0\backslash \M\neq \emptyset}\PO[\M:\M_0] = o_p\Big( \sum_{\M:\M\supset \M_0}\PO[\M:\M_0] \Big) \\
& \sum_{\M:\M_0\backslash \M\neq \emptyset}\poh[\M:\M_0]  = o_p\Big(\sum_{\M:\M\supset \M_0}\poh[\M:\M_0] \Big)
\end{align*}
Hence it is clear that the posterior consistency follows, with $q(\M_0|\D)=1+o_p(1)$ and $p(\M_0|\hat \theta) = 1+o_p(1)$. Moreover,
\begin{align*}
& q(\M:\M\neq \M_0 |\D)  = \Big(\sum_{\M:\M_0\backslash \M\neq \emptyset}\PO[\M:\M_0] + \sum_{\M:\M\supset \M_0}\PO[\M:\M_0]\Big) q(\M_0|\D) \\
& = (1+o_p(1)) \sum_{\M:\M\supset \M_0}\PO[\M:\M_0] \asymp n^{-1/2} \\
& p(\M:\M\neq \M_0 |\hat \theta)  = \Big(\sum_{\M:\M_0\backslash \M\neq \emptyset}\poh[\M:\M_0] + \sum_{\M:\M\supset \M_0}\poh[\M:\M_0]\Big) p(\M_0|\hat \theta)   \\
& = (1+o_p(1)) \sum_{\M:\M\supset \M_0}\poh[\M:\M_0] \asymp n^{-1/2} \\
& \frac{ q(\M:\M\neq \M_0 |\D)}{p(\M:\M\neq \M_0 |\hat \theta)} = \frac{(1+o_p(1))\sum_{\M:\M\supset \M_0}\PO[\M:\M_0] }{(1+o_p(1)) \sum_{\M:\M\supset \M_0}\poh[\M:\M_0] } \to 1
\end{align*}
w.p.a.1 as $n \to \infty$.\\
When $\M_0=\M_{\full}$, there is no model $\M$ with $\M\supset \M_0$. Hence
\begin{align*}
& q(\M:\M\neq \M_0 |\D)  = \sum_{\M:\M_0\backslash \M\neq \emptyset}\PO[\M:\M_0] \cdot q(\M_0|\D)\leq  \exp\big(-Cn\underline \theta^2\big);\\
& p(\M:\M\neq \M_0 |\hat \theta) = \sum_{\M:\M_0\backslash \M\neq \emptyset}\poh[\M:\M_0] \cdot p(\M_0|\hat \theta)\leq  \exp\big(-Cn\underline \theta^2\big) \vee \tau_n n^{\frac{k_0-p-1}{2}}.
\end{align*}
\hfill$\blacksquare$

\vspace{+.6cm}


\noindent {\bf Proof of Theorem 3 (ii):}\\
Because of the model selection consistency in Theorem 3 (i) for $p(\theta|\hat \theta)$ and the normal approximation in the proof of Theorem 2 (i), one can show that
\begin{align}\label{thm321}
\sup_{A \subseteq \Theta} \Bigg|\int_A p(\theta | \hat \theta)\ud \theta - \int_A  \phi\left(\theta_1; \xi_1(\theta_0), (\G_{\M_0}^\top \V_n^{-1} \G_{\M_0})^{-1}/n\right) \ud \theta_1 \Bigg|\to 0,
\end{align}
where $\theta=(\theta_1^\top,\theta_2^\top)^\top$ is decomposed according to the true model $\M_0$, and $\xi_1(\theta_0)$ is defined in Lemma \ref{matalg} also according to $\M_0$. \eqref{thm321} can be proved using similar arguments to the proof of Theorem 1 (ii), and hence we omit it here. By Lemma \ref{xi1}, $\xi_1(\theta_0) = \bar \theta_{\M_0,1} + o_p\big(1/\sqrt{n}\big)$. Therefore, using the relation between the total variation distance and the Kullback-Leibler (KL) divergence (Pinsker's inequality), we can obtain that
\begin{align}\label{thm322}
& \sup_{A \subseteq \Theta} \Bigg|\int_A \phi\left(\theta_1; \bar \theta_{\M_0,1}, (\G_{\M_0}^\top \V_n^{-1} \G_{\M_0})^{-1}/n\right) \ud \theta_1 - \int_A  \phi\left(\theta_1; \xi_1(\theta_0), (\G_{\M_0}^\top \V_n^{-1} \G_{\M_0})^{-1}/n\right) \ud \theta_1 \Bigg| \nonumber \\
&\leq \left\{\frac{1}{2}\text{KL}\left(\phi\left(\theta_1; \bar \theta_{\M_0,1}, (\G_{\M_0}^\top \V_n^{-1} \G_{\M_0})^{-1}/n\right), \phi\left(\theta_1; \xi_1(\theta_0), (\G_{\M_0}^\top \V_n^{-1} \G_{\M_0})^{-1}/n\right) \right)\right\}^{\frac{1}{2}} \nonumber \\
&= \frac{1}{2} \sqrt{n\left(\bar \theta_{\M_0,1}- \xi_1(\theta_0)\right)^\top \G_{\M_0}^\top \V_n^{-1} \G_{\M_0}\left(\bar \theta_{\M_0,1}- \xi_1(\theta_0)\right)} \nonumber \\
&\leq \frac{1}{2}\overline \lambda \left(\G_{\M_0}^\top \V_n^{-1} \G_{\M_0}\right) \cdot n \left\|\bar \theta_{\M_0,1}- \xi_1(\theta_0)\right\|^2 = o_p(1),
\end{align}
where $\text{KL}(f_1,f_2)$ denote the KL divergence between two densities $f_1$ and $f_2$.  Based on Assumption 1-11, the result of Theorem 3 (ii) immediately follows by combining Theorem 1 (ii), \eqref{thm321} and \eqref{thm322}. \hfill$\blacksquare$

\vspace{+1cm} \noindent {\bf
\Large 3. Application to Motivating Examples}
\vspace{+.5cm}

In this section, we prove the Bayesian oracle property for the three motivating examples in Section 1.3 of the main paper. \\

\noindent {\bf \large Example 1: Correlated Longitudinal Data}\\

We use the same notations as in the introduction. Without loss of generality, we assume $Y_i$ and each $X_{ijk}$ has been centered
such that $\Ep Y_i=0$ and $\Ep X_{ijk}=0$, for $i=1,\ldots,n$, $j=1,\ldots,s$ and $k=1,\ldots,p$.
For the ease of presentation and the simplification of our proofs, we assume that the working correlation matrix $\R$ is correctly specified and does not depend on $\theta$. We also plug in a preliminary consistent estimator $\tilde \theta$ for $\theta$ to the nonlinear part $\frac{\partial \mu_i(\theta)^\top}{\partial \theta} \S_i^{-1}$ of the moment function $g(D_i,\theta)$. Such consistent estimator $\tilde \theta$ exists even with growing $p$. For example, one can take $\tilde \theta$ to be the solution of estimating equations $\sum_{i=1}^n \X_i(Y_i-\mu_i(\theta))=0$. Under the assumptions given in the theorem below, one can show that $\|\tilde \theta - \theta_0\|=O_p(\sqrt{p/n})$ similar to Example 1 in \citet{wang11}.\\

The matrix $\V_n$ in BGMM can be taken as $\V_n=n^{-1}\sum_{i=1}^n (g(D_i,\tilde \theta)-\bar g(\D,\tilde \theta))(g(D_i,\tilde \theta)-\bar g(\D,\tilde \theta))^\top$, where $\tilde \theta$ is any preliminary consistent estimator of $\theta_0$. One can show that $\V_n$ converges in probability under the operator norm to $\V=\var(g(D,\theta_0))$.\\



Let $\dot \mu(x)$ and $\ddot \mu(x)$ be the first and the second derivatives of $\mu(x)$. We then have the following theorem for the BGMM based on the moment function (5) for the correlated longitudinal data.\\

\noindent {\bf Theorem S1} {\it For the moment function (5) in the main paper, suppose that Assumptions 1, 2, 3, 7 and 8 hold. Suppose $\tilde\theta$ is the preliminary estimator that solves $\sum_{i=1}^n \X_i(Y_i-\mu_i(\theta))=0$ and $\V_n=n^{-1}\sum_{i=1}^n (g(D_i,\tilde \theta)-\bar g(\D,\tilde \theta))(g(D_i,\tilde \theta)-\bar g(\D,\tilde \theta))^\top$. In addition, if \\
\noindent (1) $|X_{ijk}|\leq C_X$ almost surely for some large constant $C_X>0$ and all $i=1,\ldots,n$, $j=1,\ldots,s$, $k=1,\ldots,p$. $\sup_{1\leq j \leq s}\Ep (Y_j^4) < \infty$;\\
\noindent (2) $\Ep (\X_{i}^\top\X_{i})$ and $\R$ have eigenvalues bounded above and below by constants for all $i=1,\ldots,n$;\\
\noindent (3) $\dot \mu(X_{ij}^\top \theta)$ is bounded above and below uniformly for all possible values of $X_{ij}$ and $\theta \in \Theta$.  $\ddot \mu(X_{ij}^\top \theta)$ is bounded above, and $\phi(X_{ij}^\top \theta)$ is bounded above and below uniformly for all $X_{ij}$ and $\theta \in B_0(c\e)$ for any fixed $c>0$;\\
then Assumptions 4, 5 and 6 hold with $\alpha=1$. Therefore BGMM for the moment function (5) in the main paper satisfies the Bayesian oracle property in Theorem 1.}

\begin{remark}
In Condition (1) we impose an absolute bound on all the covariates for convenience, though this can be replaced by relaxed conditions on the tail behavior or the high order moments on $X_{ijk}$. Condition (2) for eigenvalues are standard.
Here for simplicity, we use only one working correlation matrix $\R$ such that $m=p$, though the result can be
easily extended to more than one working correlation matrices like in \citet{qu00}. Condition (3) requires certain bounds on the derivatives of $\mu$ and also $\phi$. In particular, $\mu(t)=t$ for linear regression trivially satisfies this condition. For logistic regression, $\mu(t)=e^t/(1+e^t)$ and $\phi(t)=e^t/(1+e^t)^2$. Since in our Assumption 2 the dimension $p$ is fixed, condition (3) is satisfied for the derivatives of $\mu$ and $\phi$ evaluated at $X_{ij}^\top \theta$. Similar arguments can be applied to Poisson regression, exponential regression and probit regression, etc. In \citet{liang86} and \citet{wang12}, the marginal density of $Y_{ij}$ is modeled as a canonical exponential family, with $\var(Y_{ij}|X_{ij})=\psi \dot \mu (X_{ij}^\top \theta) $, where $\psi$ is the dispersion parameter. Here we have considered a general form of the function $\phi$ and therefore our setup includes theirs as a special case.
\end{remark}

\noindent \textbf{Proof of Theorem S1:}\\
In the following, for a generic random variable $D=(Y,\X)^\top$ (independent of the sample $\D$), we omit the subscript $i$ in $X_{ij}$ and write $X_{\cdot j}$ to represent a generic $p$-dimensional covariate vector measured at time $j$, for $j=1,\ldots,s$. Define $\mu_{\cdot}(\theta) = (\mu(X_{\cdot 1}^\top\theta),\ldots,\mu(X_{\cdot s}^\top\theta))^\top$, $\B (\theta)=\frac{\partial \mu_{\cdot} (\theta)}{\partial \theta}=(\dot \mu(X_{\cdot 1}^\top  \theta)X_{\cdot 1}, \ldots, \dot \mu(X_{\cdot s}^\top \theta)X_{\cdot s})^\top$, and $\S( \theta)=\A( \theta )^{1/2} \R \A(\theta )^{1/2}$ (if $X$ has sample index $i$, then we use the notation $\S_i(\theta)$). So the generic moment function can be written as $g(D,\theta)= \B(\tilde \theta)^\top \S(\tilde \theta)^{-1} (Y-\mu_\cdot(\theta))$, where $\tilde \theta$ is a preliminary estimator that solves $\sum_{i=1}^n \X_i(Y_i-\mu_i(\theta))=0$. Similar to the proof of (3.3) in \citet{wang11}, one can show that given the conditions (1)-(3), $\|\tilde \theta - \theta_0\| = O_p\left(\sqrt{p/n}\right)$. For simplicity we omit the proof of this relation here.

Suppose that the constant upper and lower bounds for $\phi(X_{ij}^\top\theta_0)$ in Condition (3) are $\bar \phi$ and $\underline \phi$ respectively. Since w.p.a.1 as $n\to \infty$, $\|\tilde \theta - \theta_0\| \leq C\e$, $\S_i(\tilde\theta)=\A_i(\tilde\theta)^{1/2}\R\A_i(\tilde\theta)^{1/2}$ and $\A_i(\theta)=\diag \big\{\phi(X_{i1}^\top\theta),\ldots, \phi(X_{is}^\top\theta)\big\}$, we know that w.p.a.1 as $n\to \infty$, the eigenvalues of $\S_i(\tilde\theta)$ can be bounded as
\begin{align*}
& \bar \lambda (\S_i(\tilde\theta)) \leq \bar \lambda (\R) \bar \lambda (\A_i(\tilde\theta)) \leq \bar \lambda (\R) \tr(\A_i(\tilde\theta))\leq s\bar \phi\bar \lambda (\R)\\
& \underline \lambda (\S_i(\tilde\theta)) \geq \underline \lambda (\R) \underline \lambda (\A_i(\tilde\theta)) \geq  \underline \phi\underline \lambda (\R),
\end{align*}
where the upper and the lower bounds are constants that do not change with $n$.

We now check Assumptions 4 and 5. Let $\bar {\dot \mu}$ and $\underline {\dot \mu}$ be the constant upper and lower bounds for $\dot \mu(X_{ij}^\top\theta)$ in the condition (2). Then for Assumption 4(i), using the boundedness of $X_{ijk}$ in the condition (1), we have that w.p.a.1 as $n\to \infty$,
\begin{align*}
& \sup_{\|\eta\|=1 } \Ep\left\{(\eta^\top (g(D,\theta)-g(D,\theta_0)))^2\right\} \leq  \sup_{\|\eta\|=1 } \Ep\left\{(\eta^\top \bar {\dot \mu}\X ^\top \S(\tilde \theta)^{-1}\bar {\dot \mu}\X (\theta-\theta_0))^2\right\}\\
\leq {} & \bar {\dot \mu}^2 \underline \lambda(\S(\tilde \theta))^{-2} \sup_{\|\eta\|=1 } \eta^\top \Ep \left\{\X^\top \X (\theta-\theta_0) (\theta-\theta_0)^\top \X^\top \X\right\} \eta  \\
\leq {} & \bar {\dot \mu}^2 \underline \lambda(\S(\tilde \theta))^{-2}\bar \lambda (E(\X^\top\X)) \tr(\X^ \top \X)\|\theta-\theta_0\|^2 \\
\leq {} & \bar {\dot \mu}^2 \underline \lambda(\S(\tilde \theta))^{-2}\bar \lambda (E(\X^\top\X))\cdot spC_X^2\|\theta-\theta_0\|^2 = O\left((p^{1/2}\|\theta-\theta_0\|)^2\right)
\end{align*}
Therefore this implies that in Assumption 4(i) we can take $\alpha=1$, and also in Assumption 4(ii), the $L_2$ norm of the envelope function $F$ for the class $\mathcal{F}$ is of order $O(\sqrt{p})$, since the $L_2$ radius of $\Theta$ is assumed to be bounded by constant $R$ in Assumption 1. Next we estimate the $L_2$ uniform covering number of $\mathcal{F}=\{f(\eta,\theta)=\eta^\top \B(\tilde \theta)^\top\S(\tilde \theta)^{-1} (\mu_{\cdot}(\theta_0)- \mu_{\cdot}(\theta)),\theta\in \Theta, \eta\in \mathbb{R}^m, \|\eta\|=1\}$. Suppose there exists a $\epsilon$-net in $L_2(P_D)$ norm for $\mathcal{F}$:$\big\{(\eta_1,\theta_1),\ldots,(\eta_N,\theta_N)\big\}$, with $N=N(\epsilon\|F\|_{P_D,2},\mathcal{F},L_2(P_D))$. Then by definition, for any $(\eta,\theta)$, one can pick out a pair $(\eta_k,\theta_k)$, for some $1\leq k \leq N$, such that $\Ep\left|f(\eta_k,\theta_k)-f(\eta,\theta)\right|^2\leq \epsilon^2$. Then since
\begin{align*}
& \Ep\left|f(\eta_k,\theta_k)-f(\eta,\theta)\right|^2 \\
\leq {}& 2\Ep \left\{(\eta_k-\eta)^\top \B(\tilde \theta)^\top\S(\tilde \theta)^{-1} (\mu_{\cdot}(\theta_k)-\mu_{\cdot}(\theta_0))\right\}^2 +
2\Ep\left\{\eta^\top \B(\tilde \theta)^\top\S(\tilde \theta)^{-1}  (\mu_{\cdot}(\theta)-\mu_{\cdot}(\theta_k))\right\}^2 \\
\leq {}& 2 \bar {\dot \mu}^2 \underline \lambda(\S(\tilde \theta))^{-2} \cdot \bar \lambda(E(\X^\top\X)) \cdot spC_X^2\bar {\dot \mu} ^2 \cdot 4R^2\|\eta_k-\eta\|^2  \\
& + 2 \|\eta\|^2 \bar {\dot \mu}^2 \underline \lambda(\S(\tilde \theta))^{-2} \cdot \bar \lambda(E(\X^\top\X)) \cdot spC_X^2\bar {\dot \mu} ^2 \|\theta_k-\theta\|^2\\
\leq {} & \left(C_1p^{1/2}\|\eta_k-\eta\|\right)^2 + \left(C_2 p^{1/2} \|\theta_k-\theta\|\right)^2
\end{align*}
for some constants $C_1,C_2>0$ that depend on the eigenvalues, $R$, and $\bar{ \dot\mu}$. Thus we only need $\|\eta_k-\eta\|\leq \epsilon/(2C_1p^{1/2})$ and $\|\theta_k-\theta\|\leq \epsilon/(2C_2p^{1/2})$. Since $\|\eta\|=1$ and $\|\eta_k-\eta\|\leq p^{1/2}|\eta_k-\eta|_{\infty}$, we estimate the covering number on $\eta$ using $L_{\infty}$ grids and need no more than
$N_\eta=\left(\frac{2C_1p}{\epsilon}+1\right)^p$ points. Similarly since $\|\theta\|\leq R$, we need no more than $N_\theta = \left(\frac{2C_2Rp}{\epsilon}+1\right)^p$ points. Together we have shown that for small $\epsilon>0$,
$$N\left(\epsilon\|F\|_{P_D,2},\mathcal{F},L_2(P_D)\right) \leq N_\eta N_\theta \leq \left(\frac{9C_1C_2p^2R}{\epsilon^2}\right)^p,$$
which give $\ln N\left(\epsilon\|F\|_{P_D,2},\mathcal{F},L_2(P_D)\right) = O\left(p\ln(n/\epsilon)\right)$. So Assumption 4(ii) holds.

For Assumption 5(i), we have
$$ \|\Ep g(D,\theta)\|=\left\|\Ep \left(\B(\tilde \theta)^\top\S(\tilde \theta)^{-1}  (\mu_{\cdot}(\theta) - \mu_{\cdot}(\theta_0))\right)\right\|
\geq \underline {\dot \mu} ^2 \bar \lambda(\S(\tilde \theta))^{-1} \underline \lambda \left(\Ep (\X^\top \X)\right)\|\theta-\theta_0\|.$$
Therefore Assumption 5(i) holds with $\delta_1=\underline {\dot \mu} ^2 \bar \lambda(\S(\tilde \theta))^{-1} \underline \lambda (\Ep(\X^\top\X))$ and $\delta_0=R\delta_1$.

For Assumption 5(ii), $\G=\nabla _{\theta} \Ep g(D,\theta_0) = - \Ep \left\{\B(\tilde \theta)^\top\S(\tilde \theta)^{-1}\B(\theta_0)\right\}$. By conditions (2) and (3),
\begin{align*}
\bar \lambda (\G^\top\G) & \leq \bar {\dot \mu}^4 \underline \lambda (\S(\tilde \theta))^{-2} \bar \lambda\left(\Ep(\X^\top\X)\right)^2\\
\underline \lambda (\G^\top\G) & \leq \underline {\dot \mu}^4 \bar \lambda (\S(\tilde \theta))^{-2} \underline \lambda\left(\Ep (\X^\top\X)\right)^2
\end{align*}
so the eigenvalues of $\G^\top\G$ are bounded above and below as $n\to \infty$.

In Assumption 5(iii), let $\K(\theta)=\left(\ddot \mu(X_{\cdot 1}^\top\theta)(X_{\cdot 1}\otimes X_{\cdot 1}) ,\ldots,\ddot \mu(X_{\cdot s}^\top\theta)(X_{\cdot s}\otimes X_{\cdot s}) \right)^\top$. Then for any unit vectors $u,v\in \mathbb{R}^p$,
\begin{align*}
\|\H(\theta)(u,v)\| & = \left\|\Ep\left\{\B(\tilde \theta)^\top\S(\tilde \theta)^{-1} \K(\theta)\right\}\cdot \vecc(u\otimes v)\right\| \\
&\leq \bar {\dot \mu}\bar {\ddot \mu} \underline \lambda(\S(\tilde \theta))^{-1} \bar \lambda \left(\Ep(\X^\top\X)\right) \cdot \sqrt{p}C_X,
\end{align*}
where we used the upper bound on $\ddot \mu(X_{ij}^\top\theta)$ for any $\theta \in B_0(c\e)$ in condition (3). Therefore Assumption 5(iii) holds. \\

To show Assumption 6, we note that since $\V_n$ and $\V$ are symmetric positive definite matrices, $\|\V\|=\sqrt{\bar \lambda(\V^\top\V)}=\bar \lambda (\V)$ and also $\|\V_n\|=\bar \lambda (\V_n)$. If we can show $\|\V_n-\V\|\to 0$ w.p.a.1 as $n\to \infty$ and the eigenvalues of $\V$ are bounded from above and below, then we have
\begin{align*}
&\bar \lambda(\V_n) = \|\V_n\|\leq \|\V_n-\V\| + \|\V\| = \|\V_n-\V\| + \bar \lambda(\V) \\
&\underline\lambda(\V_n) = \min_{\eta\in \mathbb{R}^p} \eta^\top \V_n\eta \geq \min_{\eta\in \mathbb{R}^p} \eta^\top \V \eta -
\max_{\eta\in \mathbb{R}^p} \eta^\top (\V-\V_n) \eta \\
& \geq \underline \lambda(\V) - \|\V_n-\V\|.
\end{align*}
Therefore w.p.a.1 as $n\to \infty$, the eigenvalues of $\V_n$ are also bounded from above and below, as long as $\|\V_n-\V\|\to 0$. Next we show the boundedness for the eigenvalues of $\V$ and the convergence of $\|\V_n-\V\|$, respectively.

Since $\Ep g(D,\theta_0)=0$, we have
\begin{align*}
& \V= \var (g(D,\theta_0)) = \Ep\left\{\B(\theta_0)^\top\S(\theta_0)^{-1}(Y-\mu(\theta_0))(Y-\mu(\theta_0))^\top \S(\theta_0)^{-1}\B(\theta_0)\right\}\\
& = \Ep \left(\B(\theta_0)^\top\S(\theta_0)^{-1} \S(\theta_0) \S(\theta_0)^{-1}\B(\theta_0)\right)\\
& = \Ep \left(\B(\theta_0)^\top \S(\theta_0)^{-1}\B(\theta_0)\right).
\end{align*}
By a similar argument to the boundedness of eigenvalues of $\S(\tilde\theta)$, one can show that the eigenvalues of $\S(\theta_0)$ are also bounded from above and below by constants. Therefore,
\begin{align*}
\bar \lambda(\V) &\leq \bar \lambda\left(\Ep\left(\B(\theta_0)^\top\B(\theta_0)\right)\right) \underline \lambda\left(\S(\theta_0)\right)^{-1} \leq \bar {\dot \mu}^2 \bar \lambda\left(\Ep(\X^\top\X)\right) \underline \lambda(\S(\theta_0))^{-1} \\
 \underline \lambda(\V) &\geq \underline \lambda\left(\Ep\left(\B(\theta_0)^\top\B(\theta_0)\right)\right) \bar \lambda\left(\S(\theta_0)\right)^{-1}\geq  \underline {\dot \mu}^2 \underline \lambda\left(\Ep(\X^\top\X)\right) \bar \lambda\left(\S(\theta_0)\right)^{-1}.
\end{align*}
The boundedness of $\bar \lambda(\V)$ and $\underline \lambda(\V)$ is proved.

To show $\|\V_n-\V\|\to 0$, we first note that
\begin{align}\label{vvvv}
\V_n & = \frac{1}{n} \sum_{i=1}^n \big(g(D_i,\tilde\theta)-\bar g(\D,\tilde \theta)\big)\big(g(D_i,\tilde\theta)-\bar g(\D,\tilde \theta)\big)^\top\nonumber \\
&= \frac{1}{n} \sum_{i=1}^n g(D_i,\tilde\theta)g(D_i,\tilde\theta)^\top + \frac{1}{n} \sum_{i=1}^n \bar g(\D,\tilde\theta) \bar g(\D,\tilde\theta)^\top\nonumber\\
& =  \frac{1}{n} \sum_{i=1}^n g(D_i,\theta_0)g(D_i,\theta_0)^\top + \frac{2}{n} \sum_{i=1}^n \big(g(D_i,\tilde\theta)- g(D_i,\theta_0)\big)g(D_i,\theta_0)^\top\nonumber \\
& +  \frac{1}{n} \sum_{i=1}^n \big(g(D_i,\tilde\theta)- g(D_i,\theta_0)\big)\big(g(D_i,\tilde\theta)- g(D_i,\theta_0)\big)^\top  + \bar g(\D,\tilde\theta) \bar g(\D,\tilde\theta)^\top\nonumber\\
& := \E_1+\E_2+\E_3+\E_4.
\end{align}
We derive bounds for each term. For $\E_1$ we have
\begin{align*}
& \|\E_1-\V\|^2\leq \|\E_1-\V\|^2_F \\
=&{} \sum_{j=1}^p\sum_{k=1}^p \Big[n^{-1}\sum_{i=1}^n g_j(D_i,\theta_0)g_k(D_i,\theta_0) - \Ep[g_j(D,\theta_0)g_k(D,\theta_0)]\Big]^2
\end{align*}
Hence by Chebyshev's inequality, for any $C>0$,
\begin{align*}
&P\big(\|\E_1-\V\|^2 >C\big) \\
&\leq C^{-2}\sum_{j=1}^p\sum_{k=1}^p \Ep\Big[n^{-1}\sum_{i=1}^n g_j(D_i,\theta_0)g_k(D_i,\theta_0) - \Ep [g_j(D,\theta_0)g_k(D,\theta_0)]\Big]^2\\
&=\frac{1}{nC^2} \sum_{j=1}^p\sum_{k=1}^p \var(g_j(D,\theta_0)g_k(D,\theta_0)) \\
&\leq  \frac{p^2}{nC^2} \sup_{1\leq j,k\leq p} \Ep[g_j(D,\theta_0)^2g_k(D,\theta_0)^2] \leq \frac  {p^2}{nC^2} \sup_{1\leq j\leq p} \Ep[g_j(D,\theta_0)^4]\\
&\leq \frac  {p^2}{nC^2} \sup_{1\leq j\leq p} \bar{\dot \mu}^4 \Ep[X_{j\cdot}^\top\S(\theta_0)^{-1}(Y-\mu(\theta_0))]^4\\
&\leq \frac  {p^2}{nC^2} \bar{\dot \mu}^4 s^2C_X^4 \underline \lambda (\S(\theta_0))^{-4} \Ep\|Y-\mu(\theta_0)\|^4\\
&= \frac  {p^2}{nC^2} \bar{\dot \mu}^4 s^2C_X^4 \underline \lambda (\S(\theta_0))^{-4} \Ep\big[\sum_{j=1}^s (Y_j-\mu_j(\theta_0))^2\big]^2\\
&\leq \frac  {p^2}{nC^2} \bar{\dot \mu}^4 s^3 C_X^4 \underline \lambda (\S(\theta_0))^{-4} \Ep\big[\sum_{j=1}^s (Y_j-\mu_j(\theta_0))^4\big]\\
&\leq \frac  {8p^2}{nC^2} \bar{\dot \mu}^4 s^3 C_X^4 \underline \lambda (\S(\theta_0))^{-4} \Ep\sum_{j=1}^s \big[ Y_j^4+\mu_j(\theta_0)^4\big]\\
&\leq \frac  {16p^2}{nC^2} \bar{\dot \mu}^4 s^4 C_X^4 \underline \lambda (\S(\theta_0))^{-4} \sup_{1\leq j\leq s} \Ep(Y_j^4).
\end{align*}
Since $ \sup_{1\leq j\leq s} \Ep(Y_j^4)<\infty$ as in Condition (1), we conclude that $\|\E_1-\V\|=O_p(p/\sqrt{n})=o_p(1)$.

Next we bound $\E_3$. Because $\|\tilde\theta-\theta_0\|=O_p(\sqrt{p/n})$, we have that for any generic $D$,
\begin{align}\label{gtg0}
&\|g(D,\tilde \theta)-g(D,\theta_0)\| = \left\|\B(\tilde \theta)^\top\S(\tilde\theta)^{-1} (\mu(\tilde \theta)-\mu(\theta_0))\right\|\nonumber \\
&= \left\|\B(\tilde \theta)^\top\S(\tilde\theta)^{-1}\B(\theta')(\tilde \theta-\theta_0)\right\|\leq
\left\|\B(\tilde \theta)^\top\S(\tilde\theta)^{-1}\B(\theta')\right\|\cdot\left\|\tilde \theta-\theta_0 \right\|\nonumber \\
& \leq \bar{\dot \mu}^2 \underline\lambda (\S(\tilde\theta))\cdot  \|\X^\top\X\|\cdot \left\|\tilde\theta-\theta_0\right\|\leq  \bar{\dot \mu}^2 \underline\lambda (\S(\tilde\theta))\cdot \|\X\|^2 \cdot O_p\left(\sqrt{\frac{p}{n}}\right)\nonumber \\
&\leq  \bar{\dot \mu}^2 \underline\lambda (\S(\tilde\theta))\cdot spC_X^2 \cdot O_p\left(\sqrt{\frac{p}{n}}\right) = O_p\left(\sqrt{\frac{p^3}{n}}\right)=o_p(1),
\end{align}
where $\theta'$ is between $\tilde\theta$ and $\theta_0$, and the derivation shows that if we replace $\D$ with $\D_i$, then the upper bound is uniform over all $i=1,\ldots,n$. Therefore
\begin{align*}
&\|\E_3\|\leq n^{-1} \sum_{i=1}^n  \Big\|\big(g(D_i,\tilde\theta)- g(D_i,\theta_0)\big)\big(g(D_i,\tilde\theta)- g(D_i,\theta_0)\big)^\top\Big\| \\
&\leq  n^{-1} \sum_{i=1}^n \Big\|g(D_i,\tilde\theta)- g(D_i,\theta_0)\Big\|^2=O_p\left(\frac{p^3}{n}\right).
\end{align*}
Given the bounds for $\E_1$ and $\E_3$, we can bound $\E_2$ as
\begin{align*}
&\|\E_2\| = 2n^{-1} \Big\| \sum_{i=1}^n\big(g(D_i,\tilde\theta)- g(D_i,\theta_0)\big)g(D_i,\theta_0)^\top \Big\|\\
& \leq  2n^{-1} \sum_{i=1}^n \Big\|g(D_i,\tilde\theta)- g(D_i,\theta_0)\Big\|\cdot \Big\|g(D_i,\theta_0)\Big\|\\
&\leq O_p\left(\sqrt{\frac{p^3}{n}}\right)\cdot 2n^{-1}\sum_{i=1}^n \Big\|g(D_i,\theta_0)\Big\| \\
&\leq O_p\left(\sqrt{\frac{p^3}{n}}\right)\cdot O_p\left(\sqrt{\frac{p}{n}}\right) = O_p\left(\frac{p^2}{n}\right)=o_p(1).
\end{align*}
For $\E_4$, we use $\|\bar g(\D,\theta_0)\|=O_p\left(\sqrt{p/n}\right)$ and \eqref{gtg0}
\begin{align*}
&\|\E_4\|= \left\|\bar g(\D,\tilde\theta)\right\|^2 \leq \big(\|\bar g(\D,\tilde\theta)-\bar g(\D,\theta_0)\|+\|\bar g(\D,\theta_0)\| \big)^2\\
& \leq\left(O_p(\sqrt{p^3/n}) +O_p(\sqrt{p/n})\right)^2 =O_p\left(p^3/n\right)=o_p(1).
\end{align*}
Finally, we combine the bounds for $\E_1,\E_2,\E_3,\E_4$ and conclude that $\|\V_n-\V\|=o_p(1)$. Therefore Assumption 6 holds.
\hfill$\blacksquare$

\vspace{+1cm}

\newpage
\noindent {\bf \large Example 2: Quantile Regression}\\

Without loss of generality, we assume that the random variables $Y$ and $X$ are centered such that $\Ep Y=0$ and $\Ep X=0$. The conditional distribution $F_{Y|X}$ is assumed to be continuous, and let $f_{Y|X}$ be its conditional density. It can be calculated that $\V=\var(g(D,\theta_0))=\tau(1-\tau)\Ep(XX^\top)$ for the unconditional moments (6), and we can estimate
$\V$ by $\V_n=n^{-1}\tau(1-\tau)\sum_{i=1}^n X_iX_i^\top$. Then we have the following theorem about quantile regression.\\



\noindent {\bf Theorem S2} {\it For the moment function (6) in the main paper, suppose that Assumptions 1, 2, 3, 7 and 8 hold. Suppose that $\V_n=n^{-1}\tau(1-\tau)\sum_{i=1}^n X_iX_i^\top$. In addition, if\\
\noindent (1) For any generic random vector $X=(X_1,\ldots,X_p)^\top$, $|X_j|\leq C_X$ almost surely for some large constants $C_X>0$ and all $j=1,\ldots,p$;\\
\noindent (2) $f_{Y|X}$ is continuously differentiable with the first derivative $\dot f_{Y|X}$. $f_{Y|X}$ and $\dot f_{Y|X}$ are almost surely bounded above on the support of $Y$ for any value of $X$. $f_{Y|X}$ is further bounded below for any value of $X$.\\
\noindent (3) $\Ep(XX^\top)$ has eigenvalues bounded above and below by constants.\\
then Assumptions 4, 5 and 6 hold with $\alpha=1/2$. The BGMM for the moment function (6) in the main paper satisfies the Bayesian oracle property in Theorem 1.
}
\begin{remark}
The quantile regression example can be generalized to the instrumental variable quantile regression model (IVQR), as discussed in \citet{chh05,chh06}. In the IVQR, the predictor $X$ could contain endogenous components, and we can still consistently estimate the parameter $\theta$ using other informative and exogenous instrumental variables. The model formulation will be more complicated but can be incorporated into the BGMM framework using the unconditional moments based on IV (e.g. \citet{ch03}).
\end{remark}

\noindent \textbf{Proof of Theorem S2:}\\
We check Assumptions 4 and 5. For a generic $\theta$, let $A=\{Y \text{ is between } X^\top\theta \text{ and } X^\top\theta_0\}$. Let  $\bar f$, $\bar {\dot f}$ and $\underline f$ be the upper bounds for the conditional density $f_{Y|X}$, its derivative $\dot f_{Y|X}$ and the lower bound for $f_{Y|X}$ in condition (2). Then
\begin{align*}
& \sup_{\|\eta\|=1}\Ep\left\{\eta^\top (g(D,\theta)-g(D,\theta_0))\right\}^2 = \sup_{\|\eta\|=1}\Ep\left\{\eta^\top X(1(Y\leq X^\top\theta)-1(Y\leq X^\top\theta_0))\right\}^2 \\
= {} & \sup_{\|\eta\|=1}\eta^\top {\Ep}_{X} \left\{X\left(\int_A f_{Y|X}(y)\ud y\right)X^\top \right\} \eta\\
\leq {} & \bar \lambda \left(E(XX^\top)\right)\bar f C_X p^{1/2} \|\theta-\theta_0\|.
\end{align*}
Therefore Assumption 4(i) follows by taking $\alpha=1/2$, since the eigenvalues of $\Ep(XX^\top)$ are bounded by condition (3). It also implies that the $L_2$ norm for the envelope function $F$ of the class $\mathcal{F}$ in Assumption 4(ii) is bounded by $O(R^{1/2} p^{1/4})\leq O(p^{1/2})$. Moreover, the VC index of the class $\mathcal {F}$ is of order $O(p)$ (see Lemma 18-20 of \citealt{bci11}), and the bound on the uniform covering number follows by Theorem 2.6.7 of \citet{vw96}.\\
For Assumption 5(i), we have
\begin{align*}
& \left\|\Ep g(D,\theta)\right\|^2=\left\|\Ep\left\{X(1(Y\leq X^\top \theta)-\tau)\right\}\right\|^2 = \left\|\Ep\left\{X(F_{Y|X}(X^\top\theta)-F_{Y|X}(X^\top\theta_0))]\right\}\right\|^2 \\
& = \left\|\Ep\left\{XX^\top f_{Y|X}(X^\top\tilde \theta)\cdot (\theta-\theta_0)\right\}\right\|^2 \geq \underline f^2 \underline \lambda\left(\Ep(XX^\top)\right) \|\theta-\theta_0\|^2,
\end{align*}
where in the second equality we used the iterated expectation, in the third equality $\tilde \theta$ is between $\theta$ and $\theta_0$. This implies that $\|\Ep g(D,\theta)\| \geq \delta_1\|\theta-\theta_0\|$, with $\delta_1 =\underline f^2 \underline \lambda(\Ep(XX^\top)) $. Therefore we can simply take $\delta_0=2R\delta_1$, and Assumption 5(i) holds. \\
For Assumption 5(ii), one can calculate that $\G=\Ep\left\{XX^\top f_{Y|X}(X^\top \theta_0)\right\}$. Using the definition of the matrix operator norm, one can see that the eigenvalues of $\G^\top\G$ can be bounded as
\begin{align*}
& \bar \lambda (\G^\top\G) \leq \bar f^2 \bar \lambda \left(E(XX^\top)\right)\\
& \underline \lambda (\G^\top\G) \geq \underline f^2 \underline \lambda \left(E(XX^\top)\right)
\end{align*}
For Assumption 5(iii), for any unit vectors $u,v\in \mathbb{R}^p$,
\begin{align*}
\|\H(\theta)(u,v)\| & = \left\|\Ep\left\{XX^\top \otimes X^\top \dot f_{Y|X}(X^\top\theta)\right\}\cdot \vecc(u\otimes v)\right\| \\
& \leq \bar \lambda\left(\Ep(XX^\top)\right)\cdot \sqrt{p}C_X\bar {\dot f}.
\end{align*}
Hence Assumption 5(iii) holds.

For Assumption 6, by Chebyshev's inequality, for any $C>0$, we have
\begin{align*}
& \Pr\big(\|\V_n-\V\|\geq C\big)\leq \Pr\big(\|\V_n-\V\|_F\geq C\big) \\
\leq {}& \frac{\tau^2(1-\tau)^2}{C^2} \var \Bigg\{\sum_{j=1}^p\sum_{k=1}^p \Big(\frac{1}{n}\sum_{i=1}^nX_{ij}X_{ik}-\Ep(X_{ij}X_{ik}) \Big)^2\Bigg\}\\
\leq{} & \frac{\tau^2(1-\tau)^2p^2}{nC^2} \sup_{1\leq j,k\leq p} \Ep(X^2_{j}X^2_{k}) \leq \frac{p^2C_X^4\tau^2(1-\tau)^2}{nC^2}.
\end{align*}
Therefore $\|\V_n-\V\|=O_p\left(p/\sqrt{n}\right)=o_p(1)$. The boundedness of eigenvalues of $\V$ follows directly from the boundedness of eigenvalues of $\Ep(XX^\top)$ in Condition (3), and hence the eigenvalues of $\V_n$ are also bounded from above and below w.p.a.1 as $n\to \infty$.
\hfill$\blacksquare$

\vspace{+1cm}
\newpage

\noindent {\bf \large Example 3: Partial Correlation Selection}\\

We use the same notation as in the introduction. Suppose the true covariance matrix is $\si_0$ and the true precision matrix is $\om_0$, whose dimensions are $s\times s$. Each parameter $\theta$ we consider here corresponds to a positive definite matrix $\om$, since $\theta$ comes from the vectorization of the upper triangle of such a $\om$. The true parameter $\theta_0$ comes from $\om_0$. We can take $\V_n=n^{-1}\sum_{i=1}^n (g(D_i,\tilde \theta)-\bar g(\D,\tilde \theta))(g(D_i,\tilde \theta)-\bar g(\D,\tilde \theta))^\top$, where $\tilde \theta$ is the estimated parameter by inverting the empirical covariance matrix. Then we have the following theorem for partial correlation selection.

\noindent {\bf Theorem S3} {\it
For the moment function (7) in the main paper, suppose that Assumptions 1, 2, 3, 7 and 8 hold for $p=s(s+1)/2$. In addition, if\\
\noindent (1) Uniformly for all $\theta \in \Theta$, the corresponding $\om$ has eigenvalues bounded above and below by constants;\\
\noindent (2) For any random vector $Y=(Y_1,\ldots,Y_s)^\top$, $\sup_{1\leq j\leq s} \Ep(Y_j^8)< \infty$;\\
\noindent (3) The eigenvalues of $\V=\var(g(D,\theta_0)) $ are bounded from above and below by constant;\\
then Assumptions 4, 5 and 6 hold with $\alpha=1/2$. Therefore BGMM for the moment function (7) in the main paper satisfies the Bayesian oracle property in Theorem 1.
}

\begin{remark}
Here we restrict the space of the precision matrix $\om$ to a (possibly large) convex and compact set, with boundaries set by the smallest and the largest eigenvalues of $\om$ as in the condition (1). The boundedness of $\sup_j \Ep(Y_j^8)$ in the condition (2) is to guarantee the convergence of $\V_n$ to $\V$. Here we have directly assumed that the eigenvalues of $\V$ are bounded from above and below, mainly because this condition is not trivial and can hardly be obtained from any low level conditions.
\end{remark}

\noindent \textbf{Proof of Theorem S3:}\\
Hereafter we denote the $(i,j)$th entry in a generic $s\times s$ positive definite matrix $\si$ or $\om$ as $\sigma_{ij}$ or $\omega_{ij}$, respectively. Denote the $(i,j)$th entry in the true covariance matrix $\si_0$ and the true precision matrix $\om_0$ as $\sigma_{ij,0}$ or $\omega_{ij,0}$, respectively. For a generic parameter $\theta$, we denote the corresponding precision matrix as $\om$ and the corresponding covariance matrix as $\si=\om^{-1}$. The coordinates of $\theta$ and any other $p$-dimensional vector is subscripted by ``$ij$" with $1\leq i\leq j\leq s$. Then we can first establish an equivalence between the $L_2$ norm of $\theta$ and the Frobenius norm of $\om$. Since $\theta$ contains the entries in the upper triangle of $\om$, it is obvious that
\begin{equation}\label{l2f}
\frac{1}{2} \|\om- \om_0\|_F^2 \leq \|\theta-\theta_0\|^2 \leq \|\om- \om_0\|_F^2,
\end{equation}
so these two norms are equivalent.

Now we check Assumptions 4 and 5. For Assumption 4(i), since $m=p=s(s+1)/2$ for this example, we have that for any $\eta \in \mathbb{R}^p$ and $\|\eta\|^2=\sum_{1\leq i\leq j\leq s}\eta_{ij}^2=1$, by the Cauchy-Schwarz inequality,
\begin{align}\label{a41}
& \Ep\left\{\eta^\top (g(D,\theta)-g(D,\theta_0))\right\}^2 = \left\{\sum_{1\leq i\leq j\leq s}\eta_{ij}(\sigma_{ij}-\sigma_{ij,0}) \right\}^2
\nonumber \\
\leq {}& \sum_{1\leq i\leq j\leq s} \eta_{ij}^2 \cdot \sum_{1\leq i\leq j\leq s} (\sigma_{ij}-\sigma_{ij,0})^2  \leq \sum_{1\leq i\leq s, 1\leq j\leq s} (\sigma_{ij}-\sigma_{ij,0})^2  \nonumber \\
 = {}& \big\|\si-\si_0\big\|_F^2 = \big\|\om^{-1}(\om_0-\om)\om_0^{-1}\big\|_F^2 \leq  \big\|\om^{-1}\big\|_F^2 \big\|\om-\om_0\big\|_F^2 \big\|\om_0^{-1}\big\|_F^2 \nonumber \\
 \leq {}& \underline \lambda(\om)^{-2} \underline \lambda (\om_0)^{-2} s^2 \big\|\om-\om_0\big\|_F^2 = O\left((p^{1/2}\|\theta-\theta_0\|)^2\right),
\end{align}
where we used the submultiplicativity of the Fronbenius norm, the boundedness of eigenvalues in the condition (1), the relation $\|\A\|_F^2\leq s\bar \lambda(\A)^2$ for a $s\times s$ positive definite matrix $\A$, the relation $p=s(s+1)/2$ and \eqref{l2f}. Take supremum over $\eta$ and Assumption 4(i) is proved.

For Assumption 4(ii), we have derived above that the envelope function of $\mathcal{F}$ has $L_2$ norm of order $O(p^{1/2})$ given $\|\theta\|\leq R$. Note that in fact for the partial correlation selection example, the functions in $\mathcal{F}$ do not have any randomness. Suppose a $L_2$ $\epsilon$-net of $\mathcal{F}$ is $\{(\eta_1,\theta_1),\ldots,(\eta_N,\theta_N)\}$ for $N=N\left(\epsilon\|F\|_{P_D,2},\mathcal{F},L_2(P_D)\right)$. Then for any $f(\eta,\theta)\in \mathcal{F}$, we apply a similar procedure of \eqref{a41} and have
\begin{align*}
& \Ep\left|f(\eta_k,\theta_k)-f(\eta,\theta)\right|^2\leq 2|f(\eta_k,\theta_k)-f(\eta,\theta_k)|^2+2|f(\eta,\theta_k)-f(\eta,\theta)|^2 \\
\leq {} & 2\|\eta_k-\eta\|^2 \underline \lambda(\om_k)^{-2} \underline \lambda (\om_0)^{-2} s^2 \cdot 4R^2 + 2\underline \lambda(\om_k)^{-2} \underline \lambda (\om)^{-2}s^2 \cdot 2\|\theta_k-\theta\|^2\\
:= {} & (C_1p^{1/2}\|\eta_k-\eta\|)^2 + (C_2p^{1/2} \|\theta_k -\theta\|)^2,
\end{align*}
where $\om_k$ is the matrix $\om$ with parameter $\theta_k$. Therefore by a similar argument to the proof of Theorem S1,
$N\left(\epsilon\|F\|_{P_D,2},\mathcal{F},L_2(P_D)\right)\leq \left(9C_1C_2p^2R/\epsilon^2\right)^p,$
which give $\ln N\left(\epsilon\|F\|_{P_D,2},\mathcal{F},L_2(P_D)\right) = O\left(p\ln(n/\epsilon)\right)$ and hence Assumption 4(ii) holds.

For Assumption 5(i), we have
\begin{align}\label{a51}
& \|\Ep g(D,\theta)\|^2 = \sum_{1\leq i\leq j\leq s} (\sigma_{ij}-\sigma_{ij,0})^2 \geq \frac{1}{2} \sum_{1\leq i\leq s, 1\leq j\leq s} (\sigma_{ij}-\sigma_{ij,0})^2 \nonumber \\
={} & \frac{1}{2} \big\| \si-\si_0\big\|_F^2 = \frac{1}{2} \big\|\om^{-1}(\om_0-\om)\om_0^{-1} \big\|_F^2 \geq
\frac{1}{2} \bar \lambda(\om)^{-2} \bar \lambda(\om_0)^{-2}  \big\|\om-\om_0\big\|_F^2,
\end{align}
where we have used the fact that for two positive definite matrices $\A$ and $\B$,
$$\|\A\B\|_F^2=\tr(\B^\top\A^\top\A\B) \geq \underline \lambda (\B)^2\tr(\A^\top\A)\geq \underline \lambda(\B)^2 \|\A\|_F^2.$$
Now since we have assumed in the condition (1) that the eigenvalues of $\om$ are bounded above by constants, \eqref{a51} implies that
$\|\Ep g(D,\theta)\|\geq \delta_1 \|\theta-\theta_0\|$ with $0<\delta_1<\bar \lambda(\om)^{-1} \bar \lambda(\om_0)^{-1}/\sqrt{2}$. So Assumption 5(i) holds with this $\delta_1$ and $\delta_0=R\delta_1$.

To show Assumption 5(ii), we only need to show that for any unit vector $u \in \mathbb{R}^p$, $u^\top \G^\top\G u$ is bounded above and below by constants, where $\G=\nabla_{\theta} \Ep g(D,\theta_0)$. Define a linear operator $\partial_u := \sum_{1\leq i\leq j\leq s} u_{ij}\frac{\partial}{\partial \theta_{ij}}$ and define $u_{ji}:=u_{ij}$ for any $j<i$. Then $u^\top \G^\top\G u = \partial_u \Ep g(D,\theta_0)^\top \partial_u \Ep g(D,\theta_0)=\|\partial_u \Ep g(D,\theta_0)\|^2$, and similar to \eqref{l2f}, one can show that
$$\frac{1}{2} \|\partial_u \om_0^{-1}\|_F^2 \leq \|\partial_u Eg(D,\theta_0)\|^2 \leq \|\partial_u \om_0^{-1}\|_F^2,$$
and also
$$1 =\|u\|^2 \leq \|\partial_u \om_0\|_F^2 \leq 2\|u\|^2=2.$$
Since $\om \om^{-1}=\I$, we take first derivative and have $\partial_u \om^{-1} = -\om^{-1}(\partial_u \om) \om^{-1}$. Therefore, we have
\begin{align*}
&  \|\partial_u \om_0^{-1} \|_F^2 = \|\om_0^{-1}(\partial_u \om_0) \om_0^{-1}\|_F^2 \leq \underline \lambda(\om_0)^{-4} \|\partial_u \om_0\|_F^2 \leq 2\underline \lambda(\om_0)^{-4}\\
&  \|\partial_u \om_0^{-1} \|_F^2 = \|\om_0^{-1}(\partial_u \om_0) \om_0^{-1}\|_F^2 \geq \bar \lambda(\om_0)^{-4} \|\partial_u \om_0\|_F^2 \geq \bar \lambda(\om_0)^{-4},
\end{align*}
which implies the boundedness of eigenvalues of $\G^\top\G$, given the condition (1).

For Assumption 5(iii), we use the same technique and have that for unit vectors $u,v\in \mathbb{R}^p$, $\|\H(\theta)(u,v)\|^2=\|\partial_u\partial_v Eg(D,\theta)\|^2 \leq \|\partial_u\partial_v \om^{-1}\|_F^2$. While for any generic $\om$, using $\partial_u \om^{-1} = -\om^{-1}(\partial_u \om) \om^{-1}$, we get
$$\partial_u \partial_v \om^{-1} = -\om^{-1} (\partial_u \om) \om^{-1}(\partial_v \om) \om^{-1}-\om^{-1} (\partial_v \om) \om^{-1}(\partial_u \om) \om^{-1} - \om^{-1}(\partial_u \partial_v \om) \om^{-1}. $$
Therefore since $\partial_u \partial_v \om \equiv 0$ for any $\om$,
\begin{align*}
&\|\H(\theta)(u,v)\|^2 \leq  \| \om^{-1} (\partial_u \om) \om^{-1}(\partial_v \om) \om^{-1}\|_F^2 + \| \om^{-1} (\partial_v \om) \om^{-1}(\partial_u \om) \om^{-1}\|_F^2\\
& \leq 2\underline \lambda (\om)^{-6} + 2\underline \lambda (\om)^{-6} = 4\underline \lambda (\om)^{-6},
\end{align*}
which is bounded above by constant for any $\om$ considered here by the condition (1). Hence Assumption 5(iii) holds.

For Assumption 6, we have assumed the boundedness of eigenvalues of $\V$, and we still need to show the convergence $\|\V_n-\V\|\to 0$ w.p.a.1 as $n\to \infty$. By Chebyshev's inequality, for any $C>0$, we have
\begin{align*}
& \Pr\big(\|\V_n-\V\|\geq C\big) \leq \Pr \big(\|\V_n-\V\|_F \geq C\big)    \\
\leq {}& \frac{1}{C^2} \sum_{1\leq j_1\leq k_1\leq s} \sum_{1\leq j_2\leq k_2\leq s}
\var\Big\{\frac{1}{n}\sum_{i=1}^n(Y_{ij_1}Y_{ik_1}-\sigma_{j_1k_1,0})(Y_{ij_2}Y_{ik_2}-\sigma_{j_2k_2,0})\\
&-\Ep \Big\{(Y_{j_1}Y_{k_1}-\sigma_{j_1k_1,0})(Y_{j_2}Y_{k_2}-\sigma_{j_2k_2,0})\Big\} \Big\}\\
\leq {}& \frac{s^2(s+1)^2}{4nC^2}\sup_{j_1,k_1,j_2,k_2} \Ep\left\{(Y_{j_1}Y_{k_1}-\sigma_{j_1k_1,0})^2(Y_{j_2}Y_{k_2}-\sigma_{j_2k_2,0})^2\right\}\\
\leq {}& \frac{4s^2(s+1)^2}{nC^2}\sup_{j_1,k_1,j_2,k_2}\Ep\left\{Y_{j_1}^2Y_{k_1}^2Y_{j_2}^2Y_{k_2}^2\right\}\leq \frac{16p^2}{nC^2}\sup_{1\leq j\leq s}\Ep(Y_j^8).
\end{align*}
Hence we have $\|\V_n-\V\|=O_p\left(p/\sqrt{n}\right)=o_p(1)$. This together with Condition (3) implies the boundedness of eigenvalues of $\V_n$ w.p.a.1 as $n\to \infty$. Therefore Assumption 6 holds.\hfill$\blacksquare$


\vspace{+2.5cm} \noindent {\bf
\Large 4. Examples of the prior on models}
\vspace{+.5cm}

In this section, we verify Assumption 8 for several examples of priors on models. They are summarized in the following proposition.

\noindent {\bf Proposition} {\it
Assumption 8 is satisfied by the following three priors on models:\\
\noindent (a) Every component of $\theta$ enters the model $\M$ independently with probability $\nu\in (0,1)$, excluding the empty model: $\pi(\M) \propto \nu^{|\M|}(1-\nu)^{p-|\M|}$. Here $\nu$ is either a fixed constant or $\nu=n^{-c}$ for some constant $c>0$. \\
\noindent (b) The prior factorizes as $\pi(\M)=\pi(\M||\M|=k)\pi(|\M|=k)$ for $k=1,\ldots,p$. $\pi(\M~|~|\M|=k) = \binom{p}{k}^{-1}$, and $\pi(|\M|=k) \propto \zeta^k e^{-\zeta}/k!$ for some constant $\zeta>0$. \\
\noindent (c) The prior factorizes as $\pi(\M)=\pi(\M||\M|=k)\pi(|\M|=k)$ for $k=1,\ldots,p$. $\pi(\M~|~|\M|=k) = \binom{p}{k}^{-1}$, and $\pi(|\M|=k) \propto e^{-\zeta k }$ for some constant $\zeta>0$.
}
\vspace{.8cm}

\noindent \textbf{Proof of the proposition:}\\
We verify Assumption 8 (i) and (ii) for each of the three priors. We note that if $p$ satisfies the growth rate in Assumption 2, $p\prec \sqrt{n}/p\to \infty$ as $n\to\infty$.\\

For the prior in (a), let $C_1 = \max\left(\frac{\nu}{1-\nu},\frac{1-\nu}{\nu} \right)$. Then $1\leq C_1 \leq n^c$. We have
$$\frac{\pi(\M)}{\pi(\M_0)}\leq \frac{\nu^{|\M|}(1-\nu)^{p-|\M|}}{\nu^{|\M_0|}(1-\nu)^{p-|\M_0|}} = \left(\frac{\nu}{1-\nu}\right)^{|\M|-|\M_0|}.$$
If $\M\supset \M_0$, $|\M|-|\M_0|\geq 1$, then
$$\frac{\pi(\M)}{\pi(\M_0)}\leq \left(\frac{\nu}{1-\nu}\right)^{|\M|-|\M_0|} \prec \left(\sqrt{n}/p\right)^{|\M|-|\M_0|}.$$
If $\M_0\backslash \M \neq \emptyset$, then since $\left||\M|-|\M_0|\right|\leq p$,
$$\frac{\pi(\M)}{\pi(\M_0)}\leq C_1^{\left||\M|-|\M_0|\right|}\leq n^{cp} \preceq e^{r_1p\ln n},$$
where we can choose $r_1=c$. So Assumption 8 holds for the prior in (a).

For the prior in (b), let $|\M|=k$ and $|\M_0|=k_0$. We have
$$\frac{\pi(\M)}{\pi(\M_0)}=\frac{\pi(\M||\M|=k)\pi(|\M|=k)}{\pi(\M_0||\M_0|=k_0)\pi(|\M|=k_0)}
=\frac{\binom{p}{k}^{-1}\frac{\zeta^k e^{-\zeta}}{k!}}{\binom{p}{k_0}^{-1}\frac{\zeta^{k_0} e^{-\zeta}}{{k_0}!}}
= \frac{(p-k)!}{(p-k_0)!}\zeta^{k-k_0}.$$
If $\M\supset \M_0$, then $k-k_0\geq 1$ and $(p-k)!<(p-k_0)!$. It follows that
$$\frac{\pi(\M)}{\pi(\M_0)}\leq \zeta^{k-k_0} \prec (\sqrt{n}/p)^{k-k_0}.$$
If $\M_0\backslash \M \neq \emptyset$, then since $\left|k-k_0\right|\leq p$, $(p-k)!\leq p^p$, and $\ln p\preceq \ln n$ by Assumption 2,
$$\frac{\pi(\M)}{\pi(\M_0)}\leq p^p \zeta^p = e^{p\ln p + p\ln \zeta} \preceq e^{r_1 p\ln n},$$
for some constant $r_1>0$. So Assumption 8 holds for the prior in (b).

For the prior in (c), let $|\M|=k$ and $|\M_0|=k_0$. We have
$$\frac{\pi(\M)}{\pi(\M_0)}=\frac{\pi(\M||\M|=k)\pi(|\M|=k)}{\pi(\M_0||\M_0|=k_0)\pi(|\M|=k_0)}
=\frac{\binom{p}{k}^{-1} e^{-\zeta k}}{\binom{p}{k_0}^{-1} e^{-\zeta k_0}}
= \frac{k!(p-k)!}{k_0!(p-k_0)!}e^{-\zeta(k-k_0)}.$$
If $\M\supset \M_0$, then $k-k_0\geq 1$, $(p-k)!<(p-k_0)!$, $k!/k_0!\leq k^{k-k_0}\leq p^{k-k_0}$. It follows that
$$\frac{\pi(\M)}{\pi(\M_0)}\leq p^{k-k_0} e^{-\zeta(k-k_0)} \prec \left(\sqrt{n}/p\right)^{k-k_0}.$$
For $\M_0\backslash \M \neq \emptyset$, notice that we always have $(k!(p-k)!)/(k_0!(p-k_0)!) \leq p^{|k-k_0|}$. Therefore
$$\frac{\pi(\M)}{\pi(\M_0)}\leq p^{|k-k_0|} e^{-\zeta (k-k_0)}\leq p^p e^{\zeta p} \leq e^{p\ln p + p \zeta} \preceq e^{r_1 p\ln n},$$
for some constant $r_1>0$. So Assumption 8 holds for the prior in (c). \hfill $\blacksquare$

\bibliographystyle{plainnat}
\bibliography{bgmmbib}

\end{document}